\documentclass[12pt,notitlepage]{article}
\usepackage{amsmath,amsfonts}
\usepackage{multicol}
\usepackage{xifthen}
\usepackage{latexsym}
\usepackage[colorlinks=true, breaklinks=true, urlcolor= black, linkcolor= black, citecolor= black, bookmarksopen=true,linktocpage=true]{hyperref}

\numberwithin{equation}{section}

\newcounter{hasqed}
\newcommand{\newproof}[2]{
 \newenvironment{#1}{\setcounter{hasqed}{1}\begin{trivlist}\item{{\bf #2}\ }}%
 {\ifthenelse{\value{hasqed}=1}{\hfill$\square$}{}\setcounter{hasqed}{1}\end{trivlist}}}
 \newcommand{\qedhere}{\hfill$\square$\setcounter{hasqed}{0}}

\newcommand{\ovl}[1]{{\overline{#1}}}

\newcommand{\ra}{\rightarrow}

\newcommand{\iso}{\cong}

\newcommand{\card}[1]{{\left|#1\right|}}
\newcommand{\im}[1]{\null{{\rm Im}#1}}
\newcommand{\ggp}{{G}}
\newcommand{\hgp}{{H}}
\newcommand{\lgp}{{L}}
\newcommand{\fingp}{{F}}
\newcommand{\setdiff}{{-}}
\newcommand{\GL}[1]{\null{{\rm GL}_{#1}}}
\newcommand{\Ind}[3]{{{\rm Ind}^{#1}_{#2}{#3}}}
\newcommand{\Res}[3]{{{\rm Res}^{#1}_{#2}{#3}}}
\newcommand{\nilgp}{{\Gamma}}
\newcommand{\kgp}{{K}}

\newcommand{\ip}[3]{{\left\langle {#1},{#2} \right\rangle_{#3}}}

\newcommand{\indx}[2]{{\left|{#1} \colon {#2}\right|}}

\newcommand{\val}{\mathop{\rm val}}

\newcommand{\gp}{{G}}
\newcommand{\irr}[2]{{{\rm R}_{#1}({#2})}}
\newcommand{\cirr}[2]{{{\rm R}^{\rm ad}_{#1}({#2})}}
\newcommand{\pirr}[2]{{{\rm R}^{(p)}_{#1}({#2})}}

\newcommand{\twzeta}[3]{{\zeta_{{#1},{#2}}({#3})}}
\newcommand{\rts}[1]{{\mu_{#1}}}
\newcommand{\pair}{{\Psi}}
\newcommand{\globlin}{{\Xi}}
\newcommand{\pgp}{{\widehat{\nilgp}_p}}
\newcommand{\tfnilgp}{{\Delta}}
\newcommand{\tfpgp}{{\widehat{\tfnilgp}_p}}
\newcommand{\kr}{{\Theta}}
\newcommand{\pkr}{{\widehat{\kr}_p}}
\newcommand{\jinc}{{j}}
\newcommand{\proj}{{\pi}}
\newcommand{\pproj}{{\widehat{\proj}_p}}

\makeatletter
\def\blfootnote{\xdef\@thefnmark{}\@footnotetext}
\makeatother

\newcommand{\Latt}[1][]{{{ S}_{#1}}}
\newcommand{\Tor}[1][]{{{ T}_{#1}}}
\newcommand{\K}{{K}}
\newcommand{\res}{{k}}
\newcommand{\Valgp}{{\Gamma}}

\newtheorem{thm}{Theorem}[section]
\newenvironment{theorem}[1][]{\begin{thm}\ifthenelse{\isempty{#1}}{}{\label{thm:#1}}}{\end{thm}}
\newtheorem{prop}[thm]{Proposition}
\newenvironment{proposition}[1][]{\begin{prop}\ifthenelse{\isempty{#1}}{}{\label{prop:#1}}}{\end{prop}}
\newtheorem{lem}[thm]{Lemma}
\newenvironment{lemma}[1][]{\begin{lem}\ifthenelse{\isempty{#1}}{}{\label{lem:#1}}}{\end{lem}}
\newtheorem{cor}[thm]{Corollary}
\newenvironment{corollary}[1][]{\begin{cor}\ifthenelse{\isempty{#1}}{}{\label{cor:#1}}}{\end{cor}}
\newtheorem{conj}[thm]{Conjecture}

\newtheorem{quest}[thm]{Question}

\newtheorem{ITdefn}[thm]{Definition}
\newenvironment{definition}[1][]
 {\begin{ITdefn}\ifthenelse{\isempty{#1}}{}{\label{def:#1}}\rm}{\end{ITdefn}}
\newenvironment{defn}[1][]
 {\begin{ITdefn}\ifthenelse{\isempty{#1}}{}{\label{def:#1}}\rm}{\end{ITdefn}}
\newtheorem{ITclm}[thm]{Claim}
\newenvironment{claim}[1][]
 {\begin{ITclm}\ifthenelse{\isempty{#1}}{}{\label{claim:#1}}\rm}{\end{ITclm}}
\newtheorem{ITex}[thm]{Example}
\newenvironment{example}[1][]
 {\begin{ITex}\ifthenelse{\isempty{#1}}{}{\label{ex:#1}}\rm}{\end{ITex}}
 \newtheorem{ITrem}[thm]{Remark}
\newenvironment{remark}[1][]
 {\begin{ITrem}\ifthenelse{\isempty{#1}}{}{\label{rem:#1}}\rm}{\end{ITrem}}
\newtheorem{ITnot}[thm]{Notation}
\newenvironment{notation}[1][]
 {\begin{ITnot}\ifthenelse{\isempty{#1}}{}{\label{not:#1}}\rm}{\end{ITnot}}
\newproof{proof}{Proof}
\newproof{warning}{Warning}
\newtheorem{ITnwarn}[thm]{Warning}

\newtheorem{ITobs}[thm]{Observation}
\newenvironment{observation}
 {\begin{ITobs}\rm}{\end{ITobs}}
\newtheorem{ITconv}[thm]{Convention}

\newtheorem{ITpor}[thm]{Porism}

\newtheorem{ITguess}[thm]{Guess}

\newtheorem{ITassumption}[thm]{Assumption}

\newfont{\sBbb}{msbm7 scaled\magstep1} 

\def\Qq {\mathbb{Q}}
\def\Zz {\mathbb{Z}}
\def\Ff {\mathbb{F}}
\def\Rr {\mathbb{R}}
\def\Cc {\mathbb{C}}
\def\Nn {\mathbb{N}}
\def\Bb {\mathbb{B}}
\def\Gg {\mathbb{G}}

\def\unif{\lambda_0}
 \def\,{\ }

\def\code#1{ {< #1 >} }
\def\:{{\setminus}}

\def\acl{{\rm acl}} \def\dcl{{\rm dcl}}
\def\mod{{\rm mod}}
\def\dim{{\rm dim}}

\def\dom{{\rm dom}}
\def\rad{{\rm rad}}

\def\implies{ \Rightarrow }

\def\tensor{{\otimes}}

\def\tT{\widetilde{T}}   \def\tL{\widetilde{\LL}} \def\tP{\widetilde{p}} \def\tM{\widetilde{M}}
\def\ts{\widetilde{\sigma}}\def\ns{\sigma}
\def\nT{T}   \def\nL{\LL} \def\nM{M} \def\nP{p}

\def\RR{{\cal R}}

\newcommand\bA{\bar{A}}

\newcommand\LL{\mathcal{L}}
\newcommand\LG{\mathcal{L}_{\mathcal{G}}}
\newcommand\LGminus{\mathcal{L}_{\mathcal{G}}^{-}}
\newcommand\LGm{\mathcal{L}_{\mathcal{G}}^{\Nn}}
\newcommand\LGn{\mathcal{L}_{\mathcal{G}}^{\mathcal{N}}}
\newcommand\eq[1]{#1^{\rm eq}}
\newcommand\dcleq{\dcl^{\rm eq}}
\newcommand\DCLeq{\DCL^{\rm eq}}
\newcommand\acleq{\acl^{\rm eq}}
\newcommand\ACLeq{\ACL^{\rm eq}}
\newcommand\DCL{{\dcl}_{\nL}}
\newcommand\ACL{{\acl}_{\nL}}
\newcommand\tpL{{\rm tp}_{\nL}}
\newcommand\tpl{{\rm tp}_{\tL}}
\newcommand\tp{\mathop{\rm tp}}

\newcommand\dcll{{\dcl}_{\tL}}
\newcommand\acll{{\acl}_{\tL}}
\newcommand\Th{{\rm Th}}
\newcommand\THL{{\rm Th}_{\nL}}
\newcommand\Thl{{\rm Th}_{\tL}}
\newcommand\RV[1][]{{\rm RV}_{\!#1}}
\newcommand\rv{\mathop{\rm rv}}
\newcommand\resf{\mathop{\rm res}}
\newcommand\latt{\mathop{\Lambda}}
\newcommand\oBallSet{\dot{\Bb}}
\newcommand\cBallSet{\overline{{\Bb}}}
\newcommand\cpow[2]{#1^{#2}}
\newcommand\gpow[2]{#1^{\cdot#2}}
\newcommand\inv[1]{#1^{*}}
\newcommand\fix{\mathop{\rm fix}}

\DeclareMathOperator{\algop}{alg}
\newcommand{\alg}[1]{\overline{#1}^{\algop}}
\newcommand\aut{\mathop{\rm Aut}}
\newcommand\ACVF{\mathop{\rm ACVF}}
\newcommand\Lrg{\LL_{\mathrm{rg}}}
\newcommand{\PLO}{\mathop{\rm PL}_{0}}
\newcommand{\germ}[2]{\partial_{#1}#2}
\newcommand{\Val}{\mathcal{O}}
\newcommand{\triang}{{\rm B}}
\newcommand{\BallSet}{\Bb}
\newcommand{\HFO}{\mathop{\rm HF}_{0}}
\newcommand{\sminus}{\setdiff}
\newcommand{\restr}[2]{{\left.#1\right|_{#2}}}
\newcommand{\Langrestr}{\restr}
\newcommand{\tprestr}[2]{#1|#2}
\newcommand{\substr}{\leq}

\newcommand{\ang}{\mathop{\rm ac}}

\newcommand{\ord}{\operatorname{ord}}
\def\ac{{\overline{\rm ac}}}
\def\cL{{\mathcal L}}
\newcommand{\Vol}{\operatorname{Vol}}
\newcommand{\MVol}{\operatorname{MultiVol}}
\newcommand{\MB}{\operatorname{MB}}
\newcommand{\MNum}{\operatorname{MultiNumber}}

\renewcommand{\mid}{:}

\newenvironment{thmenum}{\begin{enumerate}}{\end{enumerate}}

\begin{document}

\title{Definable equivalence relations and zeta functions of groups}
\author{Ehud Hrushovski, Ben Martin and Silvain Rideau,\\ with an appendix by Raf Cluckers}
\date{September 1, 2017}

\maketitle

\centerline{\emph{In memory of Fritz Grunewald.}}\bigskip

\begin{abstract}
We prove that the theory of the $p$-adics $\Qq_p$ admits elimination of imaginaries provided we add  a sort for $\GL{n}(\Qq_p)/\GL{n}(\Zz_p)$ for each $n$. We also prove that the elimination of imaginaries is uniform in $p$. Using $p$-adic and motivic integration, we deduce the uniform rationality of certain formal zeta functions arising from definable equivalence relations. This also yields analogous results for definable equivalence relations over local fields of positive characteristic.  The appendix contains an alternative proof, using cell decomposition, of the rationality (for fixed $p$) of these formal zeta functions that extends to the subanalytic context.

As an application, we prove rationality and uniformity results for zeta functions obtained by counting twist isomorphism classes of irreducible representations of finitely generated nilpotent groups; these are analogous to similar results of Grunewald, Segal and Smith and of du Sautoy and Grunewald for subgroup zeta functions of finitely generated nilpotent groups.
\end{abstract}

\section{Introduction}
\label{sec:intro}

\blfootnote{2010 {\em Mathematics Subject Classification.}  03C60 (03C10 11M41 20E07 20C15).}
\blfootnote{{\em Key words and phrases.}  Elimination of imaginaries, invariant extensions of types, cell decompositions, rational zeta functions, subgroup zeta functions, representation zeta functions}

This paper concerns the model theory of the $p$-adic numbers $\Qq_p$ and applications to certain counting problems arising in group theory.  Recall that a theory (in the model-theoretic sense of the word) is said to have {\em elimination of imaginaries} (EI) if the following holds: for every model $M$ of the theory, for every \(\emptyset\)-definable subset $D$ of some $\cpow{M}{n}$ and for every \(\emptyset\)-definable equivalence relation $R$ on $D$, there exists an \(\emptyset\)-definable function $f\colon D\to \cpow{M}{m}$, for some $m$, such that the fibers of $f$ over $f(D)$ are precisely the equivalence classes of $R$. In other words, elimination of imaginaries states that every pair $(D',E')$ (consisting of an \(\emptyset\)-definable set $D'$ and an \(\emptyset\)-definable equivalence relation $E'$ on it) reduces to a pair $(D,E)$ where $E$ is equality---here, as in descriptive set theory, we say that $(D',E')$ {\em reduces to} $(D,E)$ if there exists a $\emptyset$-definable map $f\colon D' \to D$ with $xE'y \iff f(x) E f(y)$.

The theory of $\Qq_p$ (in the language of rings with a predicate for $\val(x)\geq \val(y)$) does not admit EI \cite{ScoMac}: for example, no such $f$ exists for the definable equivalence relation $R$ on $\Qq_p$ given by $x R y$ if $\val(x-y)\geq 1$, because $\Qq_p/R$ is countably infinite but any definable subset of $\cpow{\Qq_p}{m}$ is either finite or uncountable. Our first main theorem gives a $p$-adic EI result when we add for each $n$ a sort $S_n$ for the family of $\Zz_{p}$-lattices in $\Qq_p^n$. These new sorts are called the {\em geometric imaginaries}. The language $\LGminus$ consists of the valued field sort and the sorts $S_n$ (with some more structure described in Section~\ref{sec:elim}).

\begin{theorem}
\label{thm:Qp_intro}
The theory of $\Qq_p$ eliminates imaginaries in the language $\LGminus$.
\end{theorem}

\noindent To be precise, we prove a version of this that holds for any finite extension of $\Qq_p$ (Theorem~\ref{thm:Qp}).

Suppose we are given not just a single definable equivalence relation for some fixed $\Qq_p$, but one for every $\Qq_p$.  For our applications to zeta functions below, we want to control the behavior of the elimination of imaginaries as we vary the prime $p$.  Our second main result is that the theory of ultraproducts of $\Qq_{p}$ also eliminates imaginaries if we add similar sorts.

\begin{theorem}
\label{thm:unifQp_intro}
The theory of non-principal ultraproducts $\prod_{p}\Qq_p/\mathcal{U}$ eliminates imaginaries in the language $\LGminus$ provided we add some constants.
\end{theorem}

\noindent See Theorem~\ref{thm:ultraprod Qp} for a more precise statement of what constants are needed to eliminate imaginaries. This last result implies that the elimination of imaginaries in $\Qq_{p}$ is uniform in $p$; see Corollary~\ref{cor:EIunif} for a precise statement of this uniformity.

In fact, we prove a more general result (Corollary~\ref{cor:EIcrit}), which yields EI both for $\Qq_p$ and for ultraproducts: given two theories $\nT$, $\tT$ satisfying certain hypotheses, $\nT$ has EI if $\tT$ does.  In our application, $\tT$ is the theory of algebraically closed valued fields of mixed characteristic ($\ACVF_{0,p}$) or equicharacteristic zero ($\ACVF_{0,0}$) and $\nT$ is either the theory of a finite extension of $\Qq_p$ or the theory of an ultraproduct of $\Qq_p$ where $p$ varies, with appropriate extra constants in each case (in fact in the latter case Corollary~\ref{cor:EIcrit} does not apply immediately but a variant does).

The notion of an invariant extension of a type plays a key part in our proof.  If $T$ is a theory, $M\models T$, $A\subseteq M$ and $p$ is a type over $A$ then an {\em invariant extension} of $p$ is a type $q$ over $M$ such that $\tprestr{q}{A}= p$ and $q$ is $\aut(M/A)$-invariant.  The theory $\ACVF$ is not stable; in \cite{HHM,HHM2}, Haskell, the first author and Macpherson used invariant extensions of types to study the stability properties of $\ACVF$ and to define notions of forking and independence. They proved that $\ACVF$ plus some extra sorts admits EI.

As an important consequence of Theorems~\ref{thm:Qp} and \ref{thm:ultraprod Qp}, we prove
the following rationality and uniformity result for zeta functions $S_p(t)$ counting the number of equivalence classes in some uniformly definable family of equivalence relations.  (Here $t= (t_1,\ldots, t_r)$ is a tuple of indeterminates and $S_p(t)$ is a power series in the $t_i$; we obtain a zeta function in the more usual sense by setting $t_i= p^{-s_i}$, where the $s_i$ are complex variables.)

\begin{theorem}
\label{thm:rat_intro}
 The zeta functions $(S_p(t))_{p \textrm{ prime}}$ are uniformly rational.
\end{theorem}

\noindent We also give a version of Theorem~\ref{thm:rat_intro} for a uniformly definable family of equivalence relations over a local field of positive characteristic (Corollary~\ref{cor:transfer}).

See Section~\ref{sec:rat} for definitions and a precise statement (Theorem~\ref{thm:rat}).  Roughly speaking, uniform rationality means that each $S_p(t)$ can be expressed as a rational function with coefficients in $\Qq$, where the denominator is a product of functions of the form $(1- p^at^b)$ or $p^n$ with $a,b,n$ independent of $p$, and the numerator is a polynomial in $t$ such that each coefficient comes from counting the $\Ff_p$-points of a fixed variety over $\Zz$.  In particular, we obtain that $S_p(t)$ is rational not just for all sufficiently large primes, but for every prime; this is crucial for our applications to representation growth below, as well as to the following result, which deals with the abscissa of convergence.

\begin{thm}
\label{thm:abscissa_rat_intro}
 Let $S_p(t)$ be as above and suppose we are in the one-variable case ($r= 1$, $t= t_1$).  Define $\zeta_p(s)= S_p(p^{-s})$.  Assume that the constant term of $\zeta_p(s)$ is 1 for all but finitely many primes and set $\zeta(s)= \prod_p \zeta_p(s)$.  Then the abscissa of convergence of $\zeta(s)$ is rational (or $-\infty$).
\end{thm}

In fact, Theorem~\ref{thm:rat} yields a kind of ``double uniformity'': the ultraproduct formalism allows us to vary not just the prime $p$, but also the choice of an extension $L_p$ of $\Qq_p$.  For an application of this double uniformity, see the end of Section~\ref{sec:twst}.

To describe the proof of Theorem~\ref{thm:rat_intro}, let us now come back to the meaning of our elimination of imaginaries result. It shows that any $(D',E')$  can be reduced to a $(D,E)$ of a special kind---namely, the equivalence relation on $\GL{N}(\Qq_p)$ for some $N$ whose equivalence classes are the left $\GL{N}(\Zz_p)$-cosets. The quotients $D/E$ have a specific geometric meaning---but can one explain abstractly in what way they are special? One useful observation is that we have reduced an arbitrary equivalence relation to a quotient by a definable group action. Another concerns volumes: the $E$-classes have volumes that are {\em motivically invertible} (in fact, each class is equivalent to a polydisk of an appropriate dimension and size).

Indeed, it is only the latter property of the geometric imaginaries that is actually used in the proof of Theorem~\ref{thm:rat}. This proof relies on representing the number of classes of some definable equivalence relation $E$ on some definable set $D$ as an integral. The  idea, going back to Denef and Igusa, is simple:  the number of classes of $E$ on $D$ equals the volume of $D$, for any measure such that each $E$-class has measure one. The question is how to come up (definably) with such a measure. The setting is that we already have the Haar measure $\mu$ on $\Qq_{p}$ (normalized so that $\mu(\Zz_{p}) = 1$), and for simplicity---one can easily reduce to this case---let us assume each $E$-class $[x]_E\in D/E$ has finite, nonzero measure. The problem then is to show that there exists a definable function $f\colon D\to \Qq_p$ such that the measure of each $E$-class $[x]_E$ is of the form
\begin{equation}
\label{eqn:eqvcclassmeas}
\mu([x]_E)= |f(x)|,
\end{equation}
where \(|f(x)|\) denotes the \(p\)-adic norm. Then we can replace $\mu$ with $|f|^{-1} \mu$. In practice, $f$ is usually given explicitly (cf.\ \cite[Section\,2]{GSS}). For more complicated equivalence relations, however, such as the one for representation zeta functions in Section~\ref{sec:twst}, it is not clear {\em a priori} that such an $f$ can be found, even in principle.

This point is beautifully brought out in work by Raf Cluckers; we are very pleased to have his permission to include it here as an Appendix. The Appendix contains a complete proof of the rationality results in Section~\ref{sec:rat} for fixed $p$ which also extends to the analytic case, while avoiding an explicit elimination of imaginaries. It might be useful to say a word here about the  two proofs. Given EI to the geometric sorts, we can represent $E$ as the coset equivalence relation of $\GL{n}(\Zz_p)$. In this case we can take the measure in our $p$-adic integral to be the Haar measure on $\GL{n}(\Qq_p)$, where each class automatically has measure one. The density of this measure with respect to (the \(n^2\)-fold Cartesian power of) the additive Haar measure is given by $M\mapsto 1/ | \det(M)|$. In other words, for this canonical $E$, the reciprocal of the (additive) Haar measure of any $E$-class {\em is} represented by a definable function.

In the Appendix, any equivalence relation $(D',E')$ is reduced to one with motivically invertible volumes. Indeed  there exists a $\emptyset$-definable $D\subseteq D'$ such that $D \cap e$ has motivically invertible volume for each $E'$-class $e$ (in particular, with $E = \restr{E'}{D}$, the natural map $D/E \to D'/E'$ is a bijection). This result is valid in the analytic case too, unlike geometric EI in its present formulation (see \cite{HHMAn}).

We illustrate the power of Theorem~\ref{thm:rat} by using it to prove rationality results for certain zeta functions of finitely generated nilpotent groups (Theorems~\ref{thm:ratsbgp} and \ref{thm:rattwst_intro}).  Grunewald, Segal and Smith \cite[Sec.~2]{GSS} showed that subgroups of $p$-power index of such a group $\Gamma$ can be parametrized $p$-adically if $\Gamma$ is also torsion-free.  More precisely, these subgroups can be {\em interpreted}: that is, placed in bijective correspondence with the set of equivalence classes of some definable equivalence relation on a definable subset $D$ of some $\cpow{\Qq_p}{N}$.  Let $b_n<\infty$ denote the number of subgroups of $\Gamma$ of index $n$.  Using $p$-adic integration over $D$ and results of Denef and Macintyre, Grunewald, Segal and Smith showed that the $p$-local subgroup zeta function $\sum_{n=0}^\infty b_{p^n}t^n$ is a rational function of $t$, and that the degrees of the numerator and denominator of this rational function are bounded independently of $p$.  Du Sautoy and others have calculated subgroup zeta functions explicitly in many cases \cite{duSSubGp,duSWoo,VolBT} and studied uniformity questions.
For instance, du Sautoy and Grunewald proved a uniformity result by showing that the $p$-adic integrals that arise in the calculation of subgroup zeta functions fall into a special class they call {\em cone integrals} \cite{duSGruAnn}. See the start of Section~\ref{sec:zetagp} for further discussion of uniformity in the context of subgroup zeta functions.

We also consider situations where it is not clear how to find a definable function $f$ satisfying Eqn.\ (\ref{eqn:eqvcclassmeas}) and construct suitable definable $p$-adic integrals.  The main one, and the original motivation for our results, is in the area of representation growth.  This is analogous to subgroup growth: one counts not the number $b_{p^n}$ of index $p^n$ subgroups of a group $\Gamma$, but the number $a_{p^n}$ of irreducible $p^n$-dimensional complex characters of $\Gamma$ (modulo tensoring by one-dimensional characters if $\Gamma$ is nilpotent).  Here is our main result on representation zeta functions.  Let $\twzeta{\nilgp}{p}{s}= \sum_{n=0}^\infty a_{p^n} p^{-ns}$.

\begin{theorem}
\label{thm:rattwst_intro}
 The $p$-local representation zeta functions $(\twzeta{\nilgp}{p}{s})_{p \textrm{ prime}}$ of a finitely generated nilpotent group $\Gamma$ are uniformly rational.
Moreover, the global representation zeta function $\zeta_\nilgp(s):= \sum_{n=1}^\infty a_n n^{-s}$ has rational abscissa of convergence.
\end{theorem}

The results in Sections~\ref{sec:zetagp} and \ref{sec:twst} both follow the same idea: we show how to interpret (uniformly and definably) in $\Qq_{p}$ the sets we want to count. More precisely, in Section~\ref{sec:zetagp} we show how to interpret in $\Th(\Qq_{p})$ the set of finite-index subgroups $H$ of $\nilgp$ and we show that the equivalence relations that arise are uniformly definable in $p$.  This allows us to apply Theorem~\ref{thm:rat}.  The same idea is used in Section~\ref{sec:twst}, but the details are more complicated.  We show how to interpret in $\Th(\Qq_p)$ the set of pairs $(N,\sigma)$, where $N$ is a finite-index normal subgroup of $\nilgp$ and $\sigma$ is an irreducible character of $\nilgp/N$, up to twisting by one-dimensional characters.  The key idea is first to interpret triples $(H,N,\chi)$, where $H$ is a finite-index subgroup of $\nilgp$, $N$ is a finite-index normal subgroup of $H$ and $\chi$ is a one-dimensional character of $H/N$---the point is that finite nilpotent groups are monomial, so any irreducible character is induced from a one-dimensional character of a subgroup.  The equivalence relation of giving the same induced character can be formulated in terms of restriction, and shown to be definable.  Inspecting these constructions shows that they are all uniform in $p$, so again Theorem~\ref{thm:rat} applies.

Since the first draft of this paper \cite{HruMar} was circulated, there has been considerable activity in the field of representation growth.  Jaikin-Zapirain \cite{Jai} used the coadjoint orbit formalism of Howe and Kirillov to parametrize irreducible characters of $p$-adic analytic groups; rationality of the representation zeta function then follows from the usual arguments of semi-simple compact $p$-adic integration.  Voll used similar ideas to parametrize irreducible characters of finitely generated torsion-free nilpotent groups, and showed that representation zeta functions are rational and satisfy a local functional equation \cite{VolAnn} (in fact, he proved this for a very general class of zeta functions that includes representation zeta functions and subgroup zeta functions as special cases).  Stasinski and Voll proved a uniformity result for representation zeta functions and calculated these zeta functions for some families of nilpotent groups \cite[Thms.~A and B]{StaVol}.  Ezzat \cite{EzzPhD}, \cite{Ezz2}, \cite{Ezz3} and Snocken \cite{SnoPhD} calculated further examples of representation zeta functions of nilpotent groups.  For work on representation growth for other kinds of group, see \cite{LubMar}, \cite{LarLub}, \cite{AKOVContemp}, \cite{AKOVDuke}, \cite{AKOVbasechange}, \cite{Avni}, \cite{AizAvn1}, \cite{AizAvn2}, \cite{BarHar}, \cite{BarBranch}.

The Kirillov orbit method has the advantage that it linearises the problem of parametrizing irreducible representations and simplifies the form of the imaginaries that appear.  The disadvantage is that the proof of rationality only applies to $\twzeta{\nilgp}{p}{s}$ for almost all $p$---one must discard a finite set of primes.  \emph{We stress that our result Theorem~\ref{thm:rattwst_intro} is the only known proof of rationality of $\twzeta{\nilgp}{p}{s}$ that works for every $p$.}

This paper falls naturally into two parts. The first part is model-theoretic: in Section~\ref{sec:elim} we establish an abstract criterion, Proposition~\ref{prop:EI/UFIcrit}, for elimination of imaginaries and apply it in Sections~\ref{sec:padics} and \ref{sec:asymptotic} to prove Theorems~\ref{thm:Qp} and \ref{thm:ultraprod Qp}. Section~\ref{sec:D3} consists of a study of unary types in henselian valued fields, which is used extensively in Sections~\ref{sec:padics} and \ref{sec:asymptotic}. In Section~\ref{sec:rat} we establish the general rationality result Theorem~\ref{thm:rat}, we prove Theorem~\ref{thm:abscissa_rat} and we show how the techniques developed in this paper can be used to prove transfer results between local fields of positive characteristic and mixed characteristic.

In the second part (Sections~\ref{sec:zetagp} and \ref{sec:twst}), we apply Theorem~\ref{thm:rat} to prove Theorems~\ref{thm:ratsbgp} and \ref{thm:rattwst_intro}.  The main tools are results from profinite groups; no ideas from model theory are used in a significant way beyond the notion of definability.  We finish Section~\ref{sec:twst} by using the Kirillov orbit formalism and Theorem~\ref{thm:rat} to recover a double uniformity result (Theorem~\ref{thm:double}) of Stasinski and Voll \cite{StaVol} for the representation zeta functions of groups of points of a smooth unipotent group scheme.

Finally, the Appendix contains an alternative proof, Theorem~\ref{count}, of the rationality results of Section~\ref{sec:rat} for fixed $p$ that generalizes to the analytic setting. An important application of this work of Cluckers is that it gives a tool for proving rationality of certain zeta functions associated to a compact $p$-adic analytic group: see the paragraph following Remark~\ref{rem:torsion}, for example. Here the methods of the main paper do not go through because one needs to use an extended language containing symbols for analytic functions, and elimination of imaginaries in this setting is known to require more sorts than just the geometric imaginaries.  (Note, however, that various rationality and uniformity results have been obtained for representation zeta functions of certain compact $p$-adic analytic groups using the Kirillov orbit method  \cite{AKOVDuke}.)

\bigskip

\noindent {\em Note:}\! A draft of this paper \cite{HruMar} has been available for over ten years now. Alongside the previous theorems concerning $\Qq_{p}$ for fixed $p$, the present version also contains new material on the model theory of ultraproducts of $\Qq_p$, which allows us to prove the uniformity as $p$ varies of the previous elimination of imaginaries and rationality theorems, as well as a transfer result between positive equicharacteristic and mixed characteristic.  There is extra material on representation growth and a new Appendix on cell decomposition methods.


\section{Elimination of imaginaries}
\label{sec:elim}

\subsection{Definition and first properties}

We denote by $\Nn$ ($\Nn_{>0}$) the non-negative (positive) integers, respectively.  For standard model-theoretic concepts and notation such as $\dcl$ (definable closure) and $\acl$ (algebraic closure) we refer the reader to any introduction to model theory, e.g., \cite{Mar} or the first chapter of \cite{HHM2}.  We will write interchangeably $\dcl(bb')= \dcl(b,b')= \dcl(\{b, b'\})$ and $\dcl(A,b) = \dcl(Ab) = \dcl(A\cup\{b\})$, etc.

\begin{notation}
If $X$ is a definable (possibly $\infty$-definable) set in some structure
$M$ and $A \subseteq M$, we will write $X(A) := \{a\in A\mid M\models
X(a)\}$. If we want to make the parameters of $X$ explicit, we will write $X(A;b)$.
\end{notation}

We say that the definable set $X$ is \emph{coded} (in $M$) if it can
be written as $R(M; b)$, where $b$ is a tuple of elements of $M$, and
where $b \neq b'$ implies that $R(M; b) \neq R(M; b')$. In this
situation $\dcl(b)$ depends only on $X$ and is called a {\em code} for
$X$. It is denoted $\code{X}$. We say $T$ {\em eliminates imaginaries}
(EI) if every definable set in every model of $T$ is
coded. Equivalently, if there are at least two constants, $T$
eliminates imaginaries if for any $\emptyset$-definable equivalence
relation $E$ there is an $\emptyset$-definable function whose fibers
are exactly the equivalence classes of $E$ (cf.\
\cite[Lemme\,2]{PoiImGal}).

For any theory $T$, by adding sorts for every $\emptyset$-definable quotient we obtain
a theory $\eq{T}$ that has elimination of imaginaries. These new sorts are called \emph{imaginary sorts} and the old sorts from $T$ are called the \emph{real sorts}. Similarly, to any model $M$ of $T$ we can associate a (unique) model $\eq{M}$ of $\eq{T}$ that has the same real sorts as $M$. In general, we use the notation $\code{X}$ to refer to the code of $X$ with respect to $\eq{T}$. We will denote by $\dcleq$ the definable closure in $\eq{M}$ and similarly for $\acleq$.

We will consider many-sorted theories  with a distinguished collection
$\mathcal{S}$ of sorts, referred to as the \emph{dominant} sorts; we   assume that for
 any sort $S$, there exists a $\emptyset$-definable partial function from a finite product of dominant sorts  onto $S$ (and this function is viewed as part of the presentation of the theory). The set of elements of dominant sorts in a model $M$ is denoted $\dom(M)$.

The following lemma and remark---which reduce elimination of imaginaries to coding certain functions---will not be used explicitly in the $p$-adic case, but they are an essential guideline as unary functions of the kind described in the remark are central to the proof of Proposition~\ref{prop:EI/UFIcrit}.

\begin{lemma}[1](cf.\ \cite[Remark\,3.2.2]{HHM})
A theory  $T$ admits elimination of imaginaries if every function definable (with parameters) whose domain is contained in a single dominant sort is coded in any model of $T$.
\end{lemma}

\begin{proof}
  Since encoding a set is equivalent to encoding the identity function
  on this set, it suffices to show that every definable function $f$
  is coded. Pulling back by the given $\emptyset$-definable functions,
  it suffices to show that every definable function whose domain is
  contained in a product $M_1 \times \ldots \times M_n$ of dominant
  sorts is coded.  For $n=1$, this is our assumption.  For larger $n$,
  we use induction, regarding a definable function $f : M_1 \times
  \ldots \times M_n \to \cpow{M}{k}$ as the function $f'$ mapping
  $c\in M_{1}$ to the code of the function $y \mapsto f(c,y)$. By
  compactness there are $<f>$-definable functions $h_{i}$ covering
  $f'$. The codes of these $h_{i}$ allow us to code $f$.
\end{proof}

\begin{remark}[1sharp]
In  Lemma~\ref{lem:1}, we don't need to be able to encode all
definable functions whose domain is contained in a single dominant
sort. Such functions are said to be unary. For \(T\) to eliminate
imaginaries, it suffices that:
\begin{enumerate}
\item Every unary definable function \(f\) which is the identity on
  its domain is coded. This is equivalent to unary EI (i.e., the property that every
  definable subset of a single dominant sort is coded).
\item For all \(M\models T\), every \(e\in \eq{M}\), every $A\subseteq M$, every unary \(Ae\)-definable function \(f_e\) and  every non-empty \(A\)-definable set \(D\), the following holds: if \(A\supseteq\dcleq(e)\cap M\), $e
  \in \dcleq(A,f_e(c))$ for any $c \in D$, and $\tp(e/A)$ implies the
  type of $e$ over $Ac$ for any $c \in D$, then \(f_e\) restricted to
  \(D\) is \(A\)-definable.
\end{enumerate}

\noindent Indeed, let $e$ be imaginary. There exist
$c_{1},\ldots,c_{n}\in\dom(M)$ such that
$e\in\dcleq(c_1,\ldots,c_n)$. Let $A_{l} =
\dcleq(e,c_1,\ldots,c_l)\cap M$. We know that $e\in \dcleq(A_{n})$ and
we want to show that $e\in\dcleq(A_{0}) = \dcleq(\dcleq(e)\cap M)$;
i.e., $e$ is interdefinable with a tuple of real elements.

Let us proceed by reverse induction. Suppose $e\in\dcleq(A_{l+1})$,
let $A = A_{l}$ and let $c = c_{l+1}$. Then \(e \in \dcleq(A_{l+1}) =
\dcleq(\dcleq(e,c_1,\ldots,c_{l+1})\cap M) = \dcleq(\dcleq(Aec)\cap
M)\). So we can find $d=f(e,c)$ and $e = h(d)$ for
some $A$-definable functions $f,h$. By unary EI and since
\(\dcleq(Ae)\cap M = A\), any $Ae$-definable
subset of a dominant sort is already $A$-definable. Thus, by
hypothesis, $e=h(f(e,c'))$ for any $c' \models \tp(c/A)$. Let $D$ be
an $A$-definable set with $c \in D$ and such that $e=h(f(e,c'))$ for
any $c' \in D$. Note also that, by unary EI again, for any \(c\in D\),
\(\tp(c/A)\) implies \(\tp(c/Ae)\) and thus \(\tp(e/A)\) implies
\(\tp(e/Ac)\). It follows from hypothesis 2\ that the map \(f_e : x\mapsto f(e,x)\)
restricted to \(D\) is \(A\)-definable and that $e\in\dcleq(A) =
\dcleq(A_{l})$.
\end{remark}

\begin{definition}[EI/UFI]
We will say that a theory $T$ \emph{eliminates imaginaries up to uniform finite imaginaries (EI/UFI)} if for all $M\models T$ and $e\in \eq{M}$, there exists a tuple $d\in M$ such that $e\in \acleq(d)$ and $d\in\dcleq(e)$.

The theory $T$ is said to \emph{eliminate finite imaginaries} (EFI) if any $e\in\acleq(\emptyset)$ is interdefinable with a tuple from $M$.
\end{definition}

Let us now give a criterion for elimination of imaginaries from \cite{HruFin}.

\begin{lemma}[critHfin]
A theory $T$ eliminates imaginaries if it eliminates imaginaries up to uniform finite imaginaries and for every set of parameters $A$, $T_{A}$ eliminates finite imaginaries.
\end{lemma}

\begin{proof}
Let $e\in \eq{M}\models \eq{T}$.  Then by EI/UFI, there exists $d\in M$ such that $e\in \acleq(d)$ and $d\in\dcleq(e)$. Hence $e$ is a finite imaginary in $T_{d}$ and there exists $d'\in M$ such that $e\in\dcleq(dd')$ and $dd'\in\dcleq(ed) = \dcleq(e)$, i.e., $e$ is coded by $dd'$.
\end{proof}

\subsection{Valued fields}
\label{sec:Vfield}

If $F$ is a field then we denote by $\alg{F}$ the algebraic closure of $F$.  Let $L$ be a valued field, with valuation ring $\Val(L)$,  maximal ideal
$\mathcal{M}(L)$ and residue field $\res(L)$. We will find it convenient to
consider the value group $\Valgp(L)$ in both an additive notation (with
valuation $\val : L \to \Valgp(L)\cup\{\infty\}$) and a multiplicative
notation (with reverse order and absolute value $|.|$), depending on the setting. We will consider valued fields in the geometric language whose sorts (later referred to as the geometric sorts) are as follows. We take a single dominant sort $\K$, for $L$ itself. The additional sorts $\Latt[n],\Tor[n]$ for $n\in \Nn$ are given by

$$\Latt[n] :=   \GL{n}(\K) /  \GL{n}(\Val) =  \triang_{n}(\K) / \triang_n(\Val),$$
the set of  lattices in $\cpow{\K}{n}$, and
$$\Tor[n] := \GL{n}(\K) /  \GL{n,n}(\Val)   = \bigcup_{m \leq n}\triang_n(\K)/\triang_{n,m}(\Val) =   { \bigcup _{e \in \Latt[n]} e/\mathcal{M}e }.$$

\noindent Here a \emph{lattice} is a  free $\Val$-submodule of $\cpow{\K}{n}$ of rank $n$,
  $\triang_n$ is the group of invertible upper triangular matrices,
  $\GL{n,m}(\Val)$ is the group
of matrices in $\GL{n}(\Val)$ whose $m$th column reduces mod $\mathcal{M}$ to
the column vector of $\res$ having a one in the $m$th entry and zeroes elsewhere, and  $\triang_{n,m}(\Val):=\triang_n(\Val) \cap \GL{n,m}(\Val)$.  There is a canonical map from $\Tor[n]$ to $\Latt[n]$ taking $f \in e/\mathcal{M}e$ to the lattice $e$.

  It is easy to see, using
elementary matrices, that $\GL{n}(\K) = \triang_{n}(\K) \GL{n}(\Val)$, justifying the equivalence
 of the first two definitions of $\Latt[n]$.  Equivalently, it is shown in \cite[Lemma\,2.4.8]{HHM} that every lattice has a basis in triangular form.

 Note that there is an obvious injective $\emptyset$-definable function $\Latt[m] \times \Latt[m'] \to \Latt[m+m']$, namely $(\lambda,\lambda') \mapsto \lambda \times \lambda'$, so we can identify any subset of a product of $\Latt[n_{i}]$ with a subset of $\Latt[n]$, where $n = \sum_{i} n_{i}$.

Note also that $\Latt[1]$ can be identified with $\Valgp$ by sending the coset $c\inv{\Val}$ to $v(c)$. Then $\res$ can be identified with the fiber of $T_1 \to \Latt[1]$ above the zero element of $\Latt[1] = \Valgp$. More generally, let $\cBallSet = \{ \{x\mid \val(x-a) \geq \val(b) \}\mid a,b \in
\K\}$  and let $\oBallSet = \{ \{x\mid \val(x-a) > \val(b) \}\mid a,b \in \K\}$ be the
sets of \emph{closed} (respectively \emph{open}) balls with center in $\K$ and radius in $v(\K)$. Then $\cBallSet$ embeds into $\Latt[2] \cup \K$ and $\oBallSet$ into
$T_{2}$. Indeed, the set of closed balls of radius $+\infty$ is identified with $\K$.
The group $G(\K)$ of affine transformations of the line acts
transitively on the closed balls of nonzero radius; the stabilizer of
$\Val\in \cBallSet$ is $G(\Val)$. Embedding $G(\K)$ in
\(\GL{2}(\K)\) as upper triangular matrices, we get $\cBallSet \setminus \K \iso  G(\K) / G(\Val) \subseteq
\GL{2}(\K) / \GL{2}(\Val)$. The group $G(\K)$ also acts transitively on $\oBallSet$
and the stabilizer of $\mathcal{M}\in\oBallSet$ is $G(\K)\cap \GL{2,2}(\Val)$. We
will write $\BallSet := \cBallSet\cup\oBallSet$ for the set of all balls.
Note, however, that if $\Valgp$ has a smallest positive element, the open
balls are also closed balls.

In Sections~\ref{sec:D3} and \ref{sec:asymptotic}, we will also consider the sort $\RV :=
\inv{\K}/(1+\mathcal{M})$ and the canonical projection $\rv : \inv{\K} \to \RV$. The structure on $\RV$ is given by its group structure and the structure induced by the exact sequence $\inv{\res}\to \RV\to \Valgp$, where the second map is denoted
$\val_{\rv}$---i.e., as well as as the group structure on \(\RV\), we have a binary predicate interpreted as $\val_{\rv}(x)
\leq \val_{\rv}(y)$, a unary predicate interpreted as $\inv{\res}$, and the ring
structure on $\res$ (adding a zero to $\inv{\res}$). This exact
sequence induces on each fiber of $\val_{\rv}$ the structure of a $\res$-vector
space (if we add a zero to the fiber). Any $T\supseteq \HFO$ (the
theory of henselian valued fields with residue field of characteristic
zero) eliminates field quantifiers by \cite{BasHF0}. It follows from this quantifier elimination result that $\RV$ is stably embedded and the structure induced
on $\RV$ is exactly the one described above. Note that we can identify
$\RV$ with $T_{1}$ if we add a zero to each fiber of $\val_{\rv}$.

The theory of a structure is determined by the theory of the dominant
sorts; so, for any field $L$ we can speak of $\Th(L)$ in the geometric
sorts. We take the geometric language $\LG$ to include the ring
structure on the sort $\K$ and the natural projections $\GL{n}(\K) \to
\Latt[n](\K)$ and $\GL{n}(\K) \to \Tor[n](\K)$.

In \cite{HHM}, it is shown that
$\ACVF$ eliminates imaginaries in $\LG$.    Let us now give the counterpart of this theorem for $p$-adic fields.

By $\LGminus$, we denote the restriction of $\LG$ to the sorts $\K$ and $\Latt[n]$. For each subset $\mathcal{N} \subseteq\Nn_{>0}$, we will also consider an expansion $\LGn$ of $\LGminus$ by a constant $a$ and for all $n\in \mathcal{N}$ a tuple of constants $c_{n}$ of length $n$ in the field sort.

By a \emph{uniformizer} of a valued field we mean an element $a$ whose
valuation is positive and generates the value group.  By an {\em
  (unramified) $n$-Galois uniformizer} we mean a tuple $c$ of elements
of the valuation ring such that \(\sum_{i=0}^{n}\res(c_i)\omega_n^i\)
generates
\(\Gg_{m}({\res}[\omega_{n}])/\gpow{\Gg_{m}({\res}[\omega_{n}])}{n}\)
where $\omega_{n}$ is some primitive $n$th root of unity and
$\gpow{G}{n}$ denotes the $n$th powers in the group $G$. If
${\res}[\omega_{n}]$ has degree smaller than \(n\) over $\res$, take the coordinates of \(c\) beyond the \([{\res}[\omega_n]:{\res}]\)th entry to be zero.

Let $\PLO$ be the theory of pseudo-local fields of residue
characteristic $0$, i.e., henselian fields with value group a
$\Zz$-group (i.e., an ordered group elementarily equivalent to
$(\Zz,0,+,<)$) and residue field a pseudo-finite field of
characteristic $0$. By \cite[Theorem 8]{Ax}, any pseudo-finite field is elementarily equivalent to an ultraproduct of fields \(\Ff_p\) for $p$ prime, so $\PLO$ is the theory
of ultraproducts $\prod \Qq_{p}/\mathcal{U}$ of $p$-adics over
non-principal ultrafilters on the set of primes.

Furthermore, let $\mathfrak{F}_{p}$ be a set of finite extensions of $\Qq_{p}$ and let $\mathfrak{F}= \bigcup_{p}\mathfrak{F}_{p}$. Any ultraproduct $\prod_{L\in\mathfrak{F}} L/\mathcal{U}$ of residue  characteristic zero---i.e., such that the ultrafilter $\mathcal{U}$ does not contain any set included in some $\mathfrak{F}_{p_{0}}$---is a model of $\PLO$. Note that if $\mathfrak{F}_{p}$ is nonempty for infinitely many $p$ then there exists an ultrafilter $\mathcal{U}$ on $\mathfrak{F}$ such that $\prod_{L\in\mathfrak{F}} L/\mathcal{U}$ has residue characteristic zero.

Let $L$ be a valued field, regarded as an $\LGminus$-structure, and
let $p$ be its residue characteristic. Fix
\(\mathcal{N}\subseteq\Nn_{>0}\). A {\em proper expansion} of $L$ to
$\LGn$ is a choice of $a$ and tuples $c_{n}$ for each $n\in \mathcal{N}$ such that:
\begin{enumerate}
\item $a$ is a uniformizer;
\item if $n$ is prime to $p$ or $p=0$ then $c_{n}$ is an unramified $n$-Galois uniformizer if one exists, and $0$ otherwise;
\item if $p \neq 0$, $p$ divides $n$ and $n > p$ then $c_n$ is a tuple
  of zeros;
\item if $p\neq 0$, $p\in \mathcal{N}$ and $L$ is not a finite
  extension of $\Qq_p$ then $c_p$ is a tuple of zeros;
\item if $p\neq 0$, $p\in \mathcal{N}$ and $L$ is a finite extension of $\Qq_p$ then the first coordinate of $c_p$ is a generator of $L$ over $\Qq_{p}$ that is algebraic over $\Qq$, and the other coordinates are zero.
\end{enumerate}
The point of 5.\ is to ensure we have a constant that generates $L$
over $\Qq_p$ in the local field case when $p\in \mathcal{N}$.

Note that because there are only finitely many possibilities for the minimal polynomial of $\omega_{n}$ over $\res$, the class of proper expansions to $\LGn$ of models of $\PLO$ is elementary. Let us denote this class by $\PLO^{\mathcal{N}}$. Note also that a residue characteristic zero ultraproduct of proper expansions to $\LGn$ of $L\in\mathfrak{F}$ is a model of $\PLO^{\mathcal{N}}$.

Here is a precise statement of the two main elimination of imaginaries results of this paper. The first is for finite extensions of $\Qq_{p}$:

\begin{theorem}[Qp]
The theory of $\Qq_p$ eliminates imaginaries in $\LGminus$. The same is true for any finite extension $L$ of $\Qq_p$, provided one adds a constant symbol for a generator of $L\cap\alg{\Qq}$ over $\Qq_{p}\cap\alg{\Qq}$.
\end{theorem}

\noindent The second is for their ultraproducts of residue characteristic zero:

\begin{theorem}[ultraprod Qp]
 $\PLO^{\Nn}$ eliminates imaginaries in $\LGm$.
\end{theorem}

Note that elimination of imaginaries in an incomplete theory is
equivalent to elimination of imaginaries in all of its completions. It
follows that elimination of imaginaries is uniform over all
pseudo-local fields and hence over local fields of large residue
characteristic (see Corollary\,\ref{cor:EIunif}).

\begin{remark}[rem thm]
\begin{thmenum}
\item\label{no Tn} Although the $T_{n}$ are needed to obtain EI in algebraically closed fields, they are not needed here. Indeed, if a valued field $K$ has a discrete valuation (i.e., the value group has a smallest positive element $\val(\unif)$), then for any lattice $e$, $\unif e$ is itself a lattice, and a coset $h$ of $\unif e$--- a typical element of $\Tor[n]$---can be coded by the lattice in $\cpow{K}{n+1}$ generated by  $h\times\{1\}$. Hence all elements of $T_{n}(K)$ are coded in $\Latt[n+1](K)$.
\item\label{constants Qp} As we will see in Section\,\ref{sec:padics}, to
  obtain elimination of imaginaries in a finite extension \(L\) of
  \(\Qq_p\), we need to add constants for elements of a subfield
  $F\subseteq L$ with a certain number of properties:
  \begin{thmenum}
  \item $F$ contains a uniformizer;
  \item $\resf(F) = \res(L)$;
  \item $\alg{L} = \alg{F}L$ (in fact, we need that for every finite extension $K$ of $L$ there is a generator of $\Val(K)$ whose minimal polynomial is over $F$).
\end{thmenum}
Note that it suffices to take $F =  \Qq[c]$, where $c$ generates $L$
over $\Qq_p$. Note that we can choose such a $c$ that belongs to
$\alg{\Qq}$, hence the statement of Theorem\,\ref{thm:Qp}.

Note also that a proper expansion of some finite extension $L$ of
$\Qq_{p}$ to $\LGm$, contains a generator (named by a constant) of $L$ over $\Qq_p$. Hence
such proper expansions of $L$ eliminate imaginaries in $\LGm$.

\item\label{constants asymp} To prove elimination of imaginaries in a pseudo-local field \(L\), we need to name in Section\,\ref{sec:asymptotic} elements of a subfield $F\leq L$ which satisfies  (a), (c) as above and the following conditions:
\begin{thmenum}\setcounter{enumii}{3}
\item $\resf(F) \gpow{(\inv{\res(L)})}{n} = \res(L)$ for all $n$;
\item $\res(L)$ admits EI in the language of rings augmented by constants for elements of $\resf(F)$.
\end{thmenum}

Let us show that we can choose $F$ to be generated by a uniformizer $a$ and unramified $n$-Galois uniformizers $c_{n}$ for all $n$. It is clear that such an $F$ satisfies (a). Furthermore, $\alg{\res(L)} = \alg{\resf(F)}\res(L)$. Indeed, let $\omega_{n}$ be a primitive $n$th root of unity, and let $d_{n} = \sum_{i} c_{n,i}\omega_{n}^{i}$. The degree $n$ extension of $\res(L)[\omega_{n}]$ is contained in $\res(L)[\omega_{n},\sqrt[n]{d_{n}}]$ by Kummer theory and it contains the degree $n$ extension of $\res(L)$.

Now (c) is a consequence of (a) and the above statement and (e) also follows as any extension of degree $n$ is generated by an element in $\alg{\resf(F)}$, so there is an irreducible polynomial of degree $n$
with $\resf(F)$-definable parameters; this is the hypothesis of \cite[Proposition\,B.(3)]{ChaHruACFA}. Finally for any $n$, there is a $d$ such that $\{x\in \res(L)\mid x^{n} = 1\} = \{x\in \res(L)\mid x^{d} = 1\}$ and $\res(L)$ contains primitive $d$th roots of unity. Then $c_{d}\in \res(L)$ generates $\inv{\res(L)}/\gpow{(\inv{\res(L)})}{d} = \inv{\res(L)}/\gpow{(\inv{\res(L)})}{n}$, so (d) holds.
\item It would be nice to find a more precise description of the imaginaries if no constants are named. For finite extensions of $\Qq_{p}$, this is done in Remark~\ref{rem:without constant}.
\end{thmenum}
\end{remark}

Before going any further, let us show that Theorem~\ref{thm:ultraprod Qp} allows us to prove a uniform version of Theorem~\ref{thm:Qp}.

\begin{corollary}[EIunif]  Let $\mathfrak{F}_{p}$ be any set of finite
  extensions of $\Qq_{p}$ and let
  $\mathfrak{F}=\bigcup_{p}\mathfrak{F}_{p}$. Let $\phi(x,y)$ be an
  $\LGm$-formula (where $x,y$ range over \(\emptyset\)-definable sets
  $X,Y$). Then there exist integers $m$, $l$, a set of integers
  $\mathcal{N}$, a prime $p_{0}$ and some $\LGn$-formula $\psi(x,w)$
  such that the following uniform statement of elimination of
  imaginaries holds. For all $p\geq p_{0}$ and all proper expansions
  to $\LGn$ of $L_{p}\in\mathfrak{F}_{p}$, $\psi(x,w)$ defines a
  function \[ f_{L_{p}}: X \to \Latt[m](L_{p}) \times
  \cpow{\K(L_{p})}{l}   \]  and \[L_p \models (\forall x,x')
  (f_{L_{p}}(x)= f_{L_{p}}(x') \iff [(\forall y)\phi(x,y) \iff \phi(x',y)
  ]). \]
\end{corollary}

\begin{proof}
Assume $\mathfrak{F}_{p}$ is nonempty for infinitely many $p$,
otherwise the statement is trivial. The formula $\forall
y\,\phi(x,y)\iff \phi(x',y)$ defines an equivalence relation in any
ultraproduct $L$ of fields in $\mathfrak{F}$. By
Theorem~\ref{thm:ultraprod Qp}, there is a formula $\psi(x,w)$ (which
works for any proper expansion to $\LGm$ of any such ultraproduct of
residue characteristic zero) such that, in every model of \(\PLO^{\Nn}\), $\psi(x,w)$ defines a function $f$ and $f(x) = f(x')$ if and only if $\forall y\,\phi(x,y)\iff \phi(x',y)$.

Let us now assume there is an infinite set $I\subseteq\mathfrak{F}$ such that $I$ has a nonempty intersection with infinitely many $\mathfrak{F}_{p}$ and for every $L\in I$, there is a proper expansion of $L$ to $\LGm$ such that we do not have $f(x) = f(x')$ if and only if $\forall y\,\phi(x,y)\iff \phi(x',y)$ in $L$. Then there exists an ultrafilter on $\mathfrak{F}$ containing $I$ but containing no set included in some $\mathfrak{F}_{p_{0}}$ and such that $\prod_{L\in\mathfrak{F}} L/\mathcal{U} \models \PLO^{\Nn}$; but we do not have $f(x) = f(x')$ if and only if $\forall y\,\phi(x,y)\iff \phi(x',y)$ in this ultraproduct, a contradiction.

By compactness, this equivalence also holds in proper expansions to
$\LGn$, for some finite $\mathcal{N}$.
\end{proof}

\begin{remark}
\begin{thmenum}
\item In particular, whenever $\phi(x,y)$ is interpreted in $L_{p}$ as
  an equivalence relation \(xEy\), $f_{L_p}(x)$ codes the $E$-class of
  $x$.
 \item If $\mathfrak{F}_{p}$ is finite for all $p$ then, as $\bigcup_{p
     < p_{0}} \mathfrak{F}_{p}$ is finite, we can find, using
   Theorem~\ref{thm:Qp} and Remark~\ref{constants Qp}, a $\psi$ and an
   $\mathcal{N}$ that work for all $L \in \bigcup_{p} \mathfrak{F}_{p}$
   and not just for sufficiently large $p$.
 \end{thmenum}
 \end{remark}

The proof of Theorems~\ref{thm:Qp} and \ref{thm:ultraprod Qp} uses elimination of imaginaries  and the existence of invariant extensions in the theory of algebraically closed valued fields. Recall that a theory $T$ has the \emph{invariant extension property} 
if whenever $A = \acleq(A) \subseteq M \models T$
and $c \in M$,
 $\tp(c/A)$ extends to an $\aut(\eq{M}/A)$-invariant type over $M$.
 This holds trivially for any finite field, and by inspection, for $\Th(\Zz,+,<)$ and, although we will only use a weaker version of the extension property (Corollary~\ref{cor:D3.3}) in the proof of Theorem~\ref{thm:Qp}, we will show that the theory of a finite extension of $\Qq_{p}$ (with the geometric sorts) enjoys the stronger version (Remark~\ref{rem:inv ext}).

\subsection{Real elimination of imaginaries}

 To illustrate the idea of transferring imaginaries from one theory to the other, consider the following way of deducing EI for RCF (the theory of real closed fields) from EI for ACF (the theory of algebraically closed fields).

\begin{example}[RCF]
 Let $F$ be a field considered in a language extending the language of
 rings. Assume for all $M\models \Th(F)$:

 \begin{itemize}
  \item[(i)]   (Algebraic boundedness): If $A\subseteq M$ then
    $\acl(A) \subseteq \alg{A}\cap M$;
  \item[(ii)]  (Rigidity of finite sets): No automorphism of $M$ can have a  finite cycle of size $>1$.  Equivalently,
for each $n$, there exists $\emptyset$-definable functions
$r_{i,n}(x_1,...,x_n)$ that are symmetric in the $x_i$, such that for
any set \(S\) of size \(n\), $S = \{r_{1,n}(S),\ldots,r_{n,n}(S) \}$.  (Here $r_{i,n}(S)$ denotes $r_{i,n}(x_1,...,x_n)$ when $S= \{x_1,\ldots,x_n\}$, possibly with repetitions.);
  \item[(iii)] (Unary EI):   Every definable subset of $M$ is coded.
 \end{itemize}
 Then $\Th(F)$ eliminates imaginaries (in the single sort of field elements).
\end{example}

\begin{proof}
 Let $f:M\to M$ be a definable function.  By Lemma~\ref{lem:1}, it suffices to prove that $f$ is coded.
Let $H$ be the Zariski closure (over $M$) of the graph of $f$.  Since the theory is algebraically bounded,
the set $H(x):= \{y\mid (x,y) \in H \}$ is finite for any $x$, of size bounded by some $n$.  Let $U_{n,i} $ be the set of $x$ such that
 $f(x) =   r_{i,n}(H(x))$.  Then, by elimination of imaginaries in ACF, $H$---being a Zariski closed set---is coded in $\alg{M}$. But the code is definable over $M$ and hence is in the perfect closure of $M$. Replacing
 this code with some $p^{n}$th power in the characteristic $p$ case, we can suppose it belongs to $M$. Moreover, each $U_{n,i}$ (being unary) is coded; these codes together give a code for $f$.
\end{proof}

Note that RCF satisfies the hypotheses of Example~\ref{ex:RCF}, but $\Th(\Qq_p)$ (in the field sort alone) does not. More precisely, as shown in the introduction, the value group cannot be definably embedded into $\cpow{\Qq_p}{n}$. Hence hypothesis (iii) fails for $\Th(\Qq_p)$ in the field sort alone.

\begin{remark} If $F$ is a field satisfying (i), (iii), then $F$ has EI/UFI.
This is an immediate consequence of Proposition~\ref{prop:EI/UFIcrit} as
hypotheses (ii) and (iv) of Proposition~\ref{prop:EI/UFIcrit} are true if $\tT$ is the theory of algebraically closed fields in the language of rings.\end{remark}

\subsection{Criterion for elimination of imaginaries}

Let   $\tT$ be a complete theory in a language $\tL$.
Assume $\tT$  eliminates quantifiers and imaginaries. Let $\nT$ be a complete  theory in a
language $\nL \supseteq \tL$;  assume $\nT$ contains the universal part of $\tT$.

In a model $\nM$ of the theory $\nT$, three kinds of definable closure
can be considered: the usual definable closure $\DCL$; the definable
closure in $\eq{\nM}$, denoted $\DCLeq$; and the imaginary definable
closure restricted to real points (that is, if \(A\subseteq\eq{M}\),
the set $\DCLeq(A)\cap \nM$). As $\DCLeq(A)\cap\nM$ and $\DCL(A)$ take
the same value on any set of real points, we will denote them both by
$\DCL(A)$. One must take care however that if $A$ contains imaginary
elements, $A\not\subseteq\DCL(A)$.

As $\tT$ eliminates imaginaries, these
distinctions are not necessary in models of $\tT$ and we will only need $\dcll$. One
should note that, as $\tT$ eliminates quantifiers, $\dcll$ is the closure
under quantifier-free $\tL$-definable functions and hence that,
for any $A\subseteq \nM$, $\dcll(A)\cap \nM\subseteq \DCL(A)$.

Analogous statements hold for $\acll$, $\ACL$, $\ACLeq$, $\tpl$, $\tpL$, etc.

One should also be careful that if $\nM\models \nT$ is contained in some $\tM\models\tT$, there is no reason in general that $\eq{\nM}$ should be contained in $\tM$. In fact, the whole purpose of the following proof is to show that under certain hypotheses every element of $\eq{\nM}$ is interdefinable with a tuple in $\nM$.

\begin{proposition}[EI/UFIcrit] 
Assume $\tT$ and $\nT$ have the properties given above. Let $\nM$ be an $\card{L}^{+}$-saturated and $\card{L}^{+}$-homogeneous model of $\nT$ and let $\tM \models \tT$ be such that $\Langrestr{\nM}{\tL} \substr_{\tL} \tM$ and such that any automorphism of $\nM$ extends to an automorphism of $\tM$. If conditions (i)--(iv) below hold for any $A =\ACL(A) \subseteq \nM$ and any $c \in \dom(\nM)$, then  $\nT$ admits elimination of imaginaries up to uniform finite imaginaries (see Definition~\ref{def:EI/UFI}).

 \begin{itemize}
  \item[(i)] (Relative algebraic boundedness)  For every $\nM'\prec \nM$,
 $\DCL(\nM'c) \subseteq   \acll(\nM'c) $.
  \item[(ii)]    (Internalizing $\tL$-codes)  
  For all $\epsilon \in \dcll(\nM)$, there exists a tuple $\eta$ of
  elements of $\nM$ such that an automorphism of $\tM$ that stabilizes
  $\nM$ setwise fixes $\epsilon$ if and only if it fixes $\eta$.
  \item[(iii)]   (Unary EI)    Every $\nL(\nM)$-definable unary subset of $\dom(\nM)$  is coded in $\nM$.  
  \item[(iv)]   (Invariant types)    There exists an $\aut(\tM/A)$-invariant type $\tP$ over $\tM$ such that $\tprestr{\tP}{\nM}$ is consistent with $\tpL(c/A)$.    

  Moreover,  for any $\tL(\tM)$-definable function $r$ whose domain contains $\tP$, let $\germ{\tP}{r}$ be the $\tP$-germ
of $r$ (where two $\tL(\tM)$-definable functions $r,r'$ have the same $\tP$-germ if they agree on a   realization of $\tP$ over $\tM$).  Then:

$(*)$  There exists a directed order $I$ and a sequence
$(\epsilon_i)_{i\in I}$, with $\epsilon_{i}\in \dcll(A,\code{r})$ such
that $\sigma \in \aut(\tM/A)$ fixes $\germ{\tP}{r}$ if and only if
$\sigma$ fixes almost every $\epsilon_i$---i.e., $\sigma$ fixes
$\epsilon_i$ for all $i \geq i_0$, for some $i_0\in I$.
  \end{itemize}
\end{proposition}

Some comments on the proposition:
\begin{enumerate}
\item There are two ways to ensure that automorphisms of $\nM$ extend to automorphisms of $\tM$. The first is to take $\tM$ sufficiently homogeneous. The other is to take $\tM$ atomic over $\nM$; in the case of valued fields, we could take $\tM$ to be the algebraic closure of $\nM$.
\item In fact, we will only need (iv) for $\card{A} \leq \card{\nL}$.
\item If $\tP$ is definable then, for a uniformly defined family of functions $r_{b}$, $\germ{\tP}{r_{b}}$ is an imaginary (and we could take $\epsilon_i$ to be that imaginary). Nevertheless, if $\tP$ is not definable and say $(\epsilon_{i})$ is countable then Condition (iv) implies that the germ
is a $\Sigma^0_2$-hyper-imaginary, i.e., an equivalence class of sequences indexed by $I$ where the equivalence relation is given by a countable union of countable intersections of definable sets (although each definable set will involve only a finite number of indices, the countable union of countable intersections can involve them all). In the case of $\ACVF$ one also has that $\sigma$ fixes $\germ{\tP}{r}$ if and only if
$\sigma$ fixes cofinally many $\epsilon_i$; in this case the equivalence relation is also a countable intersection of countable unions of definable sets, so it is $\Delta^0_2$.

  \item  Hypotheses (ii) and (iii) are special cases of elimination of imaginaries.  It  would be nice to move (iii) from the hypotheses to the conclusion,
  i.e., assuming only (i), (ii) and (iv), to show that every imaginary is ``equivalent'' to an imaginary of $\tM$ definable over $\nM$.
  \end{enumerate}

First let us clarify how $\aut(\nM)$ acts on $\dcll(\nM)$ as this action will be
used implicitly throughout the proof. Any $\ns\in \aut(\nM)$ can be extended to
an automorphism $\ts\in \aut(\tM)$ and all these
extensions are equal on $\dcll(\nM)$, hence we have a well-defined action of
$\ns$ on $\dcll(\nM)$ and the notation $\aut(\nM/B)$ makes sense even if
$B\subseteq\dcll(\nM)$. Similarly, if $\tP$ is an $\aut(\tM/B)$-invariant
type, $\aut(\nM/B)$ acts on $\tP$-germs of $\tL(\nM)$-definable functions.

We begin our proof with the elimination of finite sets:

\begin{lemma}[EFin]
Assume (ii) holds in Proposition~\ref{prop:EI/UFIcrit}. Then every finite set $E\subseteq \nM$ is coded.
\end{lemma}

\begin{proof}
By EI for $\tT$, the finite set $E$ is coded by a tuple $\epsilon\in \tM$;  $\epsilon$ may
 consist of elements in $\dcll(\nM)$ but outside $\nM$.  By (ii), there exists a tuple $\eta$ of elements of $\nM$ such that an automorphism of $\dcll(\nM)$ leaving $\nM$ invariant fixes $E$ if and only if it fixes $\epsilon$ if and only if it fixes $\eta$. Thus $\code{E}$ and $\eta\in \nM$ are interdefinable.
\end{proof}

\begin{proof}[Proposition~\ref{prop:EI/UFIcrit}]
Let $e\in \eq{\nM}$. For some $c_1,\ldots,c_n \in \dom(\nM)$, we have  $e \in \DCLeq(c_1,\ldots,c_n)$. Let $A_i = \DCL(e,c_1,\ldots,c_i) \subseteq \nM$.  The claim is that $e \in \ACLeq(A_0)$.
We have $e \in \ACLeq(A_n)$ and show by reverse induction on $l \leq n$ that $e \in \ACLeq(A_l)$.
Assume inductively that $e \in \ACLeq(A_{l+1})$.    Let $A= A_l, c=c_{l+1}$. It is easy to check that
 $$A  = \DCL(Ae).  $$
As $e \in \ACLeq(A_{l+1})$, for some tuple $d\in A_{l+1} = \DCL(Ace)$, some $\nL(A)$-definable function $f$ and
some $\nL(A)$-definable, finite-set-valued function  $g$, we have
 $$ e \in g(d),  d=f(c,e).$$

 \noindent Let $f_e(x)=f(x,e)$.   Let $\bA = \ACL(A)$ and let $\nP = \tpL(c/\bA)$.

 Let ${\nM_0} \prec \nM$ such that $\eq{\nM_{0}}$ contains $Ae$. Note that for
 all $c'$ in the domain of $f_e$, $f_{e}(c')\in\DCL(\nM_{0}c')$. By (i), there
 exists an $\tL({\nM_0})$-definable finite-set-valued function $\phi_{c'}$
 such that $f_e(c') {\in} \phi_{c'}(c')$. By compactness, for some finite
 set $I_0$ and $\tL({\nM_0})$-definable finite-set-valued functions $(\phi_{i})_{i\in
   I_{0}}$, the following holds: for any $c'$ in the domain of $f_e$,
 $f_e(c') {\in} \phi_i(c')$ for some $i \in I_0$. Let $\phi(x) = \bigcup_{i \in I_0} \phi_i(x)$;
 so  $f_e(c') {\in} \phi(c')$ for all $c'$ in the domain of $f_e$. Hence if
 $\Phi$ is the set of all $\tL(\nM)$-definable, finite-set-valued functions $\psi$
 with a domain containing that of $f_e$ and such that for all $c'$ in the
 domain of $f_e$,  $f_e(c') {\in} \psi(c')$, then $\Phi$ is nonempty.

Let $\tP$ be an $\aut(\tM/\bA)$-invariant type over $\tM$ extending
$\nP$, as in (iv).  For $m \in \Nn$, let $\Phi_m$ be the set of  all
$\tL(\nM)$-definable functions $\phi \in\Phi$ such that for $c \models
\tP$, $\phi(c)$ is an $m$-element set. Note that, since $\tP$ is a
complete type, $\Phi_m$ does not depend on $c$. Let $m$ be minimal such that $\Phi_m$ is nonempty. All $\phi \in
\Phi_m$ share the same $\tP$-germ. Indeed, if $\phi,\phi'$ do not have
the same $\tP$-germ, let $\phi''(x) := \phi(x) \cap \phi'(x)$. Then
\(\phi''\in\Phi\) and since for all \(c\models\tP\), $\phi(c) \neq
\phi'(c)$, \(\phi''(c)\) would lie in $\Phi_{m'}$ for some $m'<m$.
Pick $F_E \in \Phi_m$, defined over some $E \subseteq \nM$. By
construction, $F_E$ covers $f_e$,  $F_E$ is $\tL(E)$-definable, and
the $\tP$-germ of $F_E$ is invariant under $\aut(\nM/\bA e)$.

\begin{claim}[claim1]
The $\tP$-germ of $F_E$ is invariant under $\aut(\nM/\bA)$.
\end{claim}

\begin{proof}
Let $(\epsilon_i)$ be a sequence as in (iv), coding the germ of $F_{E}$ on $\tP$.
  Note that $\epsilon_i \in \dcll(\nM)$ (since $F_{E}$ is $\tL(E)$-definable and $E \subseteq \nM$).
By (ii), we may  replace $\epsilon_i$ with an element of $\nM$, without changing $\aut(\nM/\epsilon_i)$; we do so.

Now, almost all $\epsilon_i$ must be in $\ACL(\bA e)$.   For otherwise, by moving to a subsequence we may assume
all $\epsilon_i$ are outside $\ACL(\bA e)$.  So $\aut(\nM/\bA e \epsilon_i)$ has infinite index in $\aut(\nM/\bA e)$.  By Neumann's Lemma,
for any finite set $X$ of indices $i$, there exists $\tau \in \aut(\nM/\bA e)$ with $\tau(\epsilon_i) \neq \epsilon_i$, for all $i \in X$.  By compactness
(and homogeneity of $\nM$), there exists $\tau \in \aut(\nM/\bA e)$  with $\tau(\epsilon_i) \neq \epsilon_i$ for all $i$.  But then $\tau$ fails
to fix the $\tP$-germ of $F$, contradicting the $\aut(\nM/\bA e)$-invariance
of this germ.

So for almost all $i$, some finite set ${\mathcal{E}_i}$ containing $\epsilon_i$ is defined over $Ae$. By Lemma~\ref{lem:EFin}, the finite set ${\mathcal{E}_i}$ is coded in $\nM$.
But  $A  = \DCL(Ae)$, so $\mathcal{E}_i$ is defined over $A$.  Hence  $\epsilon_i \in \bA$,
i.e., $\epsilon_i $ is fixed by $\aut(\nM/\bA)$.  This
 being the case for almost all $i$,
the $\tP$-germ of $F_E$ is invariant under $\aut(\nM/\bA)$.
\end{proof}

\begin{claim}[claim2]
$e \in \ACLeq(A)$.
\end{claim}
\begin{proof}
It suffices to show that if $((e_i,E_i)\mid i \in \Nn)$ is an indiscernible
sequence over $\bA$ with $e_0 = e$ and $E_0 = E$, then $e_i=e_{j}$ for some $i\neq j$.  Let $c \models \tprestr{\tP}{A(E_i)_{i\in\Nn}}$ such that $c\models \nP$.
By (iii) and because $A = \DCL(Ae)$,   $\tpL(c/A) $ implies $\tpL(c/Ae)$; hence  $\tpL(e/A)$ implies $\tpL(e/Ac)$.   So $\tpL(e_i/Ac)=\tpL(e/Ac)$.

By Claim~\ref{claim:claim1},  the $\tP$-germs of the $F_{E_i}$ are equal; so
$F_{E_i}(c)$ is a finite set $F$ that does not depend on $i$.  But
$f(c,e_i) {\in} F$,  so $f(c,e_i)$ takes the same value on some infinite
set $I'$ of indices $i$.   Hence so does the finite set
$g(f(c,e_{i}))$. As $e \in g(f(c,e))$ and \(\tpL(e/Ac) = \tpL(e_i/Ac)\),
it follows that $e_i \in g(f(c,e_i))$,  so infinitely many $e_i$ lie in
the same finite set and $e_i=e_j$ for some $i \neq j$.
\end{proof}

We have just shown that $e$ lies in  $\ACLeq(A) = \ACLeq(A_{l})$. This concludes the induction. It follows that \(e\in\ACLeq(A_0) = \ACLeq(\DCLeq(e)\cap M)\) and Proposition~\ref{prop:EI/UFIcrit} is proved.\end{proof}

Let us now show that this first criterion can be turned into a criterion for elimination of imaginaries.

\begin{corollary}[EIcrit]
Let $\nT$ and $\tT$ be as in Proposition~\ref{prop:EI/UFIcrit} and let us suppose moreover that:
\begin{itemize}
  \item[(v)]  (Weak rigidity) For all $A = \ACL(A)$ and $c\in\dom(\nM)$, $\ACL(Ac)  \subseteq \DCL(Ac)$.
\end{itemize}
Then $\nT$ eliminates imaginaries.
\end{corollary}

\begin{proof}
 Let $e \in \eq{\nM}$ be an imaginary element.
 We have  $e \in \DCLeq(c_1,\ldots,c_n)$ for some $c_1,\ldots,c_n \in \dom(\nM)$.
Let $A_i = \ACL(e,c_1,\ldots,c_i) \subseteq \nM$.    Then $e \in \DCLeq(A_n)$; we show by reverse induction on $l \leq n$ that $e \in \DCLeq(A_l)$.
We assume inductively that $e \in \DCLeq(A_{l+1})$.    Let $A= A_l, c=c_{l+1}\in\dom(\nM)$.  It is easy to check that
 $$A  = \ACL(Ae)  $$
and that, for some tuple $d$, $$d \in A_{l+1} = \ACL(Ace), e \in \DCLeq(Ad).$$
By Proposition~\ref{prop:EI/UFIcrit}, $e \in \ACLeq(A_0)$, so  $d \in \ACL(Ac)$.
By weak rigidity (v), $d \in \DCL(Ac)$.  Thus $e \in \DCLeq(Ac)$.

 Say $e = h(c)$, where $h$ is an $\eq{\nL}(A)$-definable function.  Then $h^{-1}(e)$ is an $\nL(\nM)$-definable subset of $\dom(\nM)$,
 hence by (iii) it has a code $e' \in \nM$.  Clearly $e$ and $e'$ are interdefinable over $A$.   As $e' \in \nM$, we have $e' \in \DCL(Ae)=A$.
 So $e \in \DCLeq(A)=\DCLeq(A_l)$.  This finishes the induction and shows that $e \in \DCLeq(A_0)$.

Let  $a$ be a tuple from $A_0$ such that $e$ is $\eq{\nL}(a)$-definable.  Let
 $a'$ be the (finite) set of conjugates of $a$ over $e $.  Then $\DCLeq(e) = \DCLeq( a')$ and, by Lemma~\ref{lem:EFin}, $a'$ is coded, hence $e$ is interdefinable with some sequence from $\nM$.
\end{proof}

Keeping (v) out of Proposition~\ref{prop:EI/UFIcrit} makes the proof of the EI criterion messier than strictly necessary. Nonetheless, distinguishing the case without (v) is important for ultraproducts of the $p$-adics where (v) fails.

The following lemma will be used to prove (v) in the $p$-adic case.

\begin{lemma}[por]
 Assume that for any $a \in \nM$, there exists an $\aut(\nM/\ACL(a))$-invariant type \(\nP\) over $\nM$ and an \(\emptyset\)-definable function \(f\) such that \(\nP(x)\vdash f(x) = a\). Then (v) follows from:
 \begin{itemize}
  \item[(v$'$)]  If $B \subseteq \dom(\nM)$ then $\ACL(B) \subseteq \DCL(B)$.
 \end{itemize}
\end{lemma}

\begin{proof}
 Let $A = \ACL(A) = \{a_i\mid i < \kappa\}$ and \(c\in\dom(\nM)\).  For each $i$, by hypothesis, we find an \(\aut(\nM/\ACL(a_i))\)-invariant type \(\nP_i\) and an \(\emptyset\)-definable map \(f_i\) such that \(\nP_i(x_i)\vdash f_i(x_i) = a_i\). Let $A_0=A$, and, recursively, let $A_{i+1} = A_{i} \cup \{b_i\}$, where $b_i \models\tprestr{\nP_{i}}{\ACL(A_ic)}$,
and $A_\lambda = \bigcup_{i< \lambda} A_i$ for limit $\lambda$.

\begin{claim}[por]
$\ACL(Ac) \cap \DCL(A_ic) \subseteq \DCL(Ac)$.
\end{claim}

\begin{proof}
 By induction on $i$.   The limit case is trivial.  To move from $i$ to
$i+1$, let $d\in \ACL(Ac) \cap \DCL(A_{i+1}c)$ and let $\sigma\in \aut(\nM/A_{i}c)$. As $\tpL(b_{i}/\ACL(A_{i}c))$ is invariant under $\sigma$, $d\in\ACL(A_{i}c)$ and $d$ is definable over $A_{i}cb_{i}$, we have $\sigma(d) = d$, i.e., $d$ is definable over $A_{i}c$ and hence $d\in \DCL(Ac)$ by induction.
\end{proof}
Now $A_{\kappa}\subseteq \DCL(A_{\kappa}\cap\dom(\nM))$ and so $\DCL(A_\kappa c) = \DCL(A_\kappa c \cap \dom(\nM))$. By (v$'$) this set contains $\ACL(A_\kappa c)$
and hence $\ACL(Ac)$.  Applying  Claim~\ref{claim:por} with $i=\kappa$, we obtain (v).
\end{proof}

\section{Extensible 1-types in intersections of balls}
\label{sec:D3}

The goal of this section is to establish some results about unary types in henselian fields (specifically, finite extensions of $\Qq_{p}$ and ultraproducts of such fields), which will be useful to prove that Proposition~\ref{prop:EI/UFIcrit} can be applied to these fields.

In this section, we will not be considering valued fields in the geometric language as we need quantifier elimination and not elimination of imaginaries. Let $\mathcal{R}$ be a set of symbols; we will be working in
the countable language $\LL := \{\K,+,\cdot,{^{-1}},\val:\K\to\Valgp, r: \K \to \K_r,\ldots\}_{r \in \RR}$ where the $\K_{r}$ are new sorts, each $r$ is such that $\restr{r}{\inv{\K}}$ is a surjective group homomorphism $\inv{\K} \to \K_r$ that vanishes on $1+\gpow{\mathcal{M}}{\nu}$ for some $\nu=\nu(r) \in \Nn$ and the $\ldots$ refer to additional constants on $\K$ and additional relations on the sorts $\K_r$ and $\Valgp$. Let $T$ be some theory of valued fields in this language that eliminates quantifiers. Assume that $\Valgp$ is definably well-ordered in \(T\) (every nonempty definable subset with a lower bound has a least element).

Finite extensions of the $p$-adics fit in this setting, by Prestel-Roquette \cite[Theorem\,5.6]{PR}, if we take the $r_{n}$ to be the canonical projections $\inv{\K}\to \inv{\K}/\gpow{(\inv{\K})}{n}$. Note that every element of these finite groups is in $\dcleq(\emptyset)$. In the case of  ultraproducts of $p$-adic fields of residue characteristic zero and more generally of henselian valued fields with residue characteristic zero (denoted as $\HFO$), one map $r$ suffices: the canonical projection $\rv : \inv{\K}\to\RV$.

Throughout this section $M$ will be a sufficiently saturated model of $T$ and $\unif\in \K(M)$ a uniformizer. We will write $\overline{r}$ for the (possibly infinite) tuple of all $r\in\RR$ and let \(Q_{\RR}\) be the partial \(\star\)-type of elements that are of the form \((\val(x),\overline{r}(x))\) for some \(x\in\K\). We write $\val(x)\gg\val(y)$ if $\val(x) > \val(y) + m\val(\unif)$ for all $m\in\Nn$. Observe that $\val(x-y) \gg \val(x-z)$ implies $r(x-z)=r(y-z)$ for all $r \in \RR$.  Indeed $(y-z)/(x-z) = 1+(y-x)/(x-z) \in 1+ \gpow{\mathcal{M}}{\nu(r)}$.

\begin{notation}
 If $b \in \BallSet(M)$, $x \in \dom(M)$ and $x \notin b$,  the valuation $\val(x-y)$ takes the same value  for all $y \in b$.  We denote it
$\val(x-b)$.   By $\rad(b) $ we denote the infimum of $\val(y-y')$, $y,y' \in b$.

Moreover for all $r \in \RR$, if $\val(x-b) +\nu(r) \val(\unif) \leq \rad(b)$, then $r(x-y)=r(x-y')$ for all $y,y' \in b$.  We write $r(x-b) = r(x-y)$ in this case.
\end{notation}

\begin{definition}
 Let $f = (f_{i})_{i\in I}$ be a family of $A$-definable functions for some $A\subseteq \eq{M}$.  A partial type $p$ over $A$ is
{\em complete over $A$ relative to $f$}  if the map $\tp(c/A) \mapsto \tp(f(c)/A)$
is injective on the set of complete types over $A$ that extend $p$.
\end{definition}

\begin{remark}[relcomp]
The partial type $p(x)$ is   complete over $A$ relative to $f$ if and only if for every formula $\phi(x)$ over $A$,
 there exists a formula $\theta(u)$ over $A$ such that $p(x) \vdash (\phi(x) \iff \theta(f(x)))$.
 \end{remark}


For the rest of the section we are going to study generic types of intersections of balls.
Let $\bar{b} = \{b_i\mid i\in I \}$ be a descending sequence of balls in $\BallSet(M)$.  Let
$P= \bigcap_{i\in I} b_i$.  Let $P_{\Valgp} = \{\gamma \in \Valgp\mid \forall i \in I\,\gamma > \rad(b_i) \}$.
For any $A$ with $b_i \in \dcleq(A)$, we define the generic type of $P$ over $A\subseteq\eq{M}$ to be
\[\tprestr{q_{P}}{A} := P(x) \cup \{x \notin b\mid b \in \BallSet(\acleq(A)), b \mathrm{\,strictly\,included\,in}\,P\}.\]

In Section~\ref{sec:padics}, we will also be considering the $\ACVF$-generic
of such an intersection $P$, i.e., the same notion of genericity but
considered in algebraically closed valued fields. Note that if $L$ is a valued
field, $A\subseteq L$ and $P$ is an intersection of balls in $\BallSet(A)$,
then the difference between the generic type of $P$ over $A$ in $L$ and in $\alg{L}$ is that the latter must also avoid balls that do not have a center or a radius in $L$ but in $\alg{L}$.

\begin{remark}
If $P$ is a strict intersection, i.e., $P$ is not equal to a ball or
equivalently $\bar b$ does not have a minimal element, then for an element to be generic in $P$ over $A$ it suffices to check that $x$ is not contained in any ball $b\in\BallSet(\dcleq(A))$ contained in $P$. Indeed, if $b\in\BallSet(\acleq(A))$, then the smallest ball containing all $A$-conjugates of $b$ is strictly included in $P$ and is definable over $A$.
\end{remark}

In what follows, we will consider $A\subseteq \eq{M}$ containing all constants in $\K$, and $\bar b$ a decreasing sequence of balls in $\BallSet(\dcleq(A))$ (indexed by some ordinal). Unless otherwise mentioned, until Proposition~\ref{prop:D3.1} we will suppose that $P = \bigcap_{i} b_i$ is strict.

\begin{lemma}[D3.0]
Suppose $A \subseteq \K(M)$. Fix $a \in A$ with $a \in b_i$ for each $i$.  Then
$\tprestr{q_{P}}{A} $ is complete relative to the pair of functions $(\val(x-a),\overline{r}(x-a))$.  Moreover, if $P(\dcleq(A)) = \emptyset$ then $\tprestr{q_{P}}{A} $ is complete.
\end{lemma}

\begin{proof}
Taking into account quantifier elimination,  we must show the following:
let $c,c' \in M$ be two realizations of $q:=\tprestr{{q_{P}}}{A}$ such that
$(\val(c-a),\overline{r}(c-a))$ has the same type over $A$ as
$(\val(c'-a),\overline{r}(c'-a))$; then the substructures $A(c)$, $A(c')$ generated by $c$, $c'$ over $A$ (which are simply the fields generated by $c$, $c'$ over $A$) are isomorphic over $A$.

 Extend the valuation  from $\K(M)$ to $L := \alg{\K(M)}$---the algebraic closure of $\K(M)$---and extend each $r \in \RR $ to a group homomorphism with kernel
$\ker(r) \cdot (1+ \unif^{\nu(r)} \Val(L)) \subseteq L $.  This is possible, as  for
all $a \in \ker(r)$ and $b \in \Val(L)$, if $a(1+\unif^{\nu(r)}b) \in \K(M)$, then
$a(1+\unif^{\nu(r)}b) \in \ker(r)$: indeed either $a=0$ or $1+\unif^{\nu(r)}b \in
\K(M)$, and in the latter case, $b \in \K(M)$ and $\val(b) \geq 1$, so $b \in
\Val(M)$ and thus $(1+\unif^{\nu(r)}b) \in \ker(r)$. By construction, the following still holds: for all $x,y,z\in L$, $\val(x-y) \gg \val(x-z)$ implies
$r(x-z)=r(y-z)$.

Then it suffices to show that
$\alg{A}(c),\alg{A}(c')$ are $\alg{A}$-isomorphic, by an isomorphism commuting with the
extensions of the maps $r$ (one can then restrict the isomorphism to $A(c)$). As $(\val(c-a),\overline{r}(c-a))$ and $(\val(c'-a),\overline{r}(c'-a))$ realize the same type over $A$, by taking a conjugate of $c'$ over $A$, we may assume the tuples are equal.

Take any $d \in \alg{A}$. If $d \notin b_i$ for some $i$, then $\val(c-d) = \val(c'-d)$.  Moreover, for any $k\in\Nn$, $\val(c-c') \geq \rad(b_{i+k}) \geq \rad(b_i) + k\val(\unif) > \val(c-d) +k\val(\unif)$; and it follows that $\overline{r}(c-d)=\overline{r}(c'-d)$.

If $d \in b_i$ for each $i$, then the smallest ball $b\in\BallSet(L)$ containing
$a$ and all the conjugates of $d$ over $A$ is (quantifier-free) $A$-definable
in $L$. As $\Valgp$ is definably well-ordered, the $\K(M)$-points of $b$
form a ball $b'\in \BallSet(\dcleq(A))$ which is included in $P$. Hence
$c$ and $c'$ are not in $b'$ nor, in fact, in any of the balls centered around $b'$ with radius $\rad(b') - k\val(\unif)$, for \(k\in\Nn\). It follows that $\val(c - d)=\val(c - a)$ and $\overline{r}(c - d)=\overline{r}(c - a)$, and similarly for $c'$.  Hence $\val(c'-d)=\val(c-d)$ and $\overline{r}(c-d) = \overline{r}(c'-d)$.

As any rational function $g$ over $A$ is a ratio of products of constant or linear polynomials, it follows that
 $\val(g(c))=\val(g(c'))$ and $\overline{r}(g(c))=\overline{r}(g(c'))$. This proves the first part of the lemma.

If $P$ does not contain any point in $A$, then there cannot be any $d\in \alg{A}$ such that $d\in b_{i}$ for each $i$. Indeed, let $d_{j\leq n}$ be the $\LL$-conjugates of $d$ over $A$; then $e := 1/n\sum_{j}d_{j} \in \dcleq(A)$ and for all $i$, $d_{j}\in b_{i+k}$ where $k$ is such that $k\val(\unif) \geq \val(n)$ and $\val(e - d) \geq \val(1/n) + \rad(b_{i+k}) \geq \rad(b_{i})$. It follows that $e\in P(\dcleq(A))$, a contradiction. But the hypothesis about $(\val(x-a),r(x-a))$ is only used in the case $d\in P$. Thus the second assertion follows.
\end{proof}

\begin{remark}[D3.0 HF0]
Suppose $T$ extends $\HFO$ and $A\subseteq \K(M)$.  Without any assumption on
$P$ (it can be strict, a closed ball or an open ball), if $P(A) = \emptyset$ then $P$ is a complete type. The exact same proof works as balls are convex in residue characteristic zero and the unique $r = \rv$ we need has kernel $1+\mathcal{M}$, i.e., $\val(x-y) > \val(x-z)$ alone implies $r(x-z)=r(y-z)$.
\end{remark}

We now want to prove (in Proposition~\ref{prop:D3.1}) that Lemma~\ref{lem:D3.0} is true without the assumption that $A\subseteq \K(M)$.

\begin{lemma}[D3.1.1]
Suppose $A\subseteq \eq{M}$ is such that $P$ contains no $b\in\BallSet(\dcleq(A))$.  Then $\tprestr{q_P}{A}$ is a complete type.
\end{lemma}

\begin{proof}
Let us suppose $A$ is countable. Then the partial type $P =
\bigcap_{n=1,2,\ldots} b_n$ is not isolated over $A$; for if the formula
$\theta(x)$ with parameters in $A$ implies $x\in b_{i}$ for all $i$, then,
as $\Valgp$ is definably well-ordered, there is a smallest ball $b$
containing $\theta$. This ball is strictly contained in $P$ and is $A$-definable,
a contradiction. Then by the omitting types theorem, there exists a model
$M_{0}$ such that $A\subseteq \eq{M_{0}}$ and $P(M_{0}) = \emptyset$. By
Lemma~\ref{lem:D3.0}, $\tprestr{q_{P}}{\K(M_{0})}$ is a complete type, and, as $\K$ is
dominant in $\eq{M_{0}}$, $P$ is a complete type over $\eq{M_{0}}$ and hence over $A$.

If $A$ is not countable, let $c$ and $c'$ be generic in $P$ over $A$ and let $(\eq{M_{0}},A_{0})\prec(\eq{M},A)$ be countable (in the language where we add a predicate for $A$) and contain $c$ and $c'$. Let $Q$ be the intersection of all $A_{0}$-definable balls in $M_{0}$ that contain $c$; then $Q$ is strict, it contains no $A_{0}$-definable ball and also contains $c'$ (all of this is expressed in the type of $c,c'$ in the language with the new predicate). By the countable case, $c$ and $c'$ have the same type over $A_{0}$ in $\eq{M_{0}}$ and hence, they have the same type over $A$ in $\eq{M}$.
\end{proof}

\begin{lemma}[D3.2]
Let $q_{\RR}$ be a complete type over $A$ extending $Q_{\RR}$.  Suppose $q_{\RR}$ implies both that $u\in P_{\Valgp}$ and that, for any $\gamma \in P_{\Valgp}(\dcleq(A))$, $\gamma > u$. Then
$$\tprestr{q_{P}}{M} \cup \bigcup_{a \in P(M)}  q_{\RR}(\val(x-a),\overline{r}(x-a))$$ is consistent.
\end{lemma}

\begin{proof}
 We may assume $M$ has an element $a'$ with $a' \in b_i$
for each $i$.   Note that $q_{\RR}$ is consistent with
 $\{ \gamma > u\mid \gamma \in P_{\Valgp}(M) \}$. Indeed, for any  $\gamma \in  P_{\Valgp}(M)$, if $q_{\RR} \vdash u \geq \gamma$,
then some $\psi \in q_{\RR}$ is bounded below by $\gamma$; but then
 the minimum $\gamma' \geq \gamma$  of $\psi$ in $M$ exists as $\Valgp$ is
 definably well-ordered, $\gamma'$ is in $P_{\Valgp}(\dcleq(A))$ and
$q_{\RR} \vdash \gamma' \leq u$, contradicting our hypothesis.

Let $c'$  be such that  $(\val(c'),\overline{r}(c')) \models q_{\RR} \cup
\{ \gamma > u\mid \,  \gamma \in P_{\Valgp}(M) \}$ and let $d = a'+c'$.
Clearly $d\models \tprestr{q_{P}}{M}$; indeed $\val(d-a') = \val(c')\in P_{\Valgp}$
and thus $d\in b_{i}$ for all $i$. Now, let us assume there exists $b\in\BallSet(\dcleq(M))$ included in $P$ and containing $d$. Taking a bigger ball, we can suppose that $a'\in b$, too; but then $\val(d-a') = \val(c') > \rad(b)-\val(\unif) \in
P_{\Valgp}(M)$ contradicting the choice of $c'$.  Moreover for any $a \in P(M)$,
 $\val(d-a') = \val(c') \ll \val(a-a')$. Thus $\val(d-a) = \val(d-a') = \val(c')$ and
 $\overline{r}(d-a) = \overline{r}(d-a') = \overline{r}(c')$, and $d$ realizes the given type.
\end{proof}

\begin{proposition}[D3.1]
Assume $P$ is strict and fix $a \in \BallSet(\dcleq(A))$ with $a \subseteq b_i$ for each $i$.  Then
$\tprestr{q_{P}}{A} $ is complete relative to the pair $(\val(x-a),\overline{r}(x-a))$.
Moreover, if $P$ does not contain any ball in $\BallSet(\dcleq(A))$ then $\tprestr{q_{P}}{A} $ is complete.
\end{proposition}

\begin{proof}
The second case is tackled in Lemma~\ref{lem:D3.1.1}. So we can suppose that such an $a\in \BallSet(\dcleq(A))$ exists. Let $c,c'\models \tprestr{q_{P}}{A}$ such that $q_{\RR} := \tp(\val(c-a),\overline{r}(c-a)/A) = \tp(\val(c'-a),\overline{r}(c'-a)/A)$. Let $M_{0}\prec M$ be such that $A\subseteq \eq{M_{0}}$. It follows from Lemma~\ref{lem:D3.2} that there exists $c_{0}\models \tprestr{q_{P}}{\eq{M_{0}}}\cup q_{\RR}(\val(x-a),\overline{r}(x-a))$. Taking conjugates of $c$ and $c'$ over $A$, we can suppose that $(\val(c-a),\overline{r}(c-a)) = (\val(c_{0}-a),\overline{r}(c_{0}-a)) = (\val(c'-a),\overline{r}(c'-a))$ as these three tuples have the same type over $A$.

By choice of \(c\), \(c\models \tprestr{q_P}{A}\) and hence \(c\in P\). Moreover, let \(b\in\BallSet(\dcleq(M_0))\); taking a bigger ball if necessary, we may assume that \(a\in b\) and hence that \(\val(c-a) = \val(c_0-a)  <\rad(b)\). So \(c\) is not in \(b\). It follows that \(c\models\tprestr{q_P}{\eq{M_0}}\) and, similarly, \(c'\models\tprestr{q_P}{\eq{M_0}}\). By Lemma~\ref{lem:D3.0}, $c$ and $c'$ have the same type over $\eq{M_{0}}$ and
hence over $A$.
\end{proof}

\begin{corollary}[D3.3]
Let $L$ be a finite extension of $\Qq_{p}$, $M \models \Th(L)$ and $A \subseteq M$ such that $\BallSet(\acleq(A)) \subseteq A$.
Let $c \in \dom(M)$.  Then $\tp(c/A)$  extends to a complete $\aut(M/A)$-invariant type over \(M\).
\end{corollary}

\begin{proof}
Let  $W(c;A) = \{b\in \BallSet(A)\mid  c \in b\}$,     $ P = \bigcap_{b\in W(c;A)} b$.    As the residue
field of $M$ is finite, $P$ cannot reduce to a single ball (that ball would be the union of
finitely many proper subballs, each in $\BallSet(\acleq(A))$, hence in $A$ and
one of them would contain $c$). Note that $c\models \tprestr{q_{P}}{A}$.

If there is no ball $a\in\BallSet(\dcleq(A))$ contained in $P$, then let $q_{\RR}$ be any $\aut(M/A)$-invariant type extending $Q_{\RR}$ that implies $u\in P_{\Valgp}$ and $\alpha > u$  for all $\alpha \in P_{\Valgp} (M)$. If such a ball $a$ exists, we suppose $q_{\RR}$ also extends $\tp(v(c-a),\overline{r}(c-a)/A)$. By Lemma~\ref{lem:D3.2}, $q^* := \tprestr{q_P}{M(x)} \cup \bigcup_{a \in P(M)} q_{\RR}(\val(x-a),\overline{r}(x - a))$ is consistent. Clearly $q^*$ is $\aut(M/A)$-invariant. It follows from Proposition~\ref{prop:D3.1} that $q^*$ is complete and that it extends $\tp(c/A)$.
\end{proof}

Let $N_n$ be the group of matrices of the form $I_n + b$, where $I_n$ is the identity
matrix in $\GL{n}$, and $b$ is an upper triangular matrix with all entries
having valuation $\gg 0$. Thus $N_n = \triang_n(\Val) \cap \bigcap_{m}(I_n + \unif^m \triang_n(\Val))$.

\begin{lemma}[D3.4]
There exists an $\aut(M)$-invariant type $\tprestr{p}{M}$ of matrices $a \in N_n$, invariant under right multiplication: for all $A\subseteq \eq{M}$ and $b\in N_{n}(A)$, if $c \models \tprestr{p}{A}$, then $cb \models \tprestr{p}{A}$. The type $p$ is complete relative to the absolute values and $\overline{r}$-values
of the entries. Moreover, if there exists a complete $\aut(M)$-invariant type \(t(\gamma,x)\) containing \(Q_{\RR}(\gamma,x)\cup\{\gamma>k\val(\unif)\mid k\in\Nn\}\), then $p$ can be taken to be complete.
\end{lemma}

\begin{proof}
Let $P = \bigcap_{i}(\unif^{i}\Val)$ and let $q= \tprestr{q_{P}}{M}$; then $q$ is $\aut(M)$-invariant and complete relative to $\val$ and $\overline{r}$ by Proposition~\ref{prop:D3.1} (as $P$ contains $0$). If $t$ as above exists, then take $q := \tprestr{q_{P}}{M} \cup t(\val(x),\overline{r}(x))$, which is consistent by Lemma~\ref{lem:D3.2}, complete and $\aut(M)$-invariant.

Let $p$ be the type of upper-triangular matrices obtained by taking the $\frac{n(n+1)}{2}$th tensor power of $q$ (where by tensor product, we mean the tensor product of types; see just below for a more explicit statement),
using the lexicographic order on the matrix entries, and adding $1$ on the diagonal:
 thus for all $A\subseteq \eq{M}$, if $a \in M_n$, then $I_n + a \models \tprestr{p}{A}$
if and only if $ a_{11} \models \tprestr{q}{A}, a_{12} \models \tprestr{q}{\dcleq(Aa_{11})}, \ldots,
 a_{22} \models \tprestr{q}{\dcleq(Aa_{11}\ldots a_{1,n})}, \ldots,
 a_{n,n} \models \tprestr{q}{\dcleq(Aa_{11},\ldots a_{n,n-1})} $, while $a_{ij}=0$ for $i>j$.

The fact that $p$ is an $\aut(M)$-invariant (partial) type of elements of $N_n$ is clear.  As for
the right translation invariance, let $I_n+b \in N_n(A)$ and $I_{n}+a\models \tprestr{p}{A}$; we have to show that
$(I_n+a)(I_n+b) = I_n + a+b+ab \models \tprestr{p}{A}$.  Let $d=a+b+ab$.
  Then
$d_{11} = a_{11}+b_{11}+ a_{11}b_{11}$.  We have
$$\val(a_{11}b_{11}) > \val(b_{11}) \gg \val(a_{11}) \gg 0.$$
So $\val(d_{11}) = \val(a_{11})$ and hence $d_{11}$ also realizes $\tprestr{q_{P}}{A}$. Furthermore, we also have $\overline{r}(d_{11})=\overline{r}(a_{11})$; it follows that
$d_{11} \models \tprestr{q}{A}$.  Similarly
$$d_{12} = a_{12} + b_{12} + a_{11}b_{12}+ a_{12}b_{22};$$
here $a_{12}$ has strictly bigger valuation than any of the other summands, so again
$\val(d_{12}) = \val(a_{12})$ and $\overline{r}(d_{12})=\overline{r}(a_{12})$, thus $d_{12} \models \tprestr{q}{\dcleq(Aa_{11})}$.
But since $b \in \dcleq(A)$, $d_{11} \in \dcleq(Aa_{11})$, so
$d_{12} \models \tprestr{q}{\dcleq(Ad_{11})}$.  Continuing in this way we see that $I_{n} + d \models \tprestr{p}{A}$.
\end{proof}

In the following two proofs, whenever \(A\subseteq\eq{M}\), \(\mathcal{G}(A)\) will denote the points of \(A\) that belong to a sort of the language \(\LG\subseteq\eq{\LL}\). Note that, since \(N_n\) is an interesection of quantifier-free definable groups, the elements of \(\triang_n(\K)/N_n\) can be identified with with infinite tuples in \(\mathcal{G}(\eq{M})\).

\begin{corollary}[D3.5]
 Let  $R$ be a left coset of $N_n$ in $\triang_n(\K)$. There exists an $\aut(M/R)$-invariant type of elements of $R$.
\end{corollary}

\begin{proof}
Pick $g \in R$, let $p$ be the right-$N_n$-invariant type of Lemma~\ref{lem:D3.4}, and for all $A\subseteq \eq{M}$, let $\tprestr{p^g}{A} = \tp(cg/Ag)$, where $c \models \tprestr{p}{\dcleq(Ag)}$.  Then $\tprestr{p^g}{A} = \tprestr{p^{hg}}{A}$ for $h \in N_n(\dcleq(A))$,
since $p$ is right-$N_n$-invariant.  Thus any automorphism fixing $R$ must fix the global type $\tprestr{p^g}{M}$.
\end{proof}

\begin{corollary}[D3.6]
 Let $L$ be a finite extension of $\Qq_{p}$ and let $M \models \Th(L)$, $e \in \Latt[n](M)$, $E = \mathcal{G}(\acleq(e))$. Then there exists an $\aut(M/E)$-invariant type of bases for $e$.
\end{corollary}

\begin{proof}
 It was noted in Section~\ref{sec:Vfield} that any lattice $e$ has a triangular basis; this basis can be viewed as the set of columns of a matrix in $\triang_n(\K)$.  Let $b,b'$ be two such bases, and suppose
$b' = \sigma(b)$, $\sigma \in \aut(M/E)$. Then as $e/\unif^me$ is finite and \(e/\unif^me\) is coded in \(\mathcal{G}\), the cosets of \(\unif^me\) in \(e\) are fixed by \(\aut(M/E)\). Thus, the columns of $b,b'$ must be in the same coset of $\unif^me$ for each $m$.  Thus if we write $b' = ab$ with $a \in \triang_n(\Val)$,
then $a=I_n \hbox{ modulo } \unif^m\Val$ for each $m$, so $a \in N_n$ and $\aut(M/E)$ preserves
the coset $R := N_nb$.  So it suffices to take the $\aut(M/R)$-invariant
type of elements of $R$ guaranteed by Corollary~\ref{cor:D3.5}.
\end{proof}

Let us now suppose that $T$ extends $\HFO$. Using similar techniques, we can extend the previous results to the case when $P$ is a closed ball (this case is only relevant to Section~\ref{sec:asymptotic}). For the last result, though, we will also need the residue field to be pseudo-finite.

Let $b$ be a closed ball.  We will write $\resf_{b}$ for the map that sends $x\in b$  to $x+\rad(b)\mathcal{M}$, the maximal open subball of $b$ containing $x$.

\begin{lemma}[D3.2 closed]
Let $b\in\BallSet(\dcleq(A))$ be a closed ball and $q$ a complete type over $A$ containing the formula $x\in\resf_{b}(b)$ such that $q \vdash x\neq b'$ for all $b'\in \resf_{b}(b)(\acleq(A))$.  Then $\tprestr{q_{b}}{M}\cup q(\mathrm{res}_{b}(x))$ is consistent.
\end{lemma}

\begin{proof}
Let us first show that $q$ is consistent with $\{x\neq b' \mid b'\in \resf_{b}(b)(\eq{M})\}$. If not, there is a finite number of balls $b_{i}\in \resf_{b}(b)(\eq{M})$ such that $q\vdash \bigvee_{i} x = b_{i}$. If we take a minimal number of such balls, each of them must realize $q$ and hence be algebraic over $A$, a contradiction.

Now, let $c$ be such that $\resf_{b}(c) \models q\cup \{x\neq b' \mid b'\in
\resf_{b}(b)(\eq{M})\}$; then $c\models \tprestr{q_{b}}{M}$. Indeed $c\in b$ and if $c$
is in $b'\in\BallSet(\eq{M})$ such that $b'\subseteq b$, then $c\in \resf_{b}(b')\in \resf_{b}(b)(\eq{M})$, contradicting the choice of $c$.
\end{proof}

\begin{lemma}[D3.7]
Suppose $P = b$ is a closed ball.  Then $\tprestr{q_{b}}{A}$, the generic type of $b$, is complete relative to $\resf_{b}$.
\end{lemma}

\begin{proof}
If $A \subseteq \K(M)$ then, by the same considerations as in
Lemma~\ref{lem:D3.0} (and, as $\HFO\subseteq T$, taking $\overline{r}=
\rv$ is enough), it suffices to show that if $c$ and $c'$ are realizations
of $\tprestr{q_{b}}{A}$ such that $\resf_{b}(c) = \resf_{b}(c')$ then for all $d\in
\alg{A}$, $\rv(c-d) = \rv(c'-d)$. If $d\in \resf_{b}(c)$, then $c\in
\resf_{b}(c) = \resf_{b}(d) \in\BallSet(\acleq(A))$ as $d\in \alg{A}$.  This
contradicts the fact that $c\models \tprestr{q_{b}}{A}$. Hence $d\not\in \resf_{b}(c)$. As $c$ and $c'\in \resf_{b}(c) = \resf_{b}(c')$, $\val(c-c') > \val(c-d)$ and $\rv(c-d) = \rv(c'-d)$.

If $A$ is not contained in $\K$, let $c,c'\models \tprestr{q_{b}}{A}$ such that $q := \tp(\resf_{b}(c)/A) = \tp(\resf_{b}(c')/A)$. By Lemma~\ref{lem:D3.2 closed}, there exists $c_{0}\models \tprestr{q_{b}}{M}\cup q$. Taking $A$-conjugates of $c$ and $c'$, we can suppose that $\resf_{b}(c) = \resf_{b}(c_{0}) = \resf_{b}(c')$. Then, as seen in the proof of Lemma~\ref{lem:D3.2 closed}, $c,c'\models \tprestr{q_{b}}{M}$. By the previous paragraph $c$ and $c'$ have the same type over $M$ and hence over $A$.
\end{proof}

\begin{corollary}[D3.8]
Suppose $P = b$ is a closed ball and let $a\in \BallSet(\dcleq(A))$ be contained in $b$.  Then $\tprestr{q_{b}}{A}$ is complete relative to $\rv(x-a)$.
\end{corollary}

\begin{proof}
If $c, c' \models \tprestr{q_{b}}{A}$, then $\val(c-a) = \val(c'-a) = \rad(b)$ and hence $\resf_{b}(c) = \resf_{b}(c')$ if and only if $\rv(c-a) = \rv(c'-a)$. Thus the corollary follows immediately from Lemma~\ref{lem:D3.7}.
\end{proof}

\begin{corollary}[D3.9]
Suppose $\res$ is pseudo-finite, $\res(A)$ contains the constants needed for $\res$ to have EI and $P = b$ is a closed ball that contains no ball $a\in\BallSet(\dcleq(A))$.  Then any $x\in b$ generates a complete type over $A$.
\end{corollary}

\begin{proof}
By Lemma~\ref{lem:D3.7}, it suffices to show that $\resf_b(b)$ is a complete
type over $A$. But $\resf_{b}(b)$ is a definable 1-dimensional affine space
over $\res$---i.e., a $V := \gamma\Val/\gamma\mathcal{M}$-torsor where $\gamma := \rad(b)$. Hence $H := \aut(\resf_{b}(b)/\res,A)$ is a subgroup of a semi-direct
product of $V$ and the multiplicative group
$\Gg_{m}(\res)$. The subgroup $H\cap V$ (i.e., the group of translations of
$\resf_{b}(b)$ that also are automorphisms over $A$ and $\res$) is $\infty$-definable
over $A$. Indeed it is the set
$\{u\in V\mid \forall \overline{y}\, \forall x\,(x\in\resf_b(b)\wedge \overline{y}\in \res)\Longrightarrow (\phi(x,\overline{y})\iff
\phi(x+u,\overline{y}))\ \mbox{for all}\ \mbox{$A$-formulas $\phi(x,\overline{y})$}\}.$

Since \(\Gg_a(\Ff_p)\) has no proper nontrivial subgroups and \(\res\), being pseudo-finite, is elementarily equivalent to an ultraproduct of \(\Ff_p\), it follows that \(\Gg_a(\res)\) has no proper nontrivial definable subgroups and hence neither does \(V\). Because in a pseudo-finite field any $\infty$-definable group is an intersection of definable groups, $V$ has no nontrivial proper $\infty$-definable subgroup either. If $H\cap V = V$ then $H$ acts transitively on $\resf_b(b)$ (by translation) and, as $H \leq \aut(\resf_{b}(b)/A)$, we are
done. On the contrary, if $H\cap V=\{1\}$, then $H$ contains no translations and must either have exactly one fixed point or be the trivial group and hence fix all points in $\resf_{b}(b)$.

Suppose $H$ has only one fixed point $a \in \resf_{b}(b)$ and let $\theta\in \aut(\resf_{b}(b)/A)$. For any $\sigma\in H$, $\theta^{-1}\circ\sigma\circ\theta \in H$ and hence $(\theta^{-1}\circ\sigma\circ\theta)(a) = a$, i.e., $\theta(a)$ is fixed by $\sigma$. As $a$ is the only point fixed by $H$, $\theta(a) = a$ and $a\in\dcleq(A)$: but this is a contradiction. It follows that $H$ fixes every point in $\resf_{b}(b)$ and hence, because $\res$ is stably embedded, $\resf_{b}(b) \subseteq \dcleq(\res,A)$. But then we must also have $V \subseteq \dcleq(\res,A)$. Hence $(V,\resf_{b}(b))$ is $A$-definably isomorphic to a definable (regular) homogeneous space $(G,R)$ of $\eq{\res} = \res$. As $\res$ is stably embedded, $(G,R)$ is definable over $A' := \eq{\res}(\dcleq(A)) = \res(\dcleq(A))$.

Hence we only have to show that any $A'$-definable \(\Gg_a(\res)\)-torsor in a pseudo-finite field \(\res\) has an $A'$-point to obtain a contradiction. Let us consider $\res$ elementarily embedded in the fixed field of $L\models ACFA$ and let $\overline{A'}$ be the algebraic closure of $A'$ in $L$. Note that $A$ is algebraically closed in $ACFA$ and is a model of $ACF$. By usual arguments (e.g., \cite{KowPilGpDiff}) there exists
an $ACF$ $\overline{A'}$-definable homogeneous space $(G',S')$ and interalgebraic group configurations in $(G,R)$ and $(G',S')$. Replacing $G'$ with its identity component $G_{0}'$ and $S'$ with the $G_{0}'$-orbit of any $\overline{A'}$-point in $S'$ (there is such a point because $\overline{A'}\models ACF$), we can suppose that $G'$ is connected. By some additional classical arguments (although the literature mainly concerns itself with groups and not homogeneous spaces at this point: see \cite{KowPilGpDiff} again), there is an $\overline{A'}$-definable subgroup $H$ of $G\times G'$, such that $H_0 := \{x\in G \mid (x,0)\in H\}$ and $H_0' := \{x\in G'\mid (0,x)\in H\}$ are finite central subgroups and the left and right projections of $H$ must have finite index in $G$ (respectively $G'$). But as $G$ and $G'$ are connected, these projections must be the groups themselves. As $G$ has no torsion (we are in characteristic $0$), $H_{0}$ is trivial. Taking the quotient of $(G',S')$ by $H_0'$---i.e., considering the group $G'/H_{0}'$ acting on the $H_{0}'$ orbits of $S'$---the group $H$ is in fact (the graph of) an isomorphism. In particular, as $G$ has no proper definable subgroup, this implies that the action of $G'$ on $S'$ is also regular, i.e., $S'$ is a $G'$-torsor.

Let $(a,a')$ be generic in $R\times S'$, let $X$ be the $H$-orbit of $(a,a')$ and let $P = \tp(aa'/\overline{A'})$. As $P$ and $X$ have the same dimension (equal to $1$), $P$ cannot be covered with infinitely many $H$-orbits (pseudo-finite fields have the (E) property of \cite{HruPilGpPFF}) and as $\overline{A'}$ is algebraically closed (including imaginaries), $X$ must contain $P$ and hence is  $\overline{A'}$-definable. Moreover, it is quite easy to see that $X$ is (the graph of) an isomorphism between $R$ and $S'$. As $S'$ contains $\overline{A'}$-points,
so does $R$. Let $d$ be one of these points, and let $(d_{i})_{i=1\ldots n}$ be its $A'$-conjugates. Then $1/n\sum_{i}d_{i}\in R(A')$, and we have the $A'$-point we have been looking for.
\end{proof}

To conclude this section, we summarize the classification of unary types in $\PLO$.

\begin{proposition}[classif types PLO]
Suppose $T$ extends $\PLO$ and $\res(A)$ contains the constants needed for $\res$ to have EI. Let $a \in \BallSet(\dcleq(A))$ with $a \subseteq b_i$ for each $i$.  Then $\tprestr{q_{P}}{A} $ is complete relative to $\val(x-a)$ and to $\overline{r}(x-a)$. Moreover, if $P$ does not contain any ball in $\BallSet(\dcleq(A))$ then $\tprestr{q_{P}}{A} $ is complete.
\end{proposition}

\begin{proof}
If $P$ is strict we can apply Proposition~\ref{prop:D3.1}. If not, we apply Corollary~\ref{cor:D3.8} or Corollary~\ref{cor:D3.9}.
\end{proof}

\section{The \texorpdfstring{$p$}{p}-adic case}
\label{sec:padics}

Let $L$ be a finite extension of $\Qq_{p}$. As stated in Remark~\ref{constants Qp}, it can be shown that there exists a
number field $F \subseteq L$ that contains a uniformizer $\unif$ of $L$, such
that $\resf(L) = \resf(F)$ and such that every finite extension $L'$ of $L$
is generated by an element $\alpha$ whose minimal polynomial is defined over $F$
and such that $\alpha$ also generates the valuation ring $\Val(L')$ over
$\Val(L)$. Let $T_{L}$ denote the theory of $L$ in $\LG\cup\{P_{n} :
n\in\Nn_{>0}\}\cup\{c\}$, where the predicates $P_{n}$ stand for the nonzero $n$th powers (in the sort $\K$) and $c$ generates $F$ over $\Qq$. Then $T_{L}$
is model complete (cf.\ \cite[Theorem\,5.1 and Corollary\,5.3]{PR}) and it is axiomatized by the fact that $\K$ is a henselian valued field with value group a $\Zz$-group and residue field $\Ff_{p}$, by the isomorphism type of $F$ and by the definition of the $P_{n}$ predicates.

We now check the hypotheses of Corollary~\ref{cor:EIcrit}  for $\nT=T_{L}$ and $\tT=\ACVF_{0,p,F}^{\mathcal{G}}$ (the theory of algebraically closed valued fields of mixed characteristic in the geometric language with a constant for $c$; the $F$ in the subscript is there to recall that we added a constant for a generator of $F$ to the theory). We use the same notation as in Proposition~\ref{prop:EI/UFIcrit}.

\paragraph{(i)  Relative algebraic boundedness:}  By model completeness  and the nature of the axioms---the only axioms that are not universal are the fact that the field is henselian and the definition of the predicates $P_{n}$; but both state the existence of algebraic points---$\acll(\nM'c) \cap \nM$ is an elementary submodel of $\nM$, hence certainly $\nL$-definably closed.

\paragraph{(ii)   Internalizing $\tL$-codes:}
As $\K(\nM)$ is henselian, $\K(\dcll(\nM)) = \K(\nM)$, hence if $\epsilon\in \K$, there is
nothing to do. For any element $\epsilon$ of $\Latt[n](\tM)$ let us write $\latt(\epsilon) \subseteq \cpow{\K}{n}$ for the lattice represented by $\epsilon$. If $\epsilon\in \Latt[n](\dcll(\nM))$, $\latt(\epsilon)$ has a basis in some finite extension $L_{0}$ of $L:=\K(\nM)$.   Say $[L_0:L]=m_0$; let $L'$ be the join of all field extensions of $L$ of degree $m_0$.  Then $L'$ is a finite extension of $L$ such that any $\sigma\in \aut(\tM)$ stabilizing $\nM$, stabilizes $L'$; let $[L':L]=m$. By hypothesis, there is a generator $a$ of $L'$ over $L$ whose characteristic polynomial over $L$ is defined over $F$. One has an $a$-definable isomorphism $f_a: L' \to \cpow{L}{m}$ (as vector spaces over $L$), with $f_a(\Val(L')) = \cpow{\Val(L)}{m}$  (i.e., $\Val(L')$ is a free $\Val(L)$-module of rank $m$).   The morphism $f_a$ further induces an isomorphism of the lattice $\latt(\epsilon)(L')$ with a lattice $f_a(\latt(\epsilon)(L')) = \latt(\eta)(L)$ for some $\eta\in \Latt[nm](\nM)$.  As any $a'$ of the (finitely many) that are $\aut(L'/F)$-conjugate to $a$ is also $\aut(L'/L)$-conjugate to $a$, we see that $\latt(\eta)(L) = f_{a'}(\latt(\eta)(L'))$ as well.  Thus $\epsilon$ and $\eta$ are interdefinable in the sense required in (ii).

The argument for $\Tor[n]$ is similar
(alternatively, for finite extensions $L'$ of $L$, the value group also
has a least element,  so we can apply Remark~\ref{no Tn}).

\begin{remark}[3.2] We have proved something slightly stronger than (ii): we also have $\epsilon \in \dcll(\eta)$. The inverse of $f_a$ is a linear map $\cpow{L}{m} \to L'$, say $g_a(\alpha_1,\ldots,\alpha_m) := \sum \alpha_i a^i$.   From the viewpoint of $\tM$, $g_a$ is an $a$-definable linear map with $g_a(\latt(\eta))=\latt(\epsilon)$ (as $g_a$ is $\K$-linear, this remains true for the lattices generated by the $L$- or $L'$- points of $\latt(\epsilon)$ and $\latt(\eta)$). Moreover this is also true for any of the finitely many conjugates of $a$.  Thus $\epsilon \in \dcll(\eta)$.  \end{remark}

The following corollary of this stronger version of (ii) is not needed for what follows but it does shed some light on the interaction between automorphisms of $\nM$ and $\tL$-definable sets.

\begin{corollary}[dcl2] Let $A =\DCL(A) \subseteq \nM$. Let $G$ be the group of automorphisms of $\tM$ that stabilize $\nM$ and fix $A$ point-wise. Let $\epsilon \in \tM$, and assume $g(\epsilon)=\epsilon$ for all $g \in G$.  Then $\epsilon \in \dcll(A)$.
\end{corollary}

 \begin{proof}  We have $\epsilon \in \dcll(\nM)$, since $\aut(\tM/\nM)$ fixes $\epsilon$.  Let $\eta$ be as in (ii). Then \(G\) fixes \(\eta\). Recall that we assumed that any automorphism of \(\nM\) extends to an automorphism of \(\tM\), i.e., $G$ maps surjectively to $\aut(\nM/A)$. So we have $\eta \in A$.  By Remark~\ref{rem:3.2}, $\epsilon \in \dcll(\eta)$.  So $\epsilon \in \dcll(A)$.
 \end{proof}

\paragraph{(iii) Unary EI:}
In \cite{ScowEIQp} P. Scowcroft proved a weak version of this, where the sets are classes of equivalence relations in two variables.
We prove here that every unary subset can be coded in $\BallSet$.

Let $e$ be an imaginary code for a unary subset $D \subseteq \K(\nM)$.  Let $A = \ACLeq(e)$ and let $B = \BallSet(A)$.

\begin{claim}[typeball]
 For all $c \in \K(\nM)$, $\tpL(c/B) \vdash \tpL(c/A)$.
\end{claim}

\begin{proof} Following the notation of Corollary~\ref{cor:D3.3}, recall
  that $W(c;A) = \{b \in \BallSet(A) \mid  c \in b\}$. Let $P = \bigcap W(c;A) =
  \bigcap W(c;B)$, a strict intersection. Then $\tpL(c/B) \vdash \tprestr{q_{P}}{B}
  = \tprestr{q_{P}}{A}$. By Proposition~\ref{prop:D3.1}, either $\tprestr{q_{P}}{A}$ is a
  complete type and we are done, or there is some $a\in B$ such that
  $a\subseteq P$ and $\tprestr{q_{P}}{A}$ is complete relative to $\overline{r}(x-a)$
  and $\val(x-a)$. As $\inv{\K} = F\gpow{(\inv{\K})}{n}$ for all $n$, \(r_n(c-a) \in r_n(F)\) and hence $\tpL(\overline{r}(c-a)/A)$ follows from its type over $F$, i.e., over $\DCL(\emptyset)$. Moreover   $\Valgp(\DCL(B)) = \Valgp(\DCL(A))$ (as elements of  $\Valgp$ are coded by balls). Thus, as $\Valgp$ is stably embedded and has unary EI, $\tpL(\val(c-a)/B) \vdash \tpL(\val(c-a)/A)$ and we have the expected result.
\end{proof}

As $D$ is $\nL(A)$-definable, $D$ is also $\aut(\nM/B)$-invariant, so that by
compactness $D$ is definable over $B$.  Hence  $e \in \DCLeq(B)$. We
conclude as in Corollary~\ref{cor:EIcrit}:   there is a tuple $a$ from
$\BallSet$ with $a \in \ACL(e)$ and $e \in \DCLeq(a)$; so $\DCLeq(e) =
\DCLeq(a')$, where $a'$ is the finite set of $\nL(e)$-conjugates of $a$.
We already know that finite sets are coded (e.g., by (ii) and Lemma~\ref{lem:EFin}).

 \paragraph{(iv)  Invariant types and germs:} The main ingredient for this proof is the C-minimality of $\ACVF$, i.e., the fact that every definable subset of \(L\models\ACVF\) is a finite Boolean combination of balls (and points).


Let $A = \ACL(A)$, $c\in \K(\nM)$, $W(c;A) = \{b_{i} \mid i \in I\}$ and $P = \bigcap_i b_i$. The balls $b_{i}$ are linearly ordered by inclusion, and we order $I$ correspondingly: $i\leq j$ holds if $b_{j}\subseteq b_{i}$. As seen previously, $P$ is a strict intersection. Let $\tP$ be the $\ACVF$ generic of $P$.

If $r(x,b)$ is an $\tL$-definable function,
 let $X(b',b'') = \{x  \mid r(x,b') \neq r(x,b'') \}$.   Then $X(b',b'')$ is
 a finite Boolean combination of balls and there exists $i=i(b',b'')$ such that $X(b',b'') \cap P$ is contained in a proper subball of $P$ if and only if for each $j \geq i$, $X(b',b'') \cap b_{j}$ is contained in a proper subball of $b_j$.

Define an equivalence relation $E_i$ by $b'E_{i} b''$ if and only if  $X(b',b'') \cap b_i$ is contained in a proper subball  of $b_i$ (i.e., $r(x,b')$ and $r(x,b'')$ have the same germ on the $\ACVF$-generic of $b_{i}$).  Let $e_{i} = b/E_{i}$.  Then:
\[\begin{array}{p{10pt}cl}
\multicolumn{3}{l}{\sigma\in \aut(\tM/A)\text{ fixes the \(\tP\)-germ of }r(x,b)}\\
&\iff&r(x,b)\text{ and }r(x,\sigma b)\text{ have the same }\tP\text{-germ} \\
&\iff&X(b ,\sigma b) \cap P\text{ is contained in a proper subball of }P \\
&\iff&\text{ for some }i\text{, for all }j \geq i,\,b E_j \sigma(b) \\
&\iff&\text{ for some }i\text{, for all }j \geq i,\,\sigma\text{ fixes }e_{j}.
\end{array}\]

As for the consistency of $\tprestr{\tP}{\nM}$ with $\tpL(c/A)$: by definition of the $\ACVF$ generic, $\tprestr{\tP}{\nM}$ is generated by $P$ along with all formulas $x \not\in b$, where $b\in\BallSet(\dcll(\nM))$ is a proper subball of $P$. As $P$ is  part of $\tpL(c/A)$, it suffices to show that $\tpL(c/A)$ does not imply any formula $x \in d$ with $d$ a finite union of balls $d_{j}\in\BallSet(\dcll(\nM))$ strictly included in $P$.

\begin{claim}[balls from tM]
For all $b\in\BallSet(\dcll(\nM))$ such that $b(\nM) = \{x\in b\mid x\in M\} \neq\emptyset$, there exists $b'\in\BallSet(\nM)$ such that $b(\nM) = b'(\nM)$.
\end{claim}

\begin{proof}
As $\Valgp$ is definably well-ordered, $\inf\{\val(a-c)\mid a,c\in b(\nM)\} = \gamma\in\Valgp(\nM)$. We can now take $b'$ to be the ball of radius $\gamma$ around any point in $b(\nM)$.
\end{proof}

If \(\tpL(c/A)\) implies \(x\in d\) for \(d\) as above, then it follows from the claim that $d(\nM)$ is equal to a finite union $d'$ of balls in $\BallSet(\nM)$ and $\tpL(c/A)$ implies $x\in d'\subseteq P$. But this would contradict Lemma~\ref{lem:D3.2}.

\paragraph{(v) Weak rigidity:}

We use Lemma~\ref{lem:por}.  The hypothesis that for all $a\in \nM$ there is
a tuple $c\in \K(\nM)$ such that $a\in\DCL(c)$ and $\tpL(c/\ACL(a))$ extends to
an $\aut(\nM/\ACL(a))$-invariant type, holds trivially when $a\in \K(\nM)$ and follows from Corollary~\ref{cor:D3.6} when $a\in \Latt_m(\nM)$. If $a\in T_m(\nM)$ for some $m$ then, as the value group has a least element, $a$ is coded by an element of $S_{m+1}(\nM)$ (see Remark~\ref{no Tn}) and hence, applying Corollary~\ref{cor:D3.6} to the code in $S_{m+1}(\nM)$, we are done.

The  assumption (v$'$) of Lemma~\ref{lem:por} is proved for $\Qq_{p}$ by van den Dries in \cite{vdDSko}. Let us briefly recall his proof to check that it adapts to the finite extension of $\Qq_{p}$ case.

Let $B \subseteq \K(\nM)$ (we can assume that $B =
\dcll(B)\cap \K(\nM)$ is a field and contains $F$).  Let $\sigma \in \aut(\nM/B)$  and let $B' =
\fix(\sigma)\cap(\acll(B)\cap \nM)$. It suffices to show that $B'\models T_{L}$. Indeed, by model completeness, $B'\prec \nM$ will then contain $\ACL(B)$, hence $\ACL(B)$ is rigid over $B$.

As noted in the proof of (i), in order to show that $B'\models T_{L}$, we only have to show that $B'$ is henselian and that the definition of the $P_{n}$ is preserved.

By the universal property of the henselization, $B^{h}$ is contained in
$B'$ and thus $B'$ is henselian. Moreover, let $x\in B'\cap P_{n}(\nM)$ and
let $y\in \K(\nM)$ such that $x = y^{n}$. Note first that $(y/\sigma(y))^{n} =
x/\sigma(x) = 1$ and thus that
$y/\sigma(y)\in\ACL(\emptyset)$. Furthermore, for all $m\in\Nn$, there
exists $q\in F$ such that $yq\in P_{m}(\nM)$. But then $y/\sigma(y) = yq/\sigma(yq) \in P_{m}(\nM)$ for all $m$. As $\bigcap_{m}P_{m}(\ACL(\emptyset)) = \{1\}$, it follows that $y = \sigma(y)$, i.e., $y \in B'$.

\begin{remark}
As in \cite{vdDSko}, it follows from this proof that the restriction of $T_{L}$ to the sort $\K$ has definable Skolem functions.
\end{remark}

\begin{proof}[of Theorem~\ref{thm:Qp}]
 By Corollary~\ref{cor:EIcrit}, we have EI to the sorts $\K,\Latt[n],\Tor[n]$.  But as is explained in Remark~\ref{no Tn}, the   sorts $\Tor[n]$  are not actually needed.
\end{proof}

We finish the section with some additional remarks.

\begin{remark}[without constant]
If we do not want to add a constant $c$ to the language, then it suffices to add
``Galois-twisted $\Latt[n]$'', interpreted as $\Latt[n](K')$ for $K'$ ranging over the
finite extensions of $\K(M)$.

Indeed, by Theorem~\ref{thm:Qp}, any imaginary $e$ is interdefinable over $c$ with some tuple of real elements $e'$. So we have an $e$-definable function $f_e$ with
$f_e(c)=e'$ and a $\emptyset$-definable function $h$ with $h(c,f_{e}(c)) = e$. As $c$
is algebraic over $\Qq$, restricting to $e$-conjugates of $c$, we can take
the graph of $f_{e}$ (a finite set) to be a complete type over $e$.

With the new sorts, it is clear that (ii) holds without adding a
constant and $f_{e}$ is coded by some tuple $d\in \nM$. Let us now show that $d$ is a code for $e$. If $e'$ is
$\nL(d)$-conjugate to $e$ there is some $\sigma\in\aut(\nM/d)$ such that $\sigma(e) = e'$. As $\sigma$ fixes $d$, $c' := \sigma(c)$ is also in the domain of $f_{e}$ and hence $\tpL(c'/e) = \tpL(c/e)$, i.e., $e' = \sigma(e) = \sigma(h(c,f_e(c))) = h(\sigma(c),f_{e}(\sigma(c))) = e$. This implies that $d$ is a code for $e$.
\end{remark}

\begin{remark}[inv ext]
 Let $A = \ACL(A) \subseteq \nM \models T_{L}$.  Then every type
over $A$ extends to an $\aut(\nM/A)$-invariant type.
\end{remark}

This follows immediately from \cite[Prop.\ 2.13]{HruPilNIP} and
Corollary~\ref{cor:D3.3}. But, since the more subtle considerations of
\emph{op.\ cit.}\ are not necessary in \(T_L\) as, in the relevant
case, the algebraic closure coincides with the definable closure, let
us give a more straightforward proof:

\begin{proof}
  Let $c \in \nM$; then $c = f(a_1,\ldots,a_n)$, where $a_i \in
  \dom(\nM)$, and $f$ is $\emptyset$-definable.  It suffices to extend
  $\tpL(a_1,\ldots,a_n/A)$ to an $\aut(\nM/A)$-invariant type.  If
  $\tpL(c/\nM)$ and $\tpL(d/\nM c)$ are $\aut(\nM/A)$-invariant, then
  so is $\tpL(cd/\nM)$; so it suffices to show that $\tpL(a_i/A_i) $
  extends to an $\aut(\nM/A_i)$-invariant type for each $i$, where
  $A_i := \DCL(A_{i-1}a_{i-1})$.  But (by hypothesis (v) of
  Corollary~\ref{cor:EIcrit}) we have that $A_i = \ACL(A_i)$, so
  Corollary~\ref{cor:D3.3} applies.
\end{proof}

\begin{remark}[4.2]
Rigidity of finite sets fails for the theory of a finite extension of the $p$-adics in the geometric language, i.e., $\ACL\neq \DCL$.
\end{remark}

\begin{proof}
  Note first that the angular component maps factor through the projection
  to $\K/\gpow{\K}{p}$ and hence an angular component map is just
  defined by a map between finite sets whose points are all in
  $\DCL(\emptyset)$. It follows that $L$ admits an
  $\emptyset$-definable angular component map $\ang$.
  
  As the value group is stably embedded, one can find a nontrivial
  automorphism $\sigma$ fixing the value group in a sufficiently
  saturated model. By definability of $\ang$, and since \(\sigma\)
  fixes the residue field, it follows that $x$ and $\sigma(x)$ have
  the same angular component.  Take $a\in\Val$ with $\sigma(a) \neq
  a$.  Let $\gamma = \val(\sigma(a)-a)$, $\ang(\sigma(a)-a) =:
  \alpha$.  Then $\val(\sigma^2(a)-\sigma(a)) = \gamma$,
  $\ang(\sigma^2(a) - \sigma(a)) = \alpha$, etc.  As $p \cdot \alpha =
  0$ in the residue field, $(\sigma^p(a)-a) = \sum_{i=0}^{p-1}
  (\sigma^{i+1}(a) - \sigma^i(a)) $ has valuation $\delta > \gamma$.
  Thus in the ring $\Val / \delta\Val$, the image of $a$ is not a
  fixed point, but has an orbit of size $p$ under $\sigma$.  This set
  of size $p$ is not rigid.
\end{proof}

\begin{remark}
The same techniques developed here to prove elimination of imaginaries in $\Qq_{p}$ can also be used to give an alternative proof for elimination of imaginaries in real closed valued fields (see \cite{MelRCVF}). Hypothesis (i) of Corollary~\ref{cor:EIcrit} also follows from the fact that the algebraic closure is a model, (ii) follows as in the $p$-adic case, (iii) follows from the description of 1-types given in \cite[Proposition\,4.8]{MelRCVF}; and so does the existence of the type in (iv). The rest of (iv) is proved exactly as here and so is (v).
\end{remark}

\section{The asymptotic case}
\label{sec:asymptotic}

Recall that $\HFO$ denotes the theory of henselian fields of residue
characteristic $0$ and $\PLO$ is the theory of henselian fields with
value group a $\Zz$-group and residue field a pseudo-finite field of
characteristic $0$. Our goal is now to prove that any completion
$T_{F}$ of $\PLO$ in the language $\LG$ with constants added for some
subfield $F$ of the field sort $\K$ (see Remark~\ref{constants asymp})
eliminates imaginaries. We will be using
Proposition~\ref{prop:EI/UFIcrit} with $\nT = T_{F}$ and
$\tT=\ACVF_{0,0,F}$. We still follow the notation of this proposition.

It is worth noting that we will not, in general, be able to use Corollary~\ref{cor:EIcrit} as there are some ultraproducts of $p$-adics where (v) is false. Indeed, it is shown in \cite[Theorem\,7]{BeyHru} that there exist a characteristic zero pseudo-finite field $L$, $A\subseteq L$, and $b\in L$
such that $b$ has a finite nontrivial orbit over $A$. Then $A$ can be identified with the set $A': = \{at^{0} : a\in A\}\subseteq L((t))\models\PLO$ and $b$ is algebraic but not definable over $A'$. It is easy to build a counter-example to (v) using $A'$ and $b$.

\paragraph{(i) Relative algebraic boundedness:} The proof is not as simple as in
the $p$-adic case and needs some preliminary lemmas and definitions.

\begin{definition}
We will say that $\nT$ is \emph{algebraically bounded (with respect to
$\tT$) within the sort $S$} if for all $\nM\models \nT$ and $A\subseteq \dom(\nM)$,
$S(\ACL(A))\subseteq S(\acll(A))$.

Even if $S$ is stably embedded, one must beware that this is, in general, slightly
different from saying that $\THL(S)$ (the theory induced by $\nT$ on the sort
S) is
algebraically bounded (with respect to $\Thl(S)$), as in the latter case, one
requires that $S(\ACL(A))\subseteq S(\acll(A))$ holds for all $A\subseteq S$.
\end{definition}

\begin{lemma}[k G alg bounded]
Let $T_{F} \supseteq \HFO$ be such that $\inv{\res}/\gpow{(\inv{\res})}{n}$ is finite and $\inv{\res} = \gpow{(\inv{\res})}{n}\resf(F)$.  Then:
\begin{thmenum}
\item If $A = \acll(\K(A))\cap \nM$, then $\Valgp(A) = \val(\K(A))$;
\item If $\THL(\res)$ and $\THL(\Valgp)$ are algebraically
bounded, then $\nT$ is algebraically bounded within $\res$ and
$\Valgp$.
\end{thmenum}
\end{lemma}

\begin{proof}
\begin{thmenum}
\item For any $a\in \inv{\K(\nM)}$ and $\gamma\in\Valgp(\nM)$ such that $n\gamma = \val(a)$
  for some $n\in\Nn$, there exist $x\in \K(\nM)$ such that $\val(ax^{-n}) =
  0$ and $c\in F$ such that \(\val(c) = 0\) and $\resf(ax^{-n}c^{-1}) \in \gpow{\res}{n}$. As $\nM$ is an equicharacteristic zero henselian field, $ax^{-n}c^{-1} \in \gpow{(\inv{\K(\nM)})}{n}$ and hence $ac^{-1}\in \gpow{(\inv{\K(\nM)})}{n}$. So there exists \(a'\in \alg{F(a)}\cap \nM\subseteq \acll(a)\cap \nM\) such that \((a')^n = ac^{-1}\) and hence $\val(a') = \gamma$. As $\Valgp(\acll(\K(A))) = \Qq\tensor\langle\val(\K(A))\rangle$, the
  statement follows.
\item Delon shows in \cite[Theorem\,2.1]{DelPhD} that in the three-sorted language $(\K,\res,\Valgp)$ with $\val$ and $\resf$, field quantifiers can be
  eliminated up to formulas of the form
  $$ \phi^*(x,r) = \exists
  y \in \K\,\bigwedge_i y_ix_i\in \gpow{(\inv{\K})}{n_{i}}\wedge \val(y_i) = 0\wedge
  \phi(r,\resf(y)), $$
  where $r$ is a tuple of variables from $\res$ and \(\phi\) is a formula in the ring language. It follows immediately that
  if $A\subseteq \K(\nM)$ then $\Valgp(\ACL(A)) \subseteq \ACL(\val(A)) \subseteq
  \acll(\val(A)) \subseteq \acll(A)$, where the first inclusion follows from
  field quantifier elimination and the second from algebraic boundedness of $\THL(\Valgp)$.

  The presence of the $\phi^*$ makes it a little more complicated for $\res$,
  but $\phi^*(a,r)$ implies that $a_iy_i \in \gpow{(\inv{\K})}{n_{i}}$ for some $y_i$ such that
  $\val(y_i) = 0$ and hence that $n_{i}|\val(a_i)$. By the first statement, there
  exist $b_i\in \acll(A)\cap \nM$ such that $n\val(b_i) = \val(a_i)$, thus
  $\phi^*(a,r) \iff \exists y\in \res\,\bigwedge_i y_i\resf(a_i b_i^{-n})\in\gpow{(\inv{\res})}{n_{i}} \wedge
  \phi(r,y)$. It follows that any formula with variables in \(\res\) and parameters in \(A\) can be rewritten as a formula with parameters in \(\resf(\acll(A)\cap \nM)\) and hence that $\res(\ACL(A)) \subset
  \ACL(\resf(\acll(A)\cap \nM))$. We now conclude as for $\Valgp$.\qedhere
\end{thmenum}
\end{proof}

In the next three lemmas, we will suppose that the hypotheses of the
previous lemma apply to $\nT$.

\begin{lemma}[RV alg bounded]
For all $A\subseteq \K(\nM)$, $\RV(\ACL(A))\subseteq \RV(\acll(A))$, i.e., $\nT$
is algebraically bounded within $\RV$.
\end{lemma}

\begin{proof}
Let $c\in \RV(\ACL(A))$ and let $\gamma = \val_{\rv}(c)$. Then by
Lemma~\ref{lem:k G alg bounded}, $\gamma \in \Qq \tensor \val(A)$. It follows that there exist $c'\in \K(\acll(A)\cap \nM)$ and $n\in\Nn$ such
that $\val(c') = n\gamma$. Then $c^{n}/\rv(c') \in \res(\ACL(A))\subseteq
\res(\acll(A))$---also by Lemma~\ref{lem:k G alg bounded}---and hence $c \in \acll(A)$.
\end{proof}

\begin{lemma}[ACL balls]
For any $A = \acll(\K(A))\cap \nM$, $\BallSet(\ACL(A))=\BallSet(A)$. Moreover, any ball $b\in\BallSet(\ACL(A))$ contains a point in $A$.
\end{lemma}

\begin{proof}
Let $b\in \BallSet(\ACL(A))$ and let $Q$ be the intersection of all balls in
$\BallSet(A)$ that contain $b$. As $Q$ is $\aut(\tM/A)$-invariant, it suffices to
show that $b$ contains $Q$ (and hence is equal to $Q$) to show it is  $\aut(\tM/A)$-invariant and thus in $\dcll(A)\cap \nM = A$.

If $Q(A) = \emptyset$, it follows from Remark~\ref{rem:D3.0 HF0} that $Q$
is a complete type over $A$ in $\nM$, so $Q$ is contained in $b$. Hence we
can assume that we have a point $a\in Q(A)$. We can suppose $a\not\in b$,
or, because $\rad(b) \in \Valgp(\ACL(A)) \subseteq \Valgp(\acll(A)\cap \nM) =
\Valgp(A)$, we would be done.

If $Q$ is a closed ball that strictly contains $b$, then $b$ is contained in a unique maximal open subball $b'$ of $Q$. This subball $b'$ is equal to the set \(\{x\in \K\mid \rv(x-a) = rv(b-a)\}\) and hence \(b'\) is interdefinable over $A$ (in $\tM$) with $\rv(b-a) \in \RV(\ACL(A)) \subseteq \RV(\acll(A)\cap \nM) = \RV(A)$, where the first inequality follows from Lemma~\ref{lem:RV alg bounded}. Hence $b'$ is in $\BallSet(A)$, contains $b$ and is strictly contained in $Q$, contradicting the definition of $Q$.

Finally, if $Q$ is a strict intersection or an open ball, then $\val(b-a) \in \Valgp(\ACL(A)) = \Valgp(A)$, thus the closed ball of radius $\val(b-a)$ around $a$ would be in $A$, would contain $b$ and would be strictly contained in $Q$, a contradiction.

As for the second point, once we know that $b\in\BallSet(A)$, then---since $\acll(A)$ is a model of $\ACVF$---$b$ contains a point $c$ in $\K(\acll(A)) = \alg{\K(A)}$ and---as balls are convex in residue characteristic zero---the average of the $\aut(\tM,A)$-conjugates of $c$ is in $b(\dcll(A)\cap \nM) = b(A)$.
\end{proof}

\begin{lemma}[rel alg bounded]
For any $A \subseteq \dcll(\K(A))\cap \nM$, $\ACL(A) \subseteq   \acll(A)$. In particular, for any $\nM'\prec \nM$ and $c\in \K(\nM)$, $\ACL(\nM'c) \subseteq \acll(\nM'c)$.
\end{lemma}

\begin{proof}
Let $C = \acll(A)\cap M$, so that $C = \acll(K(C))\cap M$, and let $e\in \ACL(A)$. If $e\in \K \subseteq\BallSet$, then Lemma~\ref{lem:ACL balls} applies to $e$---viewed as a
ball with an infinite radius---and we have $e\in C\subseteq \acll(A)$.

The remaining sorts $\Latt[n]$ and $\Tor[n]$ can be viewed as $\triang_n(\K)/H$ (or
a union of such in the case of $\Tor[n]$) where $H$ is an $\tL$-definable
subgroup. Note that there exists an increasing sequence of $\tL$-definable
subgroups $(G_{i})_{{i=1\ldots m}}$ of $\triang_n(\K)$ with $G_{0} = \{1\}$ and
$G_{m} = \triang_n(\K)$ such that for every $i$, there exists an $\tL$-definable
morphism $\phi_{i} : G_{i} \to G$ with kernel $G_{i-1}$, where $G$ is
either the additive group $\Gg_{a}(\K)$, or the multiplicative group $\Gg_{m}(\K)$, and such that for every point $a\in G(C)$, $\phi_{i}^{-1}(a)$ contains
a point in $G_{i}(C)$. It suffices to show by induction on $i$ that if
$H_{i}:  = G_{i}\cap H$ is an $\tL$-definable subgroup of $G_{i}$ and $e\in (G_{i}/H_{i})(\ACLeq(C))$ then $e$ is $\tL(C)$-definable.

Let $\phi_i : G_{i} \to G$, where $G = \Gg_{a}(\K)$ or $G = \Gg_{m}(\K)$, be a group
homomorphism with kernel $G_{i-1}$. Then $e\in  (G_{i}/H_{i})(\ACLeq(C))$ can be
viewed as an almost $\nL(C)$-definable coset $e H_{i} \subseteq G_{i}$---i.e., a finite union of these cosets is $\nL(C)$-definable---and $\phi_i(e H_{i})$ is an almost $\nL(C)$-definable coset of $\phi_i(H_{i})$. Moreover, the group $H:=\phi_i(H_{i})$ is an $\tL$-defined subgroup of $G$. If $G=\Gg_{a}$ then $H$ has
the form $y\Val$ or $y\mathcal{\nM}$, and its cosets are balls. If $G = \Gg_{m}$ then either $H = 1
+ I$ where $I$ is some proper ideal of $\Val$, and its cosets are balls, or $H= \inv{\Val}$, and its cosets are of the form \(y\inv{\Val} = \val^{-1}(\gamma)\) for some \(y\), \(\gamma\). In both cases, $\phi_i(eH_{i})$ has a point $a\in C$: in the ball case, apply Lemma~\ref{lem:ACL
  balls}, and in the other case, this is because we must have $\gamma\in\Valgp(\ACL(C)) = \Valgp(C) =
\val(\K(C))$, by Lemma~\ref{lem:k G alg bounded}.

Let $a'\in\phi_i^{-1}(a)\cap G_{i}(C) = (a'G_{i-1})\cap G_{i}(C)$; then
$a'^{-1}(eH_i\cap a'G_{i-1})$ is a coset of $H_{i-1} = H_i\cap G_{i-1}$ in
$G_{i-1}$ that is almost $\nL(C)$-definable. By induction,
$a'^{-1}(eH_i\cap a'G_{i-1})$ is $\tL(C)$-definable, but then $(eH_i\cap a'G_{i-1})$ is also $\tL(C)$-definable and hence $eH_i$---the only coset of $H_i$ that contains $eH_i\cap a'G_{i-1}$---is $\tL(C)$-definable.
\end{proof}

\paragraph{(ii)   Internalizing $\tL$-codes:} Let $L = \prod
\Qq_{p}/\mathcal{U}$ be a non-principal ultraproduct. Provided we have a subfield of constants $F$ such that every finite extension of $L$ is generated by an element whose minimal polynomial is over $F$ and which also generates the valuation ring over $\Val(L)$, the proof for finite extensions of $\Qq_{p}$ goes through for $\Th(L)$.

\paragraph{(iii) Unary EI:}
In the following lemmas, we will consider a theory $T_{F}$ extending $\PLO$ where we have
added constants $F$ containing a uniformizer $\unif$, such that $\resf(F)$
contains the necessary constants for $\res$ to have EI and for all
$n\in\Nn_{>0}$, $\inv{\res} = \gpow{(\inv{\res})}{n}\resf(F)$. Let $\nM\models
T_{F}$ be sufficiently saturated and homogeneous.

We will first study the imaginaries in $\RV$. For all $\gamma\in\Valgp(\nM)$, let us
write $\RV[\gamma] := \val_{\rv}^{-1}(\gamma)$. Let $H$ be a (small\footnote{With respect to the saturation and homogeneity of $\nM$.}) subgroup of $\Valgp(\nM)$
containing $1 := \val(\unif)$, and let $\RV[H] = \bigcup_{\gamma\in H} \RV[\gamma]$ where a point
$0_{\gamma}$ is added to every $\RV[\gamma]$. The structure induced by $T_{F,H}$ on
$\RV[H]$ is that of an enriched family of (1-dimensional) $\res$-vector spaces and we view it as a structure with one sort for each $\RV[\gamma]\cup
\{0_{\gamma}\}$. As $H$ is a group, $\RV[H]$ is closed under tensor products and duals.

These $\res$-linear structures are studied in
\cite{HruGpIm}. Let us recall some of the definitions there.

\begin{definition}
Let $A = (V_{i})_{i\in I}$ be a $\res$-linear structure.
\begin{thmenum}
\item We say that $A$ {\em has flags} if for any vector space $V_{i}$ in $A$ with $\dim(V_{i})>1$, there
  are vector spaces $V_{j}$ and $V_{l}$ in $A$ with $\dim(V_{j}) =
  \dim(V_{i})-1$,  $\dim(V_{l}) = 1$ and a $\emptyset$-definable exact sequence $0\to
  V_{l}\to V_{i}\to V_{j}\to 0$.
\item We say that $A$ {\em has roots} if for any $1$-dimensional $V_{i}$ and any
  $m\geq 2$, there exist $V_{j}$ and $V_{l}$ in $A$ and $\emptyset$-definable
  $\res$-linear embeddings $f : V_{j}^{\tensor m} \to V_{l}$ and $g : V_{i}\to
  V_{l}$ such that $\im(g) \subseteq \im(f)$.
\end{thmenum}
\end{definition}

\begin{lemma}[EI RVH]
The theory of $\RV[H]$ with the structure induced by $T_{F,H}$ eliminates imaginaries.
\end{lemma}

\begin{proof}
It follows from \cite[Proposition\,5.10]{HruGpIm} that it suffices to
show that $\RV[H]$ has flags and roots. As every $\RV[a]$ is 1-dimensional,
the structure trivially has flags. But it does not have roots. Let us
extend $H$ to some $H'$ such that $\RV[H']$ has roots.

Let $R = \{r\in\Nn_{>0} : \res(\nM)$ contains nontrivial $r$th roots of unity$\}$, let $L
= \K(\nM)[\unif^{1/r}\mid r\in R]$ and let $H' = \langle H, 1/r \mid r\in R\rangle
\subseteq \val(L)$. Note that $L$ is a ramified extension of $\K(\nM)$
and that $\resf(L) = \res(\nM)$, hence $\RV[H](\nM) = \RV[H](L)$. Now $\RV[1]$ has
$r$th roots in $\RV[H']$ for any $r$. Indeed, if $r\in R$ then $\RV[1/r]$
is an $r$th root and if $r\not\in R$, then as the map $\RV\to \RV : x \mapsto x^r$ is injective, $V_1$ is its own $r$th root.

Let us show that for any $\gamma\in H'$ and any $r\geq 2$,
$\RV[\gamma]$ has an $r$th root. As $\gamma\in H'$, there exists
$n\in\Nn$ such that $n\gamma\in H \subseteq \Valgp(\nM)$, a $\Zz$-group. Hence there
exist $\alpha \in H$ and $m\in\Nn$ such that $n\gamma = rn\alpha + m$. Let
$\RV[\beta]$ be an $nr$th root of $\RV[1]$; then $\RV[\alpha]\tensor
\RV[\beta]^{\tensor m}$ is an $r$th root of $\RV[\gamma]$. By \cite[Proposition\,5.10]{HruGpIm}, $\RV[H']$ has elimination of imaginaries.

Any automorphism $\ns$ of $\RV[H]$ can be extended
to an automorphism of $\RV[H']$. Indeed, if $h\in \RV[H']$ then
$\val_{\rv}(h) = \gamma + n/r$ where $\gamma\in H$, $n\in \Zz$ and $r\in R$,
and $h \rv(\unif)^{-n/r} \in \RV[H]$. Taking $\ts(h) := \ns(h
\rv(\unif)^{-n/r})\rv(\unif)^{n/r}$ will work. Moreover, we can find an
automorphism of $\RV[H']$ fixing only $\RV[H]$. Consider the homomorphism $\phi : H' \to \res(\nM)$ sending $\gamma + n/r$
to $d_r^{n}$ where $(d_r)_{r\in\Nn}\in \res(\nM)$ is such that for all $r$ and $l$, we have $d_r^r = 1$, $d_r\neq 1$ if $r\in R$ and $d_{lr}^l = d_r$. Then $\theta : h \mapsto h\phi(\val_{\rv}(h))$ is a group automorphism of $\RV[H']$ inducing the identity on both $\res$ and $H'$ hence an automorphism of the full structure of $\RV[H']$. It is easy to see that $\theta$ fixes only $\RV[H]$.

Note that because each fiber is a sort, if $X\subseteq \RV[H]^{l}$ for
some $l\in\Nn$ and $X$ is definable in $\RV[H]$, then it is defined by the same formula in $\RV[H']$. Hence it is coded by some $x\in \RV[H']$. But as there are automorphisms of $\RV[H']$ fixing only $\RV[H]$, we must have $x\in \RV[H]$, and as automorphisms of $\RV[H]$ extend to $\RV[H']$, $x$ is also a code for $X$ in $\RV[H]$.
\end{proof}

\begin{proposition}[EI RV]
The theory induced by $T_{F}$ on the sort $\RV$ (see Section~\ref{sec:Vfield}) eliminates imaginaries to the sorts $\RV$ and $\Valgp$.
\end{proposition}

\begin{proof}
First let us show that for all $n\in\Nn_{>0}$, $\RV/\gpow{\RV}{n}$ is
finite and $\RV = \gpow{\RV}{n}\rv(F)$. Let $a\in\RV$.  As $\Valgp$ is a $\Zz$-group, there exist $y\in\RV$ and
$r\in\Nn$ such that $r<n$ and $\val_{\rv}(a) = \val_{\rv}(y^{n}) +
\val(\unif^{r})$. Hence $\val_{\rv}(ay^{-n}\rv(\unif)^{-r}) = 0$,
i.e., $ay^{-n}\rv(\unif)^{-r}\in \inv{\res}$. As $\inv{\res} = \gpow{(\inv{\res})}{n}\resf(F)$,
there exists $m\in \resf(F)$ such that $a y^{-n}m^{-1}\rv(\unif^{-r})\in\gpow{(\inv{\res})}{n}$, i.e., $a \in m\rv(\unif^{r})\gpow{\RV}{n}$.

Moreover, for any $A\subseteq \RV(\nM)$, $\val_{rv}(\DCL(A)) \subseteq
\Qq\tensor\val_{\rv}(A)$. Indeed, let
$\gamma\in\Valgp(\nM)\setminus\Qq\tensor\val_{\rv}(A)$ and $d \in \bigcap
\gpow{(\inv{\res(\nM)})}{n}\setminus\{1\}$; then there exists a group homomorphism
$\phi_{d} : \Valgp(\nM) \to \inv{\res}(\nM)$ such that $\phi_{d}(\val_{\rv}(A)) =
\{1\}$, $\phi_{d}(\gamma) = d$ and $\psi_{d} : t \mapsto t \phi_{d}(\val_{\rv}(t))$ defines an automorphism of
$\RV(\nM)$ fixing $A$, $\res$ and $\Valgp$, which sends any $x\in
\val_{\rv}^{-1}(\gamma)$ to $dx\neq x$. Hence $\val_{\rv}^{-1}(\gamma)$ cannot
contain any point definable over $A$.

Let us now code finite sets. For any tuple $\gamma\in\Valgp$, let
$\RV[\gamma]$ denote $\prod_{i}\RV[\gamma_{i}]$.

\begin{claim}[RV fin set]
In the theory induced by $T_{F}$ on the sorts $\RV\cup\Valgp$, finite sets are coded.
\end{claim}

\begin{proof} Let $X \subseteq\RV^{i}\times\Valgp^{j}$ be finite. As $\Valgp$
  is ordered, we can suppose that there are tuples $\gamma$ and $\gamma'\in\Valgp$ such
  that $X\subseteq \RV[\gamma]\times\{\gamma'\}$. By Lemma~\ref{lem:EI RVH}, the
projection of $X$ on $\RV[\gamma]$ is coded (over $\gamma$) by some $x\in\RV[\langle1,\gamma\rangle]$. It is easy to see
that $x\gamma\gamma'$ is a code for $X$.
\end{proof}

To prove elimination of imaginaries in $\RV$ to the sorts $\RV$ and
$\Valgp$, it suffices, by Lemma~\ref{lem:1}, to code $\nL(A)$-definable functions
$f : \RV \to R$, where $R$ is either $\RV$ or $\Valgp$, for any $A\subseteq \RV(\nM)$. Let us first
consider the case $R = \RV$. Let $D$ be the domain of $f$ and $X$ its
graph.

\begin{lemma}[lin eq fun]
If there exist $n$ and $m\in\Zz$ such that for all $x\in D$, $n\val_{\rv}(f(x)) - m\val_{\rv}(x)$ is constant, then $f$ is coded.
\end{lemma}

\begin{proof}
Let $\gamma_{f} = n\val_{\rv}(f(x)) - m\val_{\rv}(x)
\in\Valgp(\DCL(\code{f}))$. For all $y\in\inv{\res}$ and $x,z\in\RV$, let $y\cdot (x,z) = (y^{n}x,y^{m}z)$. This defines an action of
$\inv{\res}$ on any $\RV[\gamma]$ where $\gamma$ is a 2-tuple. Let
$y\in\bigcap_{n}\gpow{(\inv{\res})}{n}$ and $\gamma\in \cpow{\Valgp(\nM)}{2}$ be such
that $n\gamma_{2}-m\gamma_{1} = \gamma_{f}$ and
$\gamma_{1} \not\in \Qq\tensor\langle\val_{\rv}(A)\rangle$. By a similar
automorphism construction as above, there is $\psi\in \aut(\RV(\nM)/A)$ such that for all
$x\in\RV[\gamma]$, $\psi(x) = y\cdot x$ and hence $x\in X$ implies $y\cdot x\in X$. By compactness, there exists
$N\in\Nn_{>0}$ such that for any $x\in\RV$ with
$\val_{\rv}(x)\not\in\Qq\tensor\langle\val_{\rv}(A)\rangle$ and for any $y\in
\gpow{(\inv{\res})}{N}$,  if $x\in X$ then $y\cdot x\in
X$. Let $X'=\{x\in X \mid \forall y\in\gpow{(\inv{\res})}{N}, y\cdot x\in X\}$. Then it suffices to
code $X'$ and $X\setminus X'$. Note that $(x,y)\in X\setminus X'$ implies
$\val_{\rv}(x) \in\Qq\tensor\langle\val_{\rv}(A)\rangle$.

\begin{claim}[k stable]
Suppose that $X$ is stable under the action of $\gpow{(\inv{\res})}{N}$.  Then $f$
is coded.
\end{claim}

\begin{proof}
Let $E \subseteq \rv(F)$ intersect all the classes of $\RV$ modulo
$\gpow{\RV}{(Nn)}$. Fix $\gamma\in\Valgp$. For any $x\in D_{\gamma} := D\cap\RV[\gamma]$, there
exist $y\in\gpow{\RV}{N}$ and $e\in E$ such that $x =
y^{n}e$. As $X$ is $\gpow{(\inv{\res})}{N}$-stable, one can check
that $g_{\gamma}(e) := y^{-m}f(x)$ depends only on $e$ and
$\gamma$. One can also check that $\val_{\rv}(g_{\gamma}(e)) = \frac{1}{n}(\gamma_{f} +
m\val_{rv}(e))\in \Valgp(\DCL(\code{f})) =: H$ and
$g_{\gamma}$ is in fact a function (with a finite graph $G_{\gamma}$) definable in
$\RV[H]$. By Lemma~\ref{lem:EI RVH} and compactness, there is a definable function $g :
\Valgp\to\RV[H]^{l}$ for some $l\in\Nn$ such that $g(\gamma)$ codes
$g_{\gamma}$ (over $H$). It is quite clear that $g$ is $\nL(\code{f})$-definable, but as
$X = \bigcup_{\gamma\in\Valgp} \gpow{(\inv{\res})}{N} G_{\gamma}$, $f$ is also $\nL(H\code{g})$-definable.

Now, as $\Valgp$ has Skolem functions, we can
definably order $\im(g)$, and, because $\RV[H]^{l}$ is internal to $\res$ and
the induced theory on $\res$ is simple, $\im(g)$
must be finite (a simple theory cannot have the strict order property). Thus $\im(g) \subseteq \ACL(\code{f})$. For any $e\in
\im(g)$, $g^{-1}(e)\subseteq\Valgp$ is coded. Let $d$ be the tuple of all
codes of fibers and corresponding images; then
$d\in \ACL(\code{f})$ and $\code{f}\in\DCLeq(\gamma d)$ for some $\gamma\in H = \Valgp(\DCL(\code{f}))$. We can
conclude by coding the finite set of $\code{f}$-conjugates of $\gamma d$ (by
Lemma~\ref{claim:RV fin set}).
\end{proof}

\begin{claim}[small dom]
Suppose that for all $x\in D$, $\val_{\rv}(x)\in
\Qq\tensor\langle\val_{\rv}(A)\rangle$.  Then $f$ is coded.
\end{claim}

\begin{proof}
By compactness, $D$ must be contained in only finitely many
$\RV[\gamma_{i}]$. All of these $\gamma_{i}$ are $\nL(\code{f})$-definable
and hence $f$ lies inside $\RV[H]$, where $H: = \Valgp(\DCL(\code{f}))$. By Lemma~\ref{lem:EI RVH}, $f$ is
coded by some $d$ over $H$, hence there is some tuple
$\gamma\in H$ such that $d\gamma$ codes $f$.
\end{proof}

Now, Claim~\ref{claim:k stable} allows us to code $X'$ and Claim~\ref{claim:small
  dom} allows us to code $X\setminus X'$. This concludes the proof of Lemma~\ref{lem:lin eq fun}.
\end{proof}

Let us now show that we can reduce to Lemma~\ref{lem:lin eq fun}. As $f(x)\in\DCL(Ax)$, we have $\val_{\rv}(f(x)) \in\Qq\tensor\langle\val_{\rv}(Ax)\rangle$. By compactness, for all $i$ in some finite set $I$, there exist $n_{i}$, $m_{i}\in\Zz$ and $\gamma_{i}\in\Qq\tensor\val_{rv}(A)\cap\Valgp(\nM)$ such that for all $x\in D$, there exists $i\in I$ with $g_{i}(x) := n_{i}\val_{\rv}(f(x)) - m_{i}\val_{\rv}(x) = \gamma_{i}$. Define $E_{i,\gamma}$ to be the fiber of $g_{i}$ above $\gamma$. Then $D\subseteq \bigcup_{i\in I} E_{i,\gamma_{i}}$. Let us assume that $\card{I}$ is minimal such that this inclusion holds.

\begin{claim}
The set $X:=\{(\gamma_{i})_{i\in I}\in\Valgp\mid D\subseteq \bigcup_{i\in I} E_{i,\gamma_{i}}\}$ is finite.
\end{claim}

\begin{proof}
We proceed by induction on $\card{I}$. Let us assume $X$ is infinite, and pick any $x\in D$. By the pigeonhole principle, there exists $i_{0}\in I$ and an infinite set $Y\subseteq X$ such that for all $(\gamma_{i})_{i\in I}\in Y$, $x\in E_{i_{0},\gamma_{i_{0}}}$, i.e., $g_{i_{0}}(x) = \gamma_{i_{0}}$. It follows that for all $(\gamma_{i})_{i\in I}$ and all $(\delta_{i})_{i\in I}\in Y$, we have $\gamma_{i_{0}} = \delta_{i_{0}}$ and $E_{i_{0},\gamma_{i_{0}}} = E_{i_{0},\delta_{i_{0}}} =: E$. By minimality of $\card{I}$, $D\setminus E$ is nonempty and the set $\{(\gamma_{i})_{i\in I\setminus\{i_{0}\}}\in\Valgp\mid D\setminus E\subseteq \bigcup_{i\in I\setminus\{i_{0}\}} E_{i,\gamma_{i}}\}$ is finite by induction, but it contains $\{(\gamma_{i})_{i\in I\setminus\{i_{0}\}}\mid (\gamma_{i})_{i\in I}\in Y\}$ which is infinite, a contradiction.
\end{proof}

Then any $(\gamma_{i})_{i\in I}\in X$ is in $\ACL(\code{f})$, $f_{i} := \restr{f}{E_{i,\gamma_{i}}}$ satisfies the conditions of Lemma~\ref{lem:lin eq fun} and it suffices to code each $f_{i}$. Indeed let $d$ be the tuple of the codes for  those functions; then $d\in\ACL(\code{f})$ and, as $f= \bigcup_{i\in I}f_{i}$, $\code{f}\in\DCLeq(d)$. The code of the finite set of $\code{f}$-conjugates of $d$---which exists by Claim~\ref{claim:RV fin set}---is a code for $f$.

Finally, if $R = \Valgp$, then for all $\gamma\in\Valgp(\nM)$, $f^{-1}(\gamma)
\subseteq \RV$ is coded by the case $R = \RV$. Hence $f$ is
interdefinable with a function from $\Valgp$ to
$\RV^{l}\times\Valgp^{m}$ for some $l$ and $m$. So we have to code functions
from $\Valgp$ to $\Valgp$ (which we already know how to code) and from
$\Valgp$ to $\RV$. Let $g : \Valgp\to\RV$ be a definable function and let
$h = g\circ\val_{\rv}$. Then $h : \RV\to\RV$ is coded as we have just
shown and, as for all $\gamma\in\Valgp$, $h(\val_{\rv}^{-1}(\gamma)) =
\{g(\gamma)\}$, a code for $h$ is also a code for $g$. This concludes the proof of Proposition~\ref{prop:EI RV}.
\end{proof}

\begin{remark}
\begin{thmenum}
\item Let $B_m = \RV / \gpow{(\inv{\res})}{m}$. We have a homomorphism $B_m \to \Valgp$ with finite kernel $\inv{\res} / \gpow{(\inv{\res})}{m}$. Hence $\gpow{B_m}{m}$ maps injectively into $\Valgp$, and our assumptions on constants imply that there is a set of $\emptyset$-definable representatives for the cosets of $\gpow{B_m}{m}$ in $B_m$.  Thus the theory (and imaginaries) of $B_m$ reduce to those of $\Valgp$.
\item On the other hand, it can be shown that every unary definable subset $D$ of $\RV$ is a finite union of pullbacks from $B_{m}$ for some $m$ and subsets of $\val_{\rv}^{-1}(a)$ for $a$ lying in some finite subset $F_D$ of $\Valgp$. This $m$ is uniform in families, and $F_D$ can be defined canonically as the set of $a \in \Valgp$ such that $\val_{\rv}^{-1}(a)$ is not a pullback from $B_m$.  This gives another proof of unary EI in $\RV$ (with the stated constants), given EI in any $\RV[H]$.

A similar (but slightly more complicated) decomposition is also true in higher dimension (e.g., adapt \cite[Lemma\,3.25]{HruKaz} to our case by replacing $\Valgp$ with a suitable $B_{m}$). Moreover, EI in $\RV$ also follows from this decomposition.
\end{thmenum}
\end{remark}

Let us come back to unary EI in $T_{F}$ (in fact, the proof given here would work in any theory $T\supseteq \HFO$ such that  $\Valgp$ is definably well-ordered and
$\RV$ has unary EI). We will proceed as in the case of finite extensions of $\Qq_{p}$. First
let us show that the analogue  of Claim~\ref{claim:typeball} is still true in this case.

\begin{claim}
Let $A = \ACLeq(A)$, let $B = \BallSet(A)$ and let $c\in \K(\nM)$.  Then $\tpL(c/B) \vdash \tpL(c/A)$.
\end{claim}

\begin{proof}
Recall from Section~\ref{sec:Vfield} that \(\RV\) is stably embedded and has unary EI. As any element in $\RV$ is coded by a ball, the claim is true if $c\in\RV(\nM)$. Recall that $W(c;A) := \{b \in \BallSet(A) \mid c \in b\}$. If $P := \bigcap
W(c; A) = \bigcap W(c; B)$ does not contain any ball in $B$ then $P$ is a complete type over $A$ and $B$ (by Proposition~\ref{prop:classif types PLO}) and we are done. If $P$ does contain a ball $b\in B$, then, by Proposition~\ref{prop:classif types PLO}, $P$ is complete relative to  $\rv(x - b)$. But $\tpL(\rv(x - b)/B) \vdash \tpL(\rv(x - b)/A)$ and we are also done.
\end{proof}

Unary EI in $T_{F}$ follows as for finite extensions of $\Qq_{p}$.

\paragraph{(iv) Invariant types and germs:} The same proof as for finite extensions of $\Qq_{p}$ (nearly) works as we only used there that $\Valgp$ is definably well-ordered. The one difference is that $P$ can be a closed ball. But in that case $\tP$, the $\ACVF$ generic of $P$, is definable, thus the $\tP$-germ of any $r$ is an imaginary element $e$, and one may take $I = \{0\}$ and $e_0=e$. Moreover, the inconsistency of $\tp(c/A)$ and $\tprestr{\tP}{\nM}$ would---by Claim~\ref{claim:balls from tM}---contradict Lemma~\ref{lem:D3.2 closed}.

\begin{corollary}[EI/UFI]
Let $T_{F}\supseteq \HFO$ be an $\nL$-theory such that $\THL(\res)$ and $\THL(\Valgp)$ are algebraically bounded, $\Valgp$ is definably well-ordered, $\RV$ has unary EI, $\K$ has a finite number of extensions of any given degree and $\inv{\res}/\gpow{(\inv{\res})}{n}$ is finite. Suppose also that we have added constants for a field $F\subseteq \K$ such that $\inv{\res} = \gpow{(\inv{\res})}{n}\resf(F)$ and any finite extension of $\K$ is generated by an element whose minimal polynomial is over $F$ and which also generates the valuation ring over $\Val(\K)$. Then $T_{F}$ has EI/UFI in the sorts $\K$ and $\Latt[n]$.

In particular this is true of ultraproducts of the $p$-adics (if we add some constants as in Remark~\ref{constants asymp}).
\end{corollary}

\begin{proof}
By Proposition~\ref{prop:EI/UFIcrit} we have EI/UFI in the sorts $\K$, $\Latt[n]$ and $T_{n}$ but as noted earlier the sorts $T_{n}$ are not needed when the value group has a smallest positive element.
\end{proof}

\paragraph{Elimination of finite imaginaries:} As we already know that
$\RV$ eliminates imaginaries, it suffices to show that every finite imaginary in $\PLO$ (over arbitrary parameters) can be coded in $\RV$ (the proof is adapted from \cite[Lemma\,2.10]{HruFin}).

\begin{definition}
If $C \subseteq C'$, we say that $C'$ is {\em stationary} over $C$ if $\dcleq(C')\cap \acleq(C) = \dcleq(C)$. A type $p = \tp(c/C)$ is {\em stationary} if $cC$ is stationary over $C$.
\end{definition}

\begin{remark}[stat]
\begin{thmenum}
\item It is clear that if $C''$ is stationary over $C'$ and $C'$ is stationary over $C$, then $C''$ is stationary over $C$.
\item\label{CS stat} If $\tp(c/C)$ generates a complete type over $\acleq(C)$, then
  $\tp(c/C)$ is stationary. Indeed, let $x\in \dcleq(Cc)\cap\acleq(C)$; then
  there is a $C$-definable function $f$ such that $f(c) = x$. As
  $\tp(c/C)$ generates a complete type over $\acleq(C)$, there is a
  $C$-definable set $D$ such that for all $c'\in D$, $f(c') = x$, hence $x\in\dcleq(C)$.
\end{thmenum}
\end{remark}

\begin{lemma}[stat res RV]
Let $\nT$ be a theory extending $\PLO$ (in the geometric language with
possibly new constants). For all $\nM\models \nT$ and $A\subseteq \nM$, there
exists $C\preceq \nM$ containing $\RV(\nM)\cup A$ and stationary over
$\RV(\nM)\cup A$.
\end{lemma}

\begin{proof} Let us first prove the following claim.
\begin{claim}
Let $B = \DCL(B) \subseteq \nM$ such that $\RV(\nM)\subseteq B$ and $b\in
\BallSet(\nM)$. Then there exists a tuple $c\in \K(\nM)$ with $\tpL(c/B)$ stationary,
$b\in\DCL(c)$ and $b(\nM)\cap c \neq\emptyset$.
\end{claim}

\begin{proof}
Let us first suppose that $b\in \RV(\nM)$, i.e., that $b$ is of the form
$c(1+\mathcal{\nM})$. Let $P\subseteq b$ be a minimal (for inclusion)
intersection of balls in $\BallSet(B)$. For any $c\models P$ we have $b = \rv(c)$,
hence it suffices to show that $P$ is a complete stationary type over $B$.

As $P$ does not strictly contain any ball in $\BallSet(B)$ by
definition, it cannot contain a ball
$b'\in\BallSet(\ACLeq(B))$. Indeed, if \(P\) is strict then, taking the
smallest ball containing the orbit of \(b'\) over \(B\), we obtain a strict subball of \(P\) which is in \(\BallSet(B)\), a
contradiction. If \(P\) is a closed ball, then we may assume that
\(b'\) is a maximal open ball in \(P\) and, since \(\res_P(P)\) is a
\(\res\)-torsor, we can take the mean of the orbit over \(B\) (we are
in residue characteristic zero) to get a strict subball of \(P\)
contained in \(\BallSet(B)\), again a contradiction. By
Proposition~\ref{prop:classif types PLO}, $P$ is a complete type over
$\ACLeq(B)$. By Remark~\ref{CS stat}, $P$ is stationary over $B$.

Now if $b\in\BallSet(\nM)$, pick any $r\in \RV(\nM)$ such that $\val_{\rv}(r) =
\rad(b)$. Applying the claim to $r$, we find $c\in \K(\nM)$ such that
$\tpL(c/B)$ is stationary and $\rad(b)\in\DCL(c)$. It now suffices to find
a point $d\in b$ whose type is stationary over $\DCL(Bc)$, but we can
proceed as in the first case. Then $b\in\DCL(cd)$ and $\tp(cd/B)$ is stationary.
\end{proof}

Starting with $B := \DCL(\RV(\nM)\cup A)$, and applying the claim iteratively, we find $C \supseteq A\cup \RV(\nM)$ such that \(C\subseteq \nM\), $C$ is stationary over $A\cup \RV(\nM)$, $\DCL(C) = C$, $\BallSet(\nM) \subseteq \DCL(\K(C))$ and every ball in $\BallSet(\nM)$ has a point in $C$.

\begin{claim}
We have $C\subseteq \DCL(\K(C))$.
\end{claim}

\begin{proof}
Let $e\in C$.  If $e\in \K$ then the result is trivial, thus we only have to consider
$e\in \Latt[n]$ or $e\in \Tor[n]$. Let us consider the same decomposition of
$\Latt[n]$ and $\Tor[n]$ as in the proof of Lemma~\ref{lem:rel alg bounded} and show by induction on $i$ that for all $e\in (G_{i}/H_{i})(\nM)$, $e$ is $\nL(\K(C))$-definable.

If we write $e$ as $eH_{i}$ then, as proved in Lemma~\ref{lem:rel alg bounded}, $\phi_{i}(eH_{i})$ is either a ball or a set of the form $y\inv{\Val}$ and hence is definable over $\BallSet(\nM)$ and has a point $a'\in \K(C)$. Let $a\in \phi_{i}^{-1}(a')(C)$. Then $a^{-1}eH_i\cap G_{i-1}$ is a coset of $H_{i-1}$ in $G_{i-1}$ which is $\nL(\K(C))$-definable by induction.  Since \(a^{-1}eH_i\) contains \(a^{-1}eH_i\cap G_{i-1}\), it follows that \(a^{-1}eH_i\) is $\nL(\K(C))$-definable, and hence, so is $eH_{i}$.
\end{proof}

As $\DCL(C) = C$, we have $\K(C) = \K(C)^{h}\models\HFO$.  Since \(\RV(\nM)\subseteq\RV(C)\), \(C\subseteq \nM\) and every ball in \(\BallSet(\nM)\)---in particular, every element of \(\RV(\nM)\))---has a point in \(\K(C)\), we have that $\rv(\K(C)) = \RV(\nM)$. It follows from field quantifier elimination in $\HFO$ in the language with sorts $\K$ and $\RV$ (see Section~\ref{sec:Vfield}), that
$\K(C)\preceq \K(\nM)$. But this implies that $C = \DCL(\K(C))\preceq \nM$.  This concludes the proof of Lemma~\ref{lem:stat res RV}.
\end{proof}

\begin{lemma}[FI bij RV]
Let $\nT$ be a theory that extends $\PLO$ (in the geometric language) and let $A\subseteq \nM\models \nT$. Then every finite imaginary sort of $\nT_{A}$ is in definable bijection with a finite imaginary sort of $\RV$ (with the structure induced by $\nT_A$).
\end{lemma}

\begin{proof}
Let $Y = D/E$ be a finite imaginary sort (in $\nT_A$) and let $\pi : D \to Y$ be
the canonical surjection. As the field sort is dominant, we can assume that
$D$ is a definable subset of $\cpow{\K}{n}$ for some $n$. Let $C\supseteq A$ be as in Lemma~\ref{lem:stat res RV}. As $Y$ is finite and
$C\prec \nM$, $Y(C) = Y(\nM)$ and there exists a finite set $H \subseteq
\cpow{\K}{n}(C)$ meeting every $E$-class. Let $W$ be some finite set in
$\RV(C)$, of bigger cardinality than $H$, and $h : W \to H$ any
surjection. Note that any such surjection is $\nL(C)$-definable. Composing,
we have an $\nL(C)$-definable surjection $\psi : W \to Y$. But there are
only finitely many maps $W \to Y$, hence they are all algebraic over $\RV(C)\cup
A = \RV(\nM)\cup A$, and by stationarity of $C$ over $\RV(\nM)\cup A$, $\psi$ is
$\nL(\RV(\nM)\cup A)$-definable. Let $e\in \RV(\nM)$ be such that $\psi$ and $W$ are $\nL(Ae)$-definable.

Let $W$ be defined by the $\nL(Ae)$-formula $\phi(x,e)$ and $\psi$ by the $\nL(Ae)$-formula
$\psi(x,y,e)$ (which implies that for any $e'$, $\psi(\nM,\nM,e')$ is the graph
of a function with domain $\phi(\nM,e')$). Then the formulas $\phi(x,z)$ and $\psi(x,y,z)$ define,
respectively, a subset $D'$ of $\RV^{|e|+1}$ and a surjection $\psi : D' \to
Y$. Let $E'$ be defined by $E'((x,z),(x',z')) \iff \forall y\,\psi(x,y,z) = \psi(x',y,z')$. Then
we have an $\nL(A)$-definable bijection $D'/E' \to Y$ and, as $\RV$ is considered
with the structure induced by $\nT_{A}$, $D'/E'$ is a finite imaginary sort of $\RV$.
\end{proof}

\begin{proof}[of Theorem~\ref{thm:ultraprod Qp}]
 Let $K \models \PLO$ and let $\nT = \Th(K)$ (with constants added as in
Corollary~\ref{cor:EI/UFI}). As we have already proved EI/UFI in Corollary~\ref{cor:EI/UFI}, by Lemma~\ref{lem:critHfin} it is enough to show that for any
$A$, $T_{A}$ eliminates finite imaginaries in the sorts $\K$, $\Latt[n]$. Let $e\in\ACLeq(A)$; then, by Lemma~\ref{lem:FI bij RV}, there exists an $\RV$-imaginary $e'$ interdefinable over $A$ with $e$. By EI in $\RV$ to the sorts $\RV$ and $\Valgp$ (Proposition~\ref{prop:EI RV}), there exists a tuple $d\in\RV\cup\Valgp$ such that $e'$ is interdefinable with $d$, hence $e$ is interdefinable with $d$ over $A$. We have shown that any
finite imaginary of $T_{A}$ is coded (over $A$) in $\RV\cup\Valgp =
\Tor[1]\cup \Latt[1]$, and the points of $\Tor[1]$ and $\Latt[1]$ are themselves coded in $\Latt[2]\cup \Latt[1]$.
\end{proof}

For a more canonical  treatment of the parameters $F$ in the pseudo-finite case,  see \cite{ChaHruACFA}---it would be interesting to adapt \emph{op.\ cit.}\ to the pseudo-local setting.

\section{Rationality}
\label{sec:rat}

Let $r\in \Nn$. For all tuples $l\in\Nn^{r}$, when $t = (t_{i})_{1\leq i \leq r}$, we write $t^{l}$ for $\prod_{i\leq r} t_{i}^{l_{i}}$. We say a power series $\sum_{l \in \Nn^r}a_l t^l \in\Qq[[t_{1},\ldots,t_{r}]]$ with each $a_l\in \Nn$ is {\em rational} if it is equal to a rational function in $t_1,\ldots, t_r$ with coefficients from $\Qq$.  In this section we prove that certain zeta functions that come from counting the equivalence classes of definable equivalence relations are rational.

For any finite extension $L_{p}$ of $\Qq_{p}$, it is natural here to consider the invariant Haar measure $\mu_{L_{p}}$ on $\GL{N}(L_p)$.  In terms of the additive Haar measure $\mu_{L_{p},+}^{N^2}$ on $\cpow{L_p}{N^{2}}$, $\mu_{L_{p}}$ can be defined thus:  for any continuous $f\colon \GL{N}(L_p)\to \Cc$ with compact support, $\int f(x) d\mu_{L_p}(x) = \int f(x) |\det(x)|^{-N} d \mu_{L_{p},+}^{N^2} (x)$.  As $\det(x)$ is uniformly definable for all $L_p$, Denef's results on definability of $p$-adic integration \cite{DenMSRI} extend immediately to $d\mu_{L_p}$ and the motivic counterpart of these results---see \cite{DeLoICM}, although the result we will be needing is already implicit in older work by Denef and Pas (see, e.g., \cite{PasEQ})---also extend to $d\mu_{L_p}$.

By left invariance, $\mu_{L_p}(A \cdot \GL{n}(\Val(L_p))) = \mu_{L_p}(\GL{n}(\Val(L_p)))$, a number that depends only on the normalization.  We choose a normalization for $\mu_{L_{p},+}$
and $\mu_{L_p}$ such that for any $A \in \GL{N}(L_p)$, we have

 \begin{equation}
 \label{eqn:ballmeas}
  \mu_{L_p}(A.\GL{N}(\Val(L_p)))=1.
 \end{equation}

Let $K$ be a number field and let $\Val_{K}$ denote its ring of integers. For each prime $p$, let $\mathfrak{F}_{p}$ be a set of finite extensions of $\Qq_{p}$, each containing $K$, and let $\mathfrak{F} = \bigcup_{p}\mathfrak{F}_{p}$. We will say that $(R_{L_{p}})_{L_{p}\in\mathfrak{F}}$ and $(E_{L_{p}})_{L_{p}\in\mathfrak{F}}$ are \emph{uniformly $K$-definable in $\mathfrak{F}$}, or just  \emph{uniformly $K$-definable}, if there exist two $\LG(K)$-formulas $\phi$ and $\theta$---i.e., $\LG$-formulas with parameters in $K$---independent of $L_{p}$ such that for all $L_{p}\in\mathfrak{F}$, $R_{L_{p}} = \phi(L_{p})$ and $E_{L_{p}} = \theta(L_{p})$.  If $K = \Qq$ then we often write \emph{uniformly $\emptyset$-definable in $\mathfrak{F}$} instead of uniformly $\Qq$-definable in $\mathfrak{F}$.  If in addition $\mathfrak{F}_p= \{\Qq_p\}$ for all $p$, then we often write \emph{uniformly $\emptyset$-definable in $p$} instead of uniformly $\Qq$-definable in $\mathfrak{F}$.

By a {\em (uniformly $K$-)definable family} $R_{L_p}=(R_{L_p,l})_{l \in \Zz^r}$ of subsets of $L_p^N$ we mean a (uniformly $K$-)definable subset $R_{L_p}$ of $L_p^N \times \Zz^r$---where $\val(\inv{L_{p}})$ is identified with $\Zz$---and we write $R_{L_p,l}$ for the fiber above $l$ of the projection from $R_{L_p}$ to $\Zz^r$.  By a {\em (uniformly $K$-)definable family} $E_{L_p}=(E_{L_p,l})_{l \in \Zz^r}$ of equivalence relations on $R_{L_p}$ we mean a (uniformly $K$-)definable equivalence relation $E_{L_p}$ on $R_{L_p}$ such that for every $x,y \in R_{L_p}$, if $xE_{L_p}y$ then there exists $l \in \Zz^r$ such that $x,y \in R_{L_p,l}$.  We then have a (uniformly $K$-)definable equivalence relation $E_{L_p,l}$ on $R_{L_p,l}$ for every $l$, and by a slight abuse of notation we can regard $(E_{L_p,l})_{l \in \Zz^r}$ as a (uniformly $K$-)definable family of subsets of $L_p^{2N}$. The set $\Nn^r$ is a (uniformly $K$-)definable subset of $\Zz^r$, so it makes sense to talk of (uniformly $K$-)definable families $R_{L_p}=(R_{L_p,l})_{l \in \Nn^r}$, etc.

Now we come to the main result of this section (cf.\ \cite[Thms.~1.3 and 1.4]{duSGruAnn}).

\begin{theorem}[rat]
Let $\mathfrak{F}_{p}$ and $\mathfrak{F}$ be as above (note that we do not assume  $\mathfrak{F}_{p}$ is nonempty for infinitely many $p$).  For all $L_{p}\in\mathfrak{F}$, let $R_{L_{p}}=(R_{L_p,l})_{l \in \Nn^r}$ be a family of subsets of $L_p^N$  and let $E_{L_{p}}= (E_{L_p,l})_{l \in \Nn^r}$ be a family of equivalence relations on $(R_{L_{p},l})_{l \in \Nn^r}$ such that $(R_{L_{p}})_{L_p\in \mathfrak{F}}$ and $(E_{L_{p}})_{L_p\in \mathfrak{F}}$ are uniformly $K$-definable in $\mathfrak{F}$. Suppose that for each $l \in \Nn^r$ and each $L_{p}\in {\mathcal{L}}$, the set of equivalence classes $R_{L_p,l}/E_{L_p,l}$ is finite.  Let $a_{L_{p},l}=\card{R_{L_p,l}/E_{L_p,l}}$. Then the power series
 \[S_{L_{p}}(t) := \sum_{l \in \Nn^r}a_{L_{p},l} t^l \in\Qq[[t_{1},\ldots,t_{r}]]\]
is rational for every $L_{p}\in\mathfrak{F}$.

Moreover, there exist $k,n,d \in\Nn$, there exist tuples $(a_{j})_{j\leq k}$ of integers and $(b_{j})_{j\leq k}$ of elements of $\Nn^{r}$, and for all tuples $l\in\Nn^r$ with $ \card{l} := \sum_{i\leq r} l_{i}\leq d$ there exist $q_{l}\in\Qq$ and varieties $X_{l}$ over $\Val_{K}$, such that the following holds:
\begin{enumerate}
\item[(1)] for all $j$, $a_{j}$ and $b_{j}$ are not both 0; and
\item[(2)] for all $p\gg 0$ and all $L_{p}\in\mathfrak{F}_{p}$, we have
\begin{equation}
\label{eqn:unifrat}
 S_{L_{p}}(t) = \frac{\sum_{\card{l}\leq d} q_{l}\card{X_{l}(\resf(L_{p}))} t^{l}}{\card{\resf(L_{p})}^{n}
\prod_{j=1}^k(1-\card{\resf(L_{p})}^{a_{j}}t^{b_{j}})}.
\end{equation}
\end{enumerate}
\end{theorem}

Suppose we are given power series $S_{L_{p}}(t)= \sum_{l \in \Nn^r}a_{L_{p},l} t^l \in\Qq[[t_{1},\ldots,t_{r}]]$ for each $L_p\in \mathfrak{F}$.  We say the power series $(S_{L_{p}}(t))_{L_p\in \mathfrak{F}}$ are {\em uniformly rational for $p\gg 0$} if there exists a prime $p_0$ such that the $S_{L_{p}}(t)$ are of the form given in (\ref{eqn:unifrat}) for all $L_p\in \mathfrak{F}$ such that $p> p_0$.

\begin{remark}
\label{rem:allprimes}
\begin{thmenum}
\item Assume $\mathfrak{F}_{p}$ is finite for all $p$ (this is the case in most of our applications in Sections~\ref{sec:zetagp} and \ref{sec:twst}). Let $\Lrg$ be the language of rings. At the cost of replacing the $X_{l}$ with quantifier-free $\Lrg(\Val_{K})$-definable sets, we can make (\ref{eqn:unifrat}) hold for every $L_p$, where $k,n,d$, the $a_j$, the $b_j$, the $q_l$ and the $X_l$ are all independent of the choice of $L_p$.  In this case, we say the power series $(S_{L_{p}}(t))_{L_p\in \mathfrak{F}}$ are {\em uniformly rational}.  In particular, suppose we are given definable $R_{p_0}$ and $E_{p_0}$ as above, but just for a single prime $p_0$ and a single $L_{p_0}$.  Then taking $\mathfrak{F} = \{L_{p_{0}}\}$, we obtain that the power series
\[S_{L_{p_0}} := \sum_{l \in \Nn^r}a_{L_{p_0},l} t^l \in\Qq[[t_{1},\ldots,t_{r}]]\]
is rational, and is of the form (\ref{eqn:unifrat}) if we allow the $X_{l}$ to be quantifier-free $\Lrg(\Val_{K})$-definable sets (in fact, we can take $X_l$ just to be a single point).
\item Often in this kind of rationality theorem, we can take $q_{l} = 1$ for all $l$. There are two reasons why more complicated rational coefficients appear here. The first reason is to turn the $X_{l}$ into varieties instead of definable sets and the other reason is to get rid of the residual constant symbols that appear due to elimination of imaginaries.
\item Given uniformly rational power series $(S_{L_{p}}(t))_{L_p\in \mathfrak{F}}$, set $\varphi_{L_p}(s)= S_{L_p}(\card{\resf(L_p)}^{-s})$, where $s$ is a complex parameter.  Then $\varphi_{L_p}(s)$ has the form
\begin{equation}
\label{eqn:zetaform}
 \varphi_{L_p}(s)= \frac{\sum_{\card{l}\leq d} q_{l}\card{X_{l}({\rm res}(L_{p}))} \card{{\rm res}(L_{p})}^{-ls}}{\card{{\rm res}(L_{p})}^{n}
\prod_{j=1}^k(1-\card{{\rm res}(L_{p})}^{a_j-sb_{j}})},
\end{equation}
where the $X_l$, etc., are as in Theorem~\ref{thm:rat}.  It then follows by change of variable that for any $s_0\in \Zz$, the function $\varphi_{L_p}(s- s_0)$ (regarded as a function of $s$) also has the form (\ref{eqn:zetaform}).  (The only slight subtlety here is that the change of variable might lead to a factor of the form $\card{{\rm res}(L_{p})}^n$ in the denominator where $n<0$; but in this case, we can delete the factor and replace each $X_l$ with $X_l\times {\mathbb A}^n$.)
\item Our applications in Sections~\ref{sec:zetagp} and \ref{sec:twst} below use only the single-variable formulation of Theorem~\ref{thm:rat} (but see Remark~\ref{rem:multivariate}).
\end{thmenum}
\end{remark}

\begin{proof}[Theorem~\ref{thm:rat}]
By uniform EI (Corollary~\ref{cor:EIunif})---and the fact that elimination of imaginaries still holds after adding new constants for $K$---there exist integers $m_{1}$ and $m_{2}$, some $\mathcal{N}\subseteq \Nn_{>0}$ and some $\LGn(K)$-formula $\phi(x,w)$ such that for all $p\gg 0$, for all proper expansions to $\LGn$ of $L_{p}\in\mathfrak{F}_{p}$, $\phi$ defines a function $f_{L_{p}}' \colon R_{L_{p}} \to L_p^{m_{2}} \times \Latt[m_{1}](L_{p})$ such that for every $x,y \in R_{L_{p}}, xE_{L_{p}}y \iff f_{L_p}'(x)=f_{L_p}'(y)$. Let $f_{L_p}' = (f_{L_p}'',f_{L_p})$ where $f_{L_p}'': R_{L_{p}} \to L_p^{m_{2}} $ and $f_{L_p}: R_{L_{p}} \to \Latt[m_{1}](L_{p})$. For $l \in \Nn^r$, let $\mathcal{E}_{L_p,l} = \{f_{L_p}'(x): x \in R_{L_p,l}\}$ and $\mathcal{E}_{L_{p}} = \bigcup_{l}\mathcal{E}_{L_{p},l}$; so $\mathcal{E}_{L_p,l} \subseteq L_p^{m_{2}} \times \Latt[m_{1}](L_p)$ is finite, and it is the series $\sum_l \card{\mathcal{E} _{L_p,l}} t^l$ we wish to understand.   Let $\pi_{L_{p}}:   \mathcal{E}_{L_{p}} \to \Latt[m_{1}](L_p)$ be the projection, and let $F_{L_{p},l}= \pi_{L_p}(\mathcal{E}_{L_p,l})$.

It follows from Lemma~\ref{lem:rel alg bounded} and the fact that on the valued field sort the model-theoretic algebraic closure in $\ACVF$ coincides with the field-theoretic algebraic closure, that the size of the fiber $e_{p}(x): = \card{(\pi_{L_{p}})^{-1}(x)}$ is bounded by some positive integer $D$ uniformly for $p\gg 0$.
 We may thus partition $F_{L_p,l}$ into finitely many pieces $F^\nu_{L_{p},l}=\{x \in F_{L_p,l}\mid e_{p}(x) = \nu\}$;
 then \[\sum_l \card{\mathcal{E}_{L_p,l}} t^l = \sum_{\nu \leq D} \nu\sum_l \card{F^\nu_{L_{p} ,l}} t^l,\] so it suffices to prove that the series for $F^\nu_{p,l}$ has the form (\ref{eqn:unifrat}).

Fix $\nu$ and let $F_{L_p,l} = F^\nu_{L_{p},l}$;
 we need to retain only the information that $(F_{L_p,l})_{L_p\in\mathfrak{F}}$ is a family of finite subsets of $\Latt[m](L_p)$, uniformly $K$-definable in $\mathfrak{F}$. We can identify each element of $\Latt[m](L_{p})$
 with an element of $\GL{m}(L_p)/\GL{m}(\Val(L_p))$, i.e., with a left coset of $\GL{m}(\Val(L_p))$; let $G_{L_{p},l}$ be the union of these cosets.  By  Eqn.\ (\ref{eqn:ballmeas}),
we have
\[ \mu_{L_p}(G_{L_{p},l}) = \card{F_{L_p,l}}. \]
Thus  \[\sum_l \card{F_{L_p,l}} t^l  =    \sum_l \mu_{L_{p}}(G_{L_{p},l}) t^l  \in \Qq[[t_1,\ldots,t_r]].\]

We can apply \cite[Theorem\,1.1 and Theorem\,3.1]{DeLoICM} to these series to obtain uniform rationality. Note that, due to the constants added for elimination of imaginaries, we need parametric versions of these results (cf.\ \cite{CluLoeMot}). So we find $n,a_j,b_j$ as in the statement of Theorem~\ref{thm:rat}, and varieties $X_{l}$ over $\Val_{K}[y]$---where $y$ is a tuple of variables specialized in $\resf(L_{p})$ to any tuple $(k_{n}\mid n\in \mathcal{N})$ of unramified $n$-Galois uniformizers---such that \ref{eqn:unifrat} holds (we can take \(q_l = 1\) for now). Let now show that we can choose the \(X_l\) over \(\Val_{K}\) at the cost of making \(q_l\) nontrivial. Let
$$ C_{n}(L_{p}) = \{k_{n}\in\resf(L_{p})\mid \mbox{$k_{n}$ is the residue of an unramified $n$-Galois uniformizer}\}. $$
If $\resf(L_{p})[\omega_{n}]$ is of degree $d = d_{n,L_{p}}$ over $\resf(L_{p})$, then
\[\card{C_{n}(L_{p})} = \frac{\phi(n)(\card{\resf(L_{p})}^{d}-1)}{n},\]
where $\phi$ is the Euler totient function. Let $C = \prod_{n\in\mathcal{N}} C_{n}$ and for all $c\in C(L_{p})$, let $X_{c,l}(L_{p})$ be the $L_{p}$-points of the specialization of $X_{l}$ to $c$ and $Y_{l} := \coprod_{c\in C} X_{c,l}$. Then \(\card{Y_l(\resf(L_p))} = \card{C_n(L_p)}\card{X_l(\resf(L_p))}\). It follows that \[\card{X_{l}(\resf(L_{p}))} = \sum_{d|n}1_{d_{n,L_{p}} = d}\;\frac{-n}{\phi(n)}\;\frac{\card{Y_{l}(\resf(L_{p}))}}{1-\card{\resf(L_{p})}^{d}},\]
where $1_{d_{n,L_{p}}=d} = 1$ if $d_{n,L_{p}} = d$ and $0$ otherwise. Note that \(Y_l\) is an $\Lrg(\Val_{K})$-definable set and hence, replacing $\card{X_{l}(\resf(L_{p}))}$ with the RHS of the above equation, we obtain a rational function of the right form where the $X_{l}$ are $\Lrg(\Val_{K})$-definable, but, by \cite[Theorem\,2.1]{DeLoICM}, $X_{l}$ may be assumed to be a $\Val_{K}$-variety for $p\gg 0$.

For $L_{p}$ such that $p$ is too small, we can still prove the rationality of $S_{L_{p}}$ by the same argument using results for finite extensions of $p$-adic fields instead of those for ultraproducts: replace Corollary~\ref{cor:EIunif} with Theorem~\ref{thm:Qp}, Lemma~\ref{lem:rel alg bounded} with the proof of (i) (relative algebraic boundedness) in Section~\ref{sec:padics} and \cite[Theorem\,1.1]{DeLoICM} with \cite[Theorem\,1.5 and Theorem\,1.6.1]{DenMSRI}.
\end{proof}

\begin{remark}
\label{rem:growth}
It follows from the uniform formula Eqn.\ (\ref{eqn:unifrat}) we gave for $S_{L_{p}}$ in Theorem~\ref{thm:rat} that there exist \(c\in \Qq\) and \(n\in\Nn\) such that we have the following uniform growth estimate on $a_{L_{p},l}$: for all \(l\), all $p\gg 0$ and all $L_{p}\in \mathfrak{F}$,
\begin{equation}
\label{eqn:growth}
 a_{L_{p},l} \leq c\card{\resf(L_{p})}^{r\card{l}}.
\end{equation}
This estimate can be obtained by applying the uniform formula Eqn.\ (\ref{eqn:unifrat}) and using a polynomial upper bound on the number of \(\Ff_q\)-points of the varieties \(X_l\).

If $\mathfrak{F}_p$ is finite for all $p$ then Eqn.\ (\ref{eqn:growth}) holds for every $L_p\in \mathfrak{F}$.
\end{remark}

Below we consider uniformly $\emptyset$-definable families that arise in the following way.  Take ${\mathcal L}_p$ to be $\{\Qq_p\}$ for all $p$.  To simplify the notation in this case, we use subscripts $p$ instead of $\Qq_p$ (hence we write ${\mathcal D}_{p}$ and $S_p(t)$ below rather than ${\mathcal D}_{\Qq_p}$ and $S_{\Qq_p}(t)$).  Let ${\mathcal D}_{p} \subseteq \Qq_p^N$, let ${\mathcal E}_{p}$ be an equivalence relation on ${\mathcal D}_{p}$ and suppose $({\mathcal D}_{p})_{p \textrm{ prime}}$ and $({\mathcal E}_{p})_{p \textrm{ prime}}$ are uniformly $\emptyset$-definable in $p$.  Suppose that $f_{p,1},\ldots,f_{p,r} \colon {\mathcal D}_{p} \to \Qq_p \sminus\{0\}$ are uniformly $\emptyset$-definable functions such that for every $l \in \Zz^r$, the subset $\{x \in {\mathcal D}_{p} \mid |f_{p,i}(x)|=p^{-l_i}\}$ is a union of $\mathcal{E}_{p}$-equivalence classes.  Set $D_p= \{(x,|f_{p,1}(x)|,\ldots, |f_{p,r}(x)|) \mid x\in {\mathcal D}_p\}\subseteq \Qq_p^N\times \Zz^r$ and define an equivalence relation $E_p\subseteq D_p\times D_p$ by $(x,s_1,\ldots, s_r)E_p(x',s_1',\ldots, s_r')$ if $x{\mathcal E}_px'$ and $s_i= s_i'$ for all $i$.  Then we can regard $(D_{p})_{p \textrm{ prime}}$ as a uniformly $\emptyset$-definable family of sets and $(E_{p})_{p \textrm{ prime}}$ as a uniformly $\emptyset$-definable family of equivalence relations on $(D_{p})_{p \textrm{ prime}}$.

We now consider the abscissa of convergence of the zeta function in the one-variable case (under the assumption that ${\mathcal L}_p= \{\Qq_p\}$ for all $p$), and give a proof of Theorem~\ref{thm:abscissa_rat_intro}.  Recall that if $\zeta(s)= \sum_{n= 1}^\infty a_nn^{-s}$ is a zeta function then the abscissa of convergence $\alpha$ of $\zeta(s)$ is the infimum of the set of $s\in \Rr$ such that the series for $\zeta(s)$ is convergent.  Moreover, if $s\in {\mathbb C}$ then $\zeta(s)$ converges if ${\rm Re}(s)> \alpha$ and diverges if ${\rm Re}(s)< \alpha$.

We give a more precise statement of Theorem~\ref{thm:abscissa_rat_intro}.

\begin{thm}
\label{thm:abscissa_rat}
 Let ${\mathcal L}_p= \{\Qq_p\}$ for every prime $p$.  Assume the notation and hypotheses of Theorem~\ref{thm:rat} and define $\zeta_p(s)= S_p(p^{-s})$ (cf.\ Remark~\ref{rem:allprimes}.3).  Assume that the constant term of $\zeta_p(s)$ is 1 for all but finitely many primes and set $\zeta(s)= \prod_p \zeta_p(s)$.  Then the abscissa of convergence of $\zeta(s)$ is rational (or $-\infty$).
\end{thm}

Let $(\zeta_p(s))_{p \textrm{ prime}}$ be a family of zeta functions each of the form $\zeta_p(s)= \sum_{n=0}^\infty a_{p,n}p^{-ns}$.  Consider the formal product $\zeta(s)$ given by $\zeta(s)= \prod_p \zeta_p(s)$.
To ensure this makes sense, we assume that the constant term $a_{p,0}$ is 1 for all but finitely many primes.  To prove Theorem~\ref{thm:abscissa_rat}, we need to control the behavior of the $p$-local factors $\zeta_p(s)$.  Our proof is similar to parts of Avni's proof that the abscissa of convergence of the representation zeta function of an arithmetic lattice in a semisimple group is rational (see \cite[proof of Thm.~6.4]{Avni}, and cf.\ also \cite[Lem.~4.6 (1)]{duSGruAnn}), but the details are slightly different because we allow the coefficients $q_l$ in Eqn.\ (\ref{eqn:unifrat}) to be negative.

We need an estimate on the size of the varieties $X_l(\Ff_p)$ in Eqn.\ (\ref{eqn:unifrat}).  Recall the concept of an {\em Artin set} \cite[Defn.~4.6]{Avni}; as noted in {\em loc.\ cit.}, an infinite Artin set $A$ has positive analytic density, which implies that $\displaystyle \prod_{p\in A} \left(1+ \frac{1}{p}\right)$ diverges.

\begin{lem}
\label{lem:lang_weil}
 Let $X$ be a variety defined over $\Zz$.  Then there exist some partition of the set of primes into \(r\) disjoint Artin sets $A_1,\ldots, A_r$, some $c>0$ and, for all \(i\leq r\), some $(d_i,\mu_i)\in \Nn\times \Qq_{>0}$ such that for every prime \(p\), if \(p\in A_i\), then
 \begin{equation}
 \label{eqn:lang_weil}
  |X(\Ff_p)- \mu_i p^{d_i}|< cp^{d_i-1/2}.
 \end{equation}
\end{lem}

\begin{proof}
 This follows from \cite[Cor.~4.7]{Avni}, taking the parameter $n$ and the formula $\phi(x,y)$ of {\em loc.\ cit.} to be 0 and a formula $\phi(x)$ that defines $X$, respectively.  Note that the quantity $N_{d,\mu}$ in {\em loc.\ cit.}\ is 1 if Eqn.\ (\ref{eqn:lang_weil}) holds for a given $p\in A_i$, and 0 if it does not, so Eqn.\ (\ref{eqn:lang_weil}) holds for sufficiently large $p$.  By increasing $c$ if necessary, we can make Eqn.\ (\ref{eqn:lang_weil}) hold for all $p$.
\end{proof}

We recall two standard facts.
\begin{itemize}
 \item[(I)] If $(x_n)$ is a sequence of non-negative real numbers then $\prod_n (1+ x_n)$ converges if and only if $\sum_n x_n$ converges.
 \item[(II)] The abscissa of convergence of a finite product of zeta functions with non-negative coefficients is the maximum of the abscissae of convergence of the factors.
\end{itemize}

Let $A$ be a set of primes with positive analytic density (in particular, this implies that $A$ is infinite).  Let $t\in \Nn_{>0}$, let $d_1,\ldots, d_t\in \Zz$, let $e_1,\ldots, e_t$ be distinct positive integers and let $q_1,\ldots, q_t$ be nonzero real numbers.  Let $k\in \Nn$, let $n\in \Nn$ and let $a_1,\ldots, a_k\in \Zz$, $b_1,\ldots, b_k\in \Nn_{>0}$.  Let $u\in \Nn$ and let $g_1,\ldots, g_u$ be nonzero integers.  Let $\epsilon_1,\ldots, \epsilon_t> 0$, let $\mu_1,\ldots, \mu_t\geq 0$ and let $f_1,\ldots, f_t\colon A\ra \Rr$ such that $|f_i(p)|\leq \mu_ip^{d_i- \epsilon_i}$ for all $p\in A$.  Consider the $p$-local zeta function
\begin{equation}
\label{eqn:p-local}
 \displaystyle \zeta_p(s):= 1+ \frac{\sum_{l= 1}^t q_l(p^{d_l}+ f_l(p))p^{-e_ls}}{p^n\prod_{m= 1}^u (1- p^{g_m})\prod_{j= 1}^k (1- p^{a_j- b_js})}.
\end{equation}
We assume that the coefficients of $\zeta_p(s)$ (as a power series in $p^{-s}$) are non-negative.

We wish to determine the abscissa of convergence $\alpha$ of $\zeta(s):= \prod_{p\in A} \zeta_p(s)$.  For each $p\in A$, the poles of $\zeta_p(s)$ lie in the set $\displaystyle \left\{ \frac{a_j}{b_j}\mid 1\leq j\leq k\right\}$; but not every $\displaystyle \frac{a_j}{b_j}$ is necessarily a pole of $\zeta_p(s)$, since the numerator of the fraction on the RHS of Eqn.\ (\ref{eqn:p-local}) might have a zero at $\displaystyle \frac{a_j}{b_j}$.  Let $\displaystyle \Xi= \left\{j\mid 1\leq j\leq k, \frac{a_j}{b_j}\ \textrm{is a pole of $\zeta_p(s)$ for some $p\in A$}\right\}$.  Set $\displaystyle M_1= {\rm max}\left\{\frac{a_j}{b_j}\mid j\in \Xi\right\}$ (we take $M_1= -\infty$ if $\Xi$ is empty).

Given $s\in \Rr$ and $i\in \{1,\ldots, t\}$, we say that $s$ is {\em $i$-dominant} if $d_i- e_is> d_l- e_ls$ for all $l\neq i$.  If $s$ is not $i$-dominant for any $i$ then we say that $s$ is {\em critical}.  The set of critical points is finite (each critical point satisfies an equation of the form $d_l- e_ls= d_{l'}- e_{l'}s$ for some distinct $l$ and $l'$, and we assume that $e_l\neq e_{l'}$).

\begin{lem}
\label{lem:abscissa_rat}
 Let the notation be as above.  Then $\alpha$ is rational or $-\infty$.
\end{lem}

\begin{proof}
 If $t= 0$ in Eqn.\ (\ref{eqn:p-local}) then $\zeta_p(s)= 1$ for all $p\in A$, so $\alpha= -\infty$ and we are done.  Hence we can suppose that $t\geq 1$; in particular, $\zeta_p(s)$ is a strictly decreasing function of $s$ for $s> M_1$.  For any $s\in \Rr$, if $\zeta(s)$ converges then standard results on infinite products of Dirichlet series imply that each $\zeta_p(s)$ converges.  Hence $\alpha\geq M_1$.

 For $s\in \Rr$, set
$$ \beta(s)= {\rm max}_{1\leq l\leq t} (d_l- e_ls)- n-  \sum_{g_m> 0} g_m+ \sum_{a_j/b_j\geq s} (b_js- a_j). $$
 Then $\beta(s)$ is piecewise linear, so it is continuous.  We show that $\beta(s)$ is a strictly decreasing function of $s$ for $s> M_1$.  To see this, let $s\in \Rr$.  If $s$ is not critical then $s$ is $i$-dominant for some $i$, and it follows that there exists $E> 0$ such that
 $$ E^{-1}p^{d_i- e_is}\leq \left|\sum_{l= 1}^t q_l(p^{d_l}+ f_l(p))p^{-e_ls}\right|\leq Ep^{d_i- e_is} $$
 for all sufficiently large $p\in A$.
 Moreover, there exists $D> 0$ such that for all $j$ and all $p\in A$, if $\displaystyle s> \frac{a_j}{b_j}$ then $1- p^{a_j- b_js}> D$, while if $\displaystyle s< \frac{a_j}{b_j}$ then $|1- p^{a_j- b_js}|> Dp^{a_j- b_js}$.  It follows from the above discussion, the definition of $\zeta_p(s)$ and the bounds on the $f_l$ that for any $s\in \Rr$ such that $s> M_1$, $s$ is not critical and $\displaystyle s\neq \frac{a_j}{b_j}$ for any $1\leq j\leq k$, there exists $C\geq 1$ such that
\begin{equation}
\label{eqn:sandwich}
 C^{-1}p^{\beta(s)}\leq \zeta_p(s)- 1\leq Cp^{\beta(s)}
\end{equation}
\noindent for all sufficiently large $p\in A$.
Since each $\zeta_p(s)$ is strictly decreasing for $s> M_1$, $\beta(s)$ must therefore also be strictly decreasing for $s> M_1$, as claimed.  Hence there is at most one point $s_0> M_1$ such that $\beta(s_0)= -1$.  Set $M_2= s_0$ if this exists; otherwise set $M_2= -\infty$ (note that in the latter case, $\beta(s)< -1$ for all $s> M_1$, as $\lim_{s\to \infty} \beta(s)= -\infty$.)  We show that $\alpha= {\rm max}(M_1, M_2)$.

Let $s> M_1$ such that $s$ is not critical.  Suppose $s< M_2$.  Then $\beta(s)> -1$, so $\sum_{p\in A} C^{-1}p^{\beta(s)}$ diverges since $A$ has positive analytic density, so $\prod_{p\in A} (1+ C^{-1}p^{\beta(s)})$ diverges by Fact I.  Hence $\zeta(s)= \prod_{p\in A} \zeta_p(s)$ diverges by Eqn.\ (\ref{eqn:sandwich}) and the comparison test, and it follows that $s\leq \alpha$.  If $s> M_2$ then $\beta(s)< -1$, and a similar argument shows that $s\geq \alpha$.  We deduce that if $M_2\leq M_1$ then $\alpha= M_1$, and also that if $M_2> M_1$ then $\alpha= M_2$.  This completes the proof.
\end{proof}

\begin{proof} (of Theorems~\ref{thm:abscissa_rat} and \ref{thm:abscissa_rat_intro})
 By Theorem~\ref{thm:rat} and Remark~\ref{rem:allprimes}, Eqn.\ (\ref{eqn:unifrat}) holds for every prime $p$ and the definable sets $X_l$ are varieties over $\Zz$ for all but finitely many $p$.   Hence for each $p$, $\zeta_p(s)$ can be written as a rational function, where the numerator is a polynomial in $p^{-s}$ and the denominator is of the form $\prod_{j=1}^{k'} (1- p^{a_j- b_js})$ with each $b_j> 0$. (Here we have ordered the factors in the denominator of Eqn.\ (\ref{eqn:unifrat}) so that $b_1,\ldots, b_{k'}> 0$ and $b_{k'+1},\ldots, b_k= 0$ for some $0\leq k'\leq k$.)  This implies that the abscissa of convergence of each $\zeta_p(s)$ is rational.

 It now follows that by Fact (II), we can disregard finitely many primes: that is, it is enough to prove that $\prod_{p>p_0} \zeta_p(s)$ has rational abscissa of convergence for some prime $p_0$.  We can assume that $\zeta_p(s)$ has constant term 1 for every $p> p_0$.  The $\zeta_p(s)$ all have non-negative coefficients by construction.  Let $S^*_p(t)= S_p(t)- 1$ and let $\zeta^*_p(s)= \zeta_p(s)- 1= S^*_p(p^{-s})$.  Then $S^*_p(t)$ is the power series that arises from counting the equivalence classes of a uniformly $\emptyset$-definable family (in $\mathfrak{F}:= \bigcup_{p> p_0} \{\Qq_p\}$) of equivalence relations---just take the family of equivalence relations corresponding to $S_p(t)$ and remove the definable piece coming from $l=0$---so the power series $(S^*_p(t))_{p \textrm{ prime}}$ are uniformly rational for $p> p_0$ by Theorem~\ref{thm:rat} and Remark~\ref{rem:allprimes}.  Hence $S^*_p(t)$ is of the form given in Eqn.\ (\ref{eqn:unifrat}) for all $p> p_0$, with the sum in the numerator beginning at $l=1$ rather than $l=0$.  Explicitly, we have
 \begin{equation}
 \label{eqn:star_S}
  S^*_p(t)= \frac{\sum_{l=1}^d q_l|X_l(\Ff_p)|t^l}{p^n \prod_{j=1}^{k} (1- p^{a_j}t^{ b_j})},
 \end{equation}
 for all $p> p_0$, where the $X_l$, etc., are as in Theorem~\ref{thm:rat}.

 We now apply Lemma~\ref{lem:lang_weil} to the varieties $X_l$.  We can choose $r$, $A_1,\ldots, A_r$, $c, d_l, \mu_l$ such that Eqn.\ (\ref{eqn:lang_weil}) holds for each of the $X_l$ (note that complements, finite unions and finite intersections of Artin sets are Artin sets).  By increasing $p_0$ if necessary, we can assume that each $A_i$ is infinite and contains no primes less than or equal to $p_0$; in particular, each $A_i$ has positive analytic density.  It is enough by Fact (II) to show that $\displaystyle \prod_{p\in A_i}  \zeta_p(s)$ has rational abscissa of convergence for each $i$.  It follows from Eqns.~(\ref{eqn:lang_weil}) and (\ref{eqn:star_S}) that the hypotheses of Lemma~\ref{lem:abscissa_rat} are satisfied for $(\zeta_p(s))_{p\in A_i}$, so the desired result follows from Lemma~\ref{lem:abscissa_rat}.
\end{proof}

\begin{remark}
 As we discuss in Section~\ref{sec:zetagp} below, one-variable zeta functions that arise from cone integrals can be meromorphically continued beyond the abscissa of convergence, so one can apply Tauberian theorems and obtain more precise growth estimates than that provided by Eqn.\ (\ref{eqn:growth}): see \cite[Thm.~1.5]{duSGruAnn}.  In particular, this applies to the subgroup zeta functions that we discuss in Section~\ref{sec:zetagp} (see the discussion following \cite[Thm.~1.1]{duSGruAnn}).  Du Sautoy and Grunewald give a simple example of a zeta function that cannot be analytically continued beyond its abscissa of convergence (see \cite[Eqn.\ (1.3)]{duSGruAnn} and the discussion that surrounds it).  Hence one should not expect the stronger growth estimates to hold for zeta functions arising from an arbitrary uniformly $\emptyset$-definable equivalence relation.
\end{remark}

Let us conclude this section with a short aside on positive characteristic local fields by explaining how Theorem\,\ref{thm:rat} also yields transfer results between positive characteristic and mixed characteristic:

\begin{corollary}
\label{cor:transfer}
Let $\phi(x,y,l)$ be an $\LG$-formula where $l$ is a tuple of variables from the value group. The following are equivalent:
\begin{thmenum}
\item For all $p\gg 0$, the formula $\phi$ defines a family of finite equivalence relations $E_{p,l}$ on some set $D_{p,l}$ in $\Qq_{p}$;
\item For all $p\gg 0$, the formula $\phi$ defines a family of finite equivalence relations $E'_{p,l}$ on some set $D'_{p,l}$ in $\Ff_{p}((t))$.
\end{thmenum}
Moreover, whenever the above statements hold, there exists a prime $p_{0}$ such that for all $p\geq p_{0}$, the series $S_{p}(t) := \sum_{l \in \Nn^r}\card{D_{p,l}/E_{p,l}} t^l$ and $S'_{p}(t) := \sum_{l \in \Nn^r}\card{D'_{p,l}/E'_{p,l}} t^l$ are uniformly rational and $S_{p} = S'_{p}$.
\end{corollary}

\begin{proof}
This follows immediately from (the proof of) Theorem\,\ref{thm:rat} and the fact that for all non-principal ultrafilters $\mathcal{U}$ on the set of primes, $\prod_{p} \Qq_{p}/\mathcal{U}$ and $\prod_{p} \Ff_{p}((t))/\mathcal{U}$ are elementarily equivalent.
\end{proof}

\section{Zeta functions of groups}
\label{sec:zetagp}

We now consider some applications to some zeta functions that arise in group theory.  From now until the final part of Section~\ref{sec:twst} we take $\mathfrak{F}_p$ to be $\{\Qq_p\}$ for all $p$.  Most of the examples in this section come from the theory of subgroup growth of finitely generated nilpotent groups.  In Section~\ref{sec:twst} we consider the representation zeta function of finitely generated nilpotent groups.  We use Theorem~\ref{thm:rat} to prove uniform rationality of these zeta functions and Theorem~\ref{thm:abscissa_rat} to prove that the abscissa of convergence of the corresponding global zeta function is rational.  In the subgroup case this gives alternative proofs of results of \cite{GSS} and \cite{duSGruAnn}.

Throughout this section $\nilgp$ is a finitely generated nilpotent group.  For any $n \in\Nn$, the number $b_n$ of index $n$ subgroups of $\nilgp$ is finite (for background on subgroup growth, see \cite{LubSeg}).  The (global) subgroup zeta function of $\nilgp$ is defined by $\xi_\nilgp(s):= \sum_{n=1}^\infty b_nn^{-s}$ and the $p$-local subgroup zeta function by $\xi_{\nilgp,p}(s):= \sum_{n=0}^\infty b_{p^n}p^{-ns}$ (the symbol $\zeta$ is commonly used to denote the subgroup zeta function but we reserve this for the representation zeta function in Section~\ref{sec:twst}).  These expressions converge if $\mathrm{Re}(s)$ is large enough.  Grunewald, Segal and Smith observed in \cite{GSS} that Euler factorization holds: we have
$$ \xi_\nilgp(s)= \prod_p \xi_{\nilgp,p}(s), $$
where $p$ ranges over all primes.  Theorem~\ref{thm:ratsbgp} below (and \cite[Theorem~1]{GSS}) says that $\xi_{\nilgp,p}(s)$ is a rational function of $p^{-s}$.  Hence $\xi_\nilgp(s)$ enjoys many of the properties of the Riemann zeta function.

To understand the behavior of the global subgroup zeta function, one needs to study the behaviour of the rational function $\xi_{\nilgp,p}(s)$ as $p$ varies (cf.\ \cite{Avni}).  Du Sautoy and Grunewald introduced a class of $p$-adic integrals they called {\em cone integrals}.  They showed \cite[Theorem~1.3]{duSGruAnn} that if $\tau_p(s):= \sum_{n=0}^\infty b_{p,n}p^{-ns}$ is the zeta function arising from an Euler product of suitable cone integrals then $\tau_p(s)$ is uniformly rational for $p\gg 0$ (in the variable $t:= p^{-s}$) in the sense of Section~\ref{sec:rat}.  In fact, they proved a considerably stronger result  \cite[Theorem~1.4]{duSGruAnn} and deduced various analytic properties of $\tau(s)$ \cite[Theorem~1.5]{duSGruAnn}: for instance, they showed that $\tau(s)$ can be meromorphically continued a short distance to the left of its abscissa of convergence.  It follows from these results on cone integrals that $\xi_{\nilgp,p}(s)$ is uniformly rational for $p\gg 0$ \cite[Section~5]{duSGruAnn}.  For $\nilgp$ a finitely generated free nilpotent group of class 2, a stronger uniformity result holds: there is a polynomial $W(X,Y)\in \Qq[X,Y]$ such that $\xi_{\nilgp,p}(s)= W(p,p^{-s})$ for every prime $p$ \cite[Theorem~2]{GSS}.  Du Sautoy, however, has given an example showing that this stronger result does not hold for $\nilgp$ of arbitrary nilpotency class \cite{duSSubGp}.

Theorem~\ref{thm:ratsbgp} below deals with some variations on the subgroup zeta function.  In order to formulate the problem in terms of definable equivalence relations, we need to recall some facts about nilpotent pro-$p$ groups, including the notion of a good basis for a subgroup of a torsion-free nilpotent group \cite[Section~2]{GSS}; we will need these ideas in Section~\ref{sec:twst} as well.  We write $\widehat{\gp}_p$ for the pro-$p$ completion of a group $\gp$.  Let $\jinc \colon \nilgp \to \pgp$ be the canonical map.  Then $\pgp$ is finitely generated as a pro-$p$ group, so every finite-index subgroup of $\pgp$ is open (cf.\ \cite[Theorem~1.17]{DDMS}) and has $p$-power index (cf.\ \cite[Lemma~1.18]{DDMS}).  Since $\nilgp$ is finitely generated nilpotent, every subgroup of $p$-power index is open in the pro-$p$ topology on $\nilgp$;
in particular, there is a bijection $\hgp \mapsto \ovl{\jinc(\hgp)}$ between index $p^n$ subgroups of $\nilgp$ and index $p^n$ subgroups of $\pgp$, and $\ovl{\jinc(\hgp)} \iso \widehat{\hgp}_p$ (see \cite[Proposition~1.2]{GSS}).  For any $\hgp \unlhd \nilgp$ of index $p^n$, we have $\nilgp/\hgp \iso \pgp/\ovl{\jinc(\hgp)}$.

Let $\tfnilgp$ be a finitely generated torsion-free nilpotent group.  A {\em Mal'cev basis} is a tuple $a_1,\ldots,a_R$ of elements of $\tfnilgp$ such that any element of $\tfnilgp$ can be written uniquely in the form $a_1^{\lambda_1}\cdots a_R^{\lambda_R}$, where the $\lambda_i \in \Zz$.  We call the $\lambda_i$ {\em Mal'cev coordinates}.  Moreover, we require that group multiplication and inversion in $\tfnilgp$ are given by polynomials in the $\lambda_i$ with coefficients in $\Qq$, and likewise for the map $\tfnilgp \times \Zz \to \tfnilgp$, $(g,\lambda) \mapsto g^\lambda$.  We may regard the $a_i$ as elements of the pro-$p$ completion $\tfpgp$, and analogous statements hold, except that $\lambda$ and the Mal'cev coordinates $\lambda_i$ now belong to $\Zz_p$ (see  \cite[Section~2]{GSS}).  In particular, the map $\jinc \colon \tfnilgp \to \tfpgp$ is injective and we may identify $\tfpgp$ with $\Zz_p^R$.

Now let $\hgp$ be a finite-index subgroup of $\tfpgp$, of index $p^n$, say. In \cite{GSS}, a {\em good basis for $\hgp$} is defined as an $R$-tuple $h_1,\ldots,h_R \in \hgp$ such that every element of $\hgp$ can be written uniquely in the form $h_1^{\lambda_1}\cdots h_R^{\lambda_R}$ ($\lambda_i \in \Zz_p$), and satisfying an extra property which does not concern us here.  We say that $h_1,\ldots,h_R \in \tfpgp$ is a {\em good basis} if it is a good basis for some finite-index subgroup $\hgp$ of $\tfpgp$.  For each $i$, we can write
\begin{equation}
\label{eqn:sbgpcoords}
 h_i=a_1^{\lambda_{i1}}\cdots a_R^{\lambda_{iR}}
\end{equation}
and we recover $\indx{\tfpgp}{\hgp} = p^n$ from the formula
\begin{equation}
\label{eqn:indx}
|\lambda_{11} \lambda_{22}\cdots \lambda_{RR}|=p^{-n}.
\end{equation}
Any finite-index subgroup of $\tfpgp$ admits a good basis.  Often we will identify a good basis $h_1,\ldots,h_R$ with the $R^2$-tuple of coordinates $(\lambda_{ij})$.

\begin{proposition}
\label{propn:gddef}
Let ${\mathcal D}_p\subseteq \Zz_p^{R^2}$ be the set of good bases $(\lambda_{ij})$ of $\tfpgp$.  Then the sets $({\mathcal D}_p)_{p \textrm{ prime}}$ are uniformly $\emptyset$-definable in $p$ (in the structure \(\Qq_p\)).
\end{proposition}

\begin{proof}
 This follows from the proof of \cite[Lemma~2.3]{GSS}.
\end{proof}

For each non-negative $n$ consider the following:\\
(a) the number of index $p^n$ subgroups of $\tfnilgp$;\\
(b) the number of normal index $p^n$ subgroups of $\tfnilgp$;\\
(c) the number of index $p^n$ subgroups $A$ of $\tfnilgp$ such that $\widehat{A}_p \iso \tfpgp$;\\
(d) the number of conjugacy classes of index $p^n$ subgroups of $\tfnilgp$;\\
(e) the number of equivalence classes of index $p^n$ subgroups of $\tfnilgp$, where we define $A \sim B$ if $\widehat{A}_p \iso \widehat{B}_p$.\\
The rationality of $\sum_{n=0}^\infty b_{p,n}t^n$ in (a)--(d) of the following result are due to Grunewald, Segal, and Smith \cite[Theorem~1]{GSS}; for uniformity statements and the rationality of the abscissa of convergence in (a)--(d), see \cite[Section~1]{duSGruAnn} and the start of this section.  Here we give a different proof.  Observe that Theorem~\ref{thm:ratsbgp} for case (e) is new; here the equivalence relation does not arise from any obvious group action, and Theorem~\ref{thm:rat} gives a genuinely new way of proving uniform rationality.  This illustrates the robustness of our methods, which are not sensitive to the precise details of how the objects to be counted are interpreted.

\begin{theorem}
\label{thm:ratsbgp}
 Let  $b_{p,n}$ be as described in any of (a)--(e) above.  Set $S_p(t)= \sum_{n=0}^\infty b_{p,n}t^n$.  Then the power series $(S_p(t))_{p \textrm{ prime}}$ are uniformly rational.  Moreover, the zeta function $\xi(s):= \prod_p \sum_{n=0}^\infty b_{p,n}p^{-ns}$ has rational abscissa of convergence.
\end{theorem}

\begin{proof}
 Clearly $b_{p,0}= 1$ for all $p$, so rationality of the abscissa of convergence of $\xi(s)$ will follows from Theorem~\ref{thm:abscissa_rat_intro} once we have proved the other assertions of Theorem~\ref{thm:ratsbgp}.  To prove the rest of the theorem, we show how to interpret the objects that we are counting in a uniformly $\emptyset$-definable way, then apply Theorem~\ref{thm:rat}.  Consider case (a).  Let ${\mathcal D}_p$ be as in Proposition~\ref{propn:gddef}.  Define $f_p \colon {\mathcal D}_p \to \Zz_p$ by $f_p(\lambda_{ij})=\lambda_{11} \cdots \lambda_{RR}$; note that the functions $(f_p)_{p \textrm{ prime}}$ are uniformly $\emptyset$-definable in $p$.  Define an equivalence relation ${\mathcal E}_p$ on ${\mathcal D}_p$ as follows: two $R$-tuples $(\lambda_{ij})$, $(\mu_{ij})$, representing good bases $h_1,\ldots,h_R$ and $k_1,\ldots,k_R$ for subgroups $\hgp$, $\kgp$ respectively, are equivalent if and only if $\hgp=\kgp$.

 Now the equivalence relations $({\mathcal E}_p)_{p \textrm{ prime}}$ are uniformly $\emptyset$-definable in $p$: each ${\mathcal E}_p$ is the subset of ${\mathcal D}_p \times {\mathcal D}_p$ given by the conjunction for $1 \leq i,j \leq R$ of the formulae
 $$ (\exists \sigma^{(i)}_1,\ldots,\sigma^{(i)}_R \in \Zz_p)\; k_i=h_1^{\sigma^{(i)}_1} \cdots h_R^{\sigma^{(i)}_R}$$
 and
 $$ (\exists \tau^{(j)}_1,\ldots,\tau^{(j)}_R \in \Zz_p)\; h_j=k_1^{\tau^{(j)}_1} \cdots k_R^{\tau^{(j)}_R},$$
 and these become polynomial equations independent of $p$ over $\Qq$ in the $\lambda_{ij}$, the $\mu_{ij}$, the $\sigma_i$ and the $\tau_j$ when we write the $h_i$ and $k_j$ in terms of their Mal'cev coordinates (Eqn.\ (\ref{eqn:sbgpcoords})).

 Construct $D_p$ and $E_p$ from ${\mathcal E}_p$, ${\mathcal D}_p$ and $f_p$ as in the paragraph after Eqn.\ (\ref{eqn:growth}).  Using Eqn.\ (\ref{eqn:indx}), we see that for each $n \in \Nn$, $D_{p,n}/E_{p,n}$ consists of precisely $b_{p,n}$ equivalence classes.  We now deduce the rationality and uniform rationality assertions from Theorem~\ref{thm:rat} (taking $\mathfrak{F}_p= \{\Qq_p\}$ for all $p$) and Remark~\ref{rem:allprimes}.

 The proofs in cases (b)--(e) are similar, modifying the definitions of ${\mathcal D}_p$ and ${\mathcal E}_p$ appropriately.  For example, in (b) we replace ${\mathcal D}_p$ with the set ${\mathcal D}_p^\unlhd$ of tuples $(\lambda_{ij})$ that define a normal finite-index subgroup $\hgp$; a tuple $(\lambda_{ij})$ corresponding to a finite-index subgroup $\hgp$ belongs to ${\mathcal D}_p^\unlhd$ if and only if it satisfies the formula
 $$ (\forall g \in \tfpgp) (\forall h \in \hgp) (\exists \nu_1,\ldots,\nu_R \in \Zz_p)\; ghg^{-1} = h_1^{\nu_1} \cdots h_R^{\nu_R}, $$
 which is made up of polynomial equations independent of $p$ over $\Qq$ in the $\nu_i$, the $\lambda_{ij}$ and the Mal'cev coordinates of $g$ and $h$.  In case (d), the equivalence relation is the subset of ${\mathcal D}_p \times {\mathcal D}_p$ given by the formula:\smallskip\\
 there exists $g \in \tfpgp$, there exist $\sigma_i^{(j)}, \tau_i^{(j)} \in \Zz_p$ for $1 \leq j \leq R$ such that $gh_jg^{-1} = k_1^{\sigma_1^{(j)}} \cdots k_R^{\sigma_R^{(j)}}$ and $g^{-1}k_jg = h_1^{\tau_1^{(j)}} \cdots h_R^{\tau_R^{(j)}}$ for $1 \leq j \leq R$.\smallskip\\
 This is made up of polynomial equations independent of $p$ over $\Qq$ in the Mal'cev coordinates of $g$ and of the $h_i$ and the $k_i$.  In cases (c) and (e), we can express the isomorphism condition in terms of polynomials in the Mal'cev coordinates; compare the proof of Proposition~\ref{propn:enumeration} below.
\end{proof}

\begin{remark}
\label{rem:torsion}
Du Sautoy and Grunewald prove that Theorem~\ref{thm:ratsbgp} (a) and (b) actually hold for an arbitrary finitely generated nilpotent group $\nilgp$, possibly with torsion.  To prove this in our setting, write $\nilgp$ as a quotient $\tfnilgp/\kr$ of a finitely generated torsion-free nilpotent group $\tfnilgp$.  Theorem~\ref{thm:ratsbgp} now follow for cases (a)--(e) from our arguments above with suitable modifications: for example, for case (a), we count not all index $p^n$ subgroups of $\tfnilgp$, but only the ones that contain $\kr$.  For details, compare the argument of the last two paragraphs of Lemma~\ref{lem:induction}.
\end{remark}

The proof for case (d) of Theorem~\ref{thm:ratsbgp} is not given explicitly in \cite{GSS}, but the appropriate definable integral can be constructed using the methods in the proof of \cite[Theorem~1.2]{duSConj}; what makes this work is that the equivalence classes are the orbits of a group action.  The language of \cite{duSConj} contains symbols for analytic functions, but our methods still apply there because we can use the results of Cluckers from the Appendix, which do hold in the analytic setting.

Here is another application, to the problem of counting finite $p$-groups.
\begin{proposition}[du Sautoy]
\label{propn:enumeration}
 Fix positive integers $c,d$.  Let $c_{p,n}$ be the number of finite $p$-groups of order $p^n$ and nilpotency class at most $c$, generated by at most $d$ elements.  Set $S_p(t)= \sum_{n=0}^\infty c_{p,n}t^n$.  Then the power series $(S_p(t))_{p \textrm{ prime}}$ are uniformly rational.  Moreover, the zeta function $\chi(s):= \prod_p \sum_{n=0}^\infty c_{p,n}p^{-ns}$ has rational abscissa of convergence.
\end{proposition}

\begin{proof}
 As in Theorem~\ref{thm:ratsbgp}, the rationality of the abscissa of convergence will follow from Theorem~\ref{thm:abscissa_rat_intro}, and to prove the rest it is enough to interpret the objects we are counting in a uniformly $\emptyset$-definable way.  Let $\tfnilgp$ be the free nilpotent group of class $c$ on $d$ generators (note that $\tfnilgp$ is torsion-free).  Any finite $p$-group of order $p^n$ and nilpotency class at most $c$ and generated by at most $d$ elements is a quotient of $\tfpgp$ by some normal subgroup of index $p^n$.  Let ${\mathcal D}_p^\unlhd$ and $f_p$ be as in the proof of Theorem~\ref{thm:ratsbgp}.  Define an equivalence relation ${\mathcal E}_p$ on ${\mathcal D}_p^\unlhd$ as follows: two $R$-tuples $(\lambda_{ij})$, $(\mu_{ij})$, representing good bases $h_1,\ldots,h_R$ and $k_1,\ldots,k_R$ for subgroups $\hgp$, $\kgp$ respectively, are equivalent if and only if $\tfpgp/\hgp \iso \tfpgp/\kgp$.

 The result will follow as in Theorem~\ref{thm:ratsbgp} if we can show that the equivalence relations $({\mathcal E}_p)_{p \textrm{ prime}}$ are uniformly $\emptyset$-definable in $p$.  Let $a_1,\ldots, a_R$ be the Mal'cev basis of $\tfpgp$, as before.  We claim that ${\mathcal E}_p \subseteq {\mathcal D}_p^\unlhd \times {\mathcal D}_p^\unlhd$ is given by the following conditions:
 \begin{equation}
 \label{eqn:isoquot0}
  |f_p(\lambda_{ij})|=|f_p(\mu_{ij})|,
 \end{equation}
 \begin{equation}
\label{eqn:isoquot1}
 (\exists b_1,\ldots,b_r \in \tfpgp)(\forall \nu_1,\ldots,\nu_r \in \Zz_p)\; a_1^{\nu_1} \cdots a_R^{\nu_R} \in \hgp \iff b_1^{\nu_1} \cdots b_R^{\nu_R} \in \kgp
 \end{equation}
and
$$ (\forall \sigma_1,\ldots,\sigma_r, \tau_1,\ldots, \tau_r \in \Zz_p) (\exists  \nu_1,\ldots, \nu_r \in \Zz_p) $$
\begin{equation}
\label{eqn:isoquot2}
 (a_1^{\sigma_1}\cdots a_R^{\sigma_R} a_1^{\tau_1}\cdots a_R^{\tau_R}= a_1^{\nu_1}\cdots a_R^{\nu_R})\wedge (b_1^{\sigma_1}\cdots b_R^{\sigma_R} b_1^{\tau_1}\cdots b_R^{\tau_R}\in b_1^{\nu_1}\cdots b_R^{\nu_R} K).
\end{equation}
 To prove this, suppose Eqns.\ (\ref{eqn:isoquot0}), (\ref{eqn:isoquot1}) and (\ref{eqn:isoquot2}) hold.  Then $\indx{\tfpgp}{\hgp}=\indx{\tfpgp}{\kgp}$ and the map $a_i\hgp \mapsto b_i\kgp$ defines an isomorphism from $\tfpgp/\hgp$ onto $\tfpgp/\kgp$.  Conversely, if $g$ is an isomorphism from $\tfpgp/\hgp$ onto $\tfpgp/\kgp$ then $\indx{\tfpgp}{\hgp}=\indx{\tfpgp}{\kgp}$, so $|f_p(\lambda_{ij})|=|f_p(\mu_{ij})|$.  Moreover, we can choose $b_i \in \tfpgp$ such that $g(a_i\hgp) = b_i\kgp$ for $1 \leq i \leq R$.  Then for all $\nu_1,\ldots,\nu_r \in \Zz$ we have
 $$ (*)\ \ \ \  a_1^{\nu_1} \cdots a_R^{\nu_R} \in \hgp \iff b_1^{\nu_1} \cdots b_R^{\nu_R} \in \kgp $$
 and for all $\sigma_1,\ldots,\sigma_r, \tau_1,\ldots, \tau_r \in \Zz$ there exist $ \nu_1,\ldots, \nu_r \in \Zz$ such that
  $$ (**)\ \ \ \  (a_1^{\sigma_1}\cdots a_R^{\sigma_R} a_1^{\tau_1}\cdots a_R^{\tau_R}= a_1^{\nu_1}\cdots a_R^{\nu_R})\wedge (b_1^{\sigma_1}\cdots b_R^{\sigma_R} b_1^{\tau_1}\cdots b_R^{\tau_R}\in b_1^{\nu_1}\cdots b_R^{\nu_R}K); $$
 since $\hgp,\kgp$ are closed and the group operations are continuous, $(*)$ and $(**)$ hold with $\Zz$ replaced by $\Zz_p$.  This proves the claim.  The formulae above involve only the functions $(f_p)_{p \textrm{ prime}}$---which are uniformly $\emptyset$-definable in $p$---and polynomials independent of $p$ over $\Qq$ in the Mal'cev coordinates, so the equivalence relations $({\mathcal E}_p)_{p \textrm{ prime}}$ are uniformly $\emptyset$-definable in $p$, as required.
\end{proof}

Du Sautoy's proof \cite[Theorem~2.2]{duSpGpPrelim}, \cite[Theorems~1.6 and 1.8]{duSpGpFull} uses the fact that an isomorphism $\tfpgp/\hgp \to \tfpgp/\kgp$ lifts to an automorphism of $\tfpgp$, which implies that the equivalence relation ${\mathcal E}_p$ arises from the action of the group $\aut(\tfpgp)$, a compact $p$-adic analytic group.  This allows one to express the power series $\sum_{n=0}^\infty c_{p,n}t^n$ as a cone integral, from which uniform rationality follows (see the start of this section).  Our proof is simpler in its algebraic input,
as elimination of imaginaries allows us to use  less information about ${\mathcal E}_p$.

\begin{remark}
 Let $\nilgp$ be a finitely generated nilpotent group and let $c_{p,n}$ be the number of isomorphism classes of quotients of $\nilgp$ of order $p^n$.  Then the power series $(\sum_{n=0}^\infty c_{p,n}t^n)_{p \textrm{ prime}}$ are uniformly rational.  If $\nilgp$ is torsion-free then this follows immediately from the proof of Proposition~\ref{propn:enumeration}.  If $\nilgp$ has torsion then we write $\nilgp$ as a quotient $\tfnilgp/\kr$ of a finitely generated torsion-free nilpotent group $\tfnilgp$ and modify the proof of Proposition~\ref{propn:enumeration} accordingly (cf.\ Remark~\ref{rem:torsion}).
\end{remark}

\section{Twist isoclasses of characters of nilpotent groups}
\label{sec:twst}

By a \emph{representation} of a group $\gp$ we shall mean a finite-dimensional complex representation, and by a \emph{character} of $\gp$ we shall mean the character of such a representation.  A character is said to be \emph{linear} if its degree is one.  We write $\ip{\,}{\;\!}{\gp}$ for the usual inner product of characters of $\gp$.  If $\chi$ is linear then we have
\begin{equation}
\label{eqn:twstinvc}
 \ip{\chi\sigma_1}{\chi\sigma_2}{\gp} = \ip{\sigma_1}{\sigma_2}{\gp}
\end{equation}
for all characters $\sigma_1$ and $\sigma_2$.  If $\gp' \leq \gp$ has finite index then we write $\Ind{\gp}{\gp'}{\cdot}$ and $\Res{\gp}{\gp'}{\cdot}$ for the induced character and restriction of a character respectively.  For background on representation theory, see \cite{CurRei}.  Below when we apply results from the representation theory of finite groups to representations of an infinite group, the representations concerned always factor through finite quotients.

We denote the set of irreducible $n$-dimensional characters of $\gp$ by $\irr{n}{\gp}$.  If $N \unlhd \gp$ then we say the character $\chi$ of an irreducible representation $\rho$ factors through $\gp/N$ if $\rho$ factors through $\ggp/N$ (this depends only on $\chi$, not on $\rho$).
\begin{notation}
 We say a character $\sigma$ of $\gp$ is {\em admissible} if $\sigma$ factors through a finite quotient of $\gp$.  If $p$ is prime then we say $\sigma$ is {\em $p$-admissible} if $\sigma$ factors through a finite $p$-group quotient of $\gp$.  We write $\cirr{n}{\gp}$ ($\pirr{n}{\gp}$) for the set of admissible ($p$-admissible) characters in $\irr{n}{\gp}$.  Note that $\pirr{n}{\gp}$ is empty if $n$ is not a $p$-power \cite[(9.3.2)~Proposition]{CurRei}.
\end{notation}

Given $\sigma_1,\sigma_2 \in \irr{n}{\gp}$, we follow \cite{LubMag} and say that $\sigma_1$ and $\sigma_2$ are {\em twist-equivalent} if $\sigma_1=\chi\sigma_2$ for some linear character $\chi$ of $\gp$.  Clearly this defines an equivalence relation on $\irr{n}{\gp}$; we call the equivalence classes {\em twist isoclasses}.

\begin{observation}
\label{obs:twstquot}
 Let $\sigma_1,\sigma_2$ be two irreducible degree $n$ characters of $\gp$ that are twist-equivalent: say $\sigma_2=\chi\sigma_1$.  If $N \unlhd \gp$ such that $\sigma_1,\sigma_2$ both factor through $\gp/N$, then $\chi$ also factors through $\gp/N$.
\end{observation}
If $N_1,N_2 \unlhd \gp$ have finite ($p$-power) index then $N_1 \cap N_2$ also has finite ($p$-power) index.  This implies that when we are working with twist isoclasses in $\cirr{n}{\gp}$ ($\pirr{n}{\gp}$), we need only consider twisting by admissible ($p$-admissible) linear characters.

Fix a finitely generated nilpotent group $\nilgp$.  The set $\irr{n}{\nilgp}$ can be given the structure of a quasi-affine complex algebraic variety.  Lubotzky and Magid analyzed the geometry of this variety and proved the following result \cite[Theorem~6.6]{LubMag}.
\begin{theorem}
\label{thm:fintwst}
 There exists a finite quotient $\nilgp(n)$ of $\nilgp$ such that every irreducible $n$-dimensional representation of $\nilgp$ factors through $\nilgp(n)$ up to twisting.  In particular, there are only finitely many twist isoclasses of irreducible $n$-dimensional characters.
\end{theorem}
Thus the number of degree $n$ twist isoclasses is a finite number $a_n$.

\begin{definition}
 We define the (global) representation zeta function $\zeta_\nilgp(s)$ by $\zeta_\nilgp(s)= \sum_{n= 1}^\infty a_nn^{-s}$ and the $p$-local representation zeta function $\twzeta{\nilgp}{p}{s}$ by $\twzeta{\nilgp}{p}{s}= \sum_{n= 0}^\infty a_{p^n}p^{-ns}$.
\end{definition}

\noindent It is shown in \cite[Lemma~2.1]{StaVol} that $\zeta_\nilgp(s)$ converges on some right-half plane.  Voll noted \cite[Section~3.2.1]{VollNewcomer} that $\zeta_\nilgp(s)$ has an Euler factorization

$$ \zeta_\nilgp(s)= \prod_{p}\twzeta{\nilgp}{p}{s} $$

\noindent for any finitely generated nilpotent group (cf.\ the proof of Lemma~\ref{lem:completion}).

We now turn to the proof of Theorem~\ref{thm:rattwst_intro}.  Clearly $a_1= 1$ by construction, so the rationality of the abscissa of convergence of $\zeta_\nilgp(s)$ will follow as usual from Theorem~\ref{thm:abscissa_rat_intro}.   To prove the rest of Theorem~\ref{thm:rattwst_intro}, we show how to interpret twist isoclasses in a uniformly $\emptyset$-definable way.  The equivalence relation in the parametrization is not simply the relation of twist-equivalence, which arises from the action of a group---the group of linear characters of $\nilgp$---but a more complicated equivalence relation.

The correspondence between index $p^n$ subgroups of $\nilgp$ and index $p^n$ subgroups of $\pgp$ gives a canonical bijection between $\pirr{p^n}{\nilgp}$ and $\pirr{p^n}{\pgp}$, and it is clear that this respects twisting by $p$-admissible characters.
\begin{lemma}
\label{lem:completion}
 For every non-negative integer $n$, there is a bijective correspondence between the sets $\irr{p^n}{\nilgp}/{(\mathrm{twisting})}$ and $\pirr{p^n}{\pgp}/{(\mathrm{twisting})}$.
\end{lemma}
\begin{proof}
 It suffices to show that given any $\sigma \in \irr{p^n}{\nilgp}$, some twist of $\sigma$ factors through a finite $p$-group quotient of $\nilgp$.  By Theorem~\ref{thm:fintwst}, we can assume that $\sigma$ factors through some finite quotient $\fingp$ of $\nilgp$. Let us also denote by \(\sigma\) the corresponding character of \(\fingp\). Then $\fingp$, being a finite nilpotent group, is the direct product of its Sylow $l$-subgroups $\fingp_l$, where $l$ ranges over all the primes dividing $|\fingp|$.  Moreover \cite[Theorem~10.33]{CurRei}, $\sigma$ is a product of irreducible characters $\sigma_l$, where each $\sigma_l$ is a character of $\fingp_l$.  Since the degree of an irreducible character of a finite group divides the order of the group \cite[Proposition~9.3.2]{CurRei}, all of the $\sigma_l$ for $l \neq p$ are linear.  We may therefore twist $\sigma$ by a linear character of $\fingp$ to obtain a character that kills $\fingp_l$ for $l \neq p$, and this linear character is admissible by Observation~\ref{obs:twstquot}.  The new character factors through $\fingp_p$, and we are done.
\end{proof}

The key idea is that finite $p$-groups are monomial: that is, every irreducible character is induced from a linear character of some subgroup.  We parametrize $p$-admissible irreducible characters of $\pgp$ by certain pairs $(\hgp,\chi)$, where $\hgp$ is a finite-index subgroup of $\pgp$ and $\chi$ is a $p$-admissible linear character of $\hgp$: to a pair we associate the induced character $\Ind{\pgp}{\hgp}{\chi}$.  We can parametrize these pairs using the theory of good bases for subgroups of $\pgp$, and this description is well-behaved with respect to twisting.  Two distinct pairs $(\hgp,\chi)$ and $(\hgp',\chi')$ may give the same induced character; this gives rise to a definable equivalence relation on the set of pairs.

If $\psi$ is a character of $\hgp$ and $g\in \pgp$ then we denote by $g.\psi$ the character of $g.\hgp:=g\hgp g^{-1}$ defined by $(g.\psi)(ghg^{-1})=\psi(h)$.
\begin{lemma}
\label{lem:induction}
 (a) Let $\sigma \in \pirr{p^n}{\pgp}$.  Then there exists $\hgp \leq \pgp$ such that $\indx{\pgp}{\hgp}=p^n$, together with a $p$-admissible linear character $\chi$ of $\hgp$ such that $\sigma=\Ind{\pgp}{\hgp}{\chi}$.\\
 (b) Let $\hgp$ be a $p$-power index subgroup of $\pgp$ and let $\chi$ be a $p$-admissible linear character of $\hgp$.  Then $\Ind{\pgp}{\hgp}{\chi}$ is a $p$-admissible character of $\pgp$, and $\Ind{\pgp}{\hgp}{\chi}$ is irreducible if and only if for all $g \in \pgp \sminus \hgp$, $\Res{g.\hgp}{g.\hgp \cap \hgp}{g.\chi} \neq \Res{\hgp}{g.\hgp \cap \hgp}{\chi}$.  Moreover, if $\psi$ is a $p$-admissible linear character of $\pgp$ and $\Ind{\pgp}{\hgp}{\chi}$ is irreducible then $\Ind{\pgp}{\hgp}{\left(\left(\Res{\pgp}{\hgp}{\psi}\right)\chi \right)}=\psi\;\Ind{\pgp}{\hgp}{\chi}$.\\
 (c) Let $\hgp,\hgp' \leq \pgp$ have index $p^n$, and let $\chi,\chi'$ be $p$-admissible linear characters of $\hgp,\hgp'$ respectively such that $\Ind{\pgp}{\hgp}{\chi}$ and $\Ind{\pgp}{\hgp'}{\chi'}$ are irreducible.  Then $\Ind{\pgp}{\hgp}{\chi}=\Ind{\pgp}{\hgp'}{\chi'}$ if and only if there exists $g \in \pgp$ such that $\Res{g.\hgp}{g.\hgp \cap \hgp'}{g.\chi}=\Res{\hgp'}{g.\hgp \cap \hgp'}{\chi'}$.
\end{lemma}
\begin{proof}
 (a) Since $\sigma$ is $p$-admissible, it factors through some finite $p$-group $\fingp$.  Since finite $p$-groups are monomial \cite[Theorem~11.3]{CurRei}, there exist $\lgp \leq \fingp$ of index $p^n$ and a linear character $\chi$ of $\lgp$ such that $\sigma$---regarded as a character of $\fingp$---equals $\Ind{\fingp}{\lgp}{\chi}$.  Let $\hgp$ be the pre-image of $\lgp$ under the canonical projection $\pgp \to \fingp$.  Regarding $\chi$ as a character of $\hgp$, it is easily checked that $\indx{\pgp}{\hgp}=p^n$ and $\sigma=\Ind{\pgp}{\hgp}{\chi}$.\\
 (b) Since $\chi$ is $p$-admissible, the kernel $\kgp$ of $\chi$ has $p$-power index in $\pgp$, so $\kgp$ contains a $p$-power index subgroup $N$ such that $N \unlhd \pgp$.  Clearly $N \leq \ker{(\Ind{\pgp}{\hgp}{\chi})}$, so $\Ind{\pgp}{\hgp}{\chi}$ is $p$-admissible.  The irreducibility criterion follows immediately from \cite[Theorem~10.25]{CurRei}.  By Frobenius reciprocity,
 \begin{eqnarray*}
   &   & \ip{\Ind{\pgp}{\hgp}{\left(\left(\Res{\pgp}{\hgp}{\psi}\right)\chi\right)}}{\psi \;\Ind{\pgp}{\hgp}{\chi}}{\pgp} \\
   & = & \ip{\left(\Res{\pgp}{\hgp}{\psi}\right)\chi}{\Res{\pgp}{\hgp}{\left(\psi\; \Ind{\pgp}{\hgp}{\chi}\right)}}{\hgp} \\
   & = & \ip{\left(\Res{\pgp}{\hgp}{\psi}\right)\chi}{\left(\Res{\pgp}{\hgp}{\psi}\right){\Res{\pgp}{\hgp} {\left(\Ind{\pgp}{\hgp}{\chi}\right)}}}{\hgp} \\
   & = & \ip{\chi}{{\Res{\pgp}{\hgp} {\left(\Ind{\pgp}{\hgp}{\chi}\right)}}}{\hgp}\ \mbox{by Eqn.\ (\ref{eqn:twstinvc})} \\
   & = & \ip{\Ind{\pgp}{\hgp}{\chi}}{\Ind{\pgp}{\hgp}{\chi}}{\pgp} \\
   & = & 1.
 \end{eqnarray*}
Now $\psi\;\Ind{\pgp}{\hgp}{\chi}$ is irreducible, because $\Ind{\pgp}{\hgp}{\chi}$ is, and the degrees of $\Ind{\pgp}{\hgp}{\left(\left(\Res{\pgp}{\hgp}{\psi}\right)\chi\right)}$ and $\psi\;\Ind{\pgp}{\hgp}{\chi}$ are equal.  We deduce that $\psi\;\Ind{\pgp}{\hgp}{\chi}= \Ind{\pgp}{\hgp}{\left(\left(\Res{\pgp}{\hgp}{\psi}\right)\chi\right)}$.\\
(c) The Mackey Subgroup Theorem~\cite[Theorem~10.13]{CurRei} gives
 \begin{equation}
 \Res{\pgp}{\hgp'}{\left(\Ind{\pgp}{\hgp}{\chi}\right)}= \sum_{\ovl{g} \in \hgp'\backslash \pgp/\hgp} \Ind{\hgp'}{g.\hgp \cap \hgp'}{\left(\Res{g.\hgp}{g.\hgp \cap\hgp'}{\;g.\chi}\right)}.
 \end{equation}
Here the sum is over a set of double coset representatives $g$ for $\hgp'\backslash \pgp/\hgp$ (the characters on the RHS of the formula are independent of choice of representative).  Since $\Ind{\pgp}{\hgp}{\chi}$ and $\Ind{\pgp}{\hgp'}{\chi'}$ are irreducible, they are distinct if and only if their inner product is zero.  We have
 \begin{eqnarray*}
   & & \ip{\Ind{\pgp}{\hgp}{\chi}}{\Ind{\pgp}{\hgp'}{\chi'}}{\pgp} \\
   & = & \ip{\Res{\pgp}{\hgp'}{(\Ind{\pgp}{\hgp}{\chi})}}{\chi'}{\hgp'}\ \mbox{by Frobenius reciprocity} \\
   & = & \sum_{\ovl{g} \in \hgp'\backslash \pgp/\hgp} \ip{\Ind{\hgp'}{g.\hgp \cap \hgp'}{\left(\Res{g.\hgp}{g.\hgp \cap\hgp'}{\;g.\chi}\right)}}{\chi'}{\hgp'}\ \mbox{by the Mackey Subgroup Theorem} \\
   & = & \sum_{\ovl{g} \in \hgp'\backslash \pgp/\hgp} \ip{\Res{g.\hgp}{g.\hgp \cap \hgp'}{g.\chi}}{\Res{\hgp'}{g.\hgp \cap \hgp'}{\chi'}}{g.\hgp \cap \hgp'}\ \mbox{by Frobenius reciprocity.}
 \end{eqnarray*}
This vanishes if and only if each of the summands vanishes, which happens if and only if $\Res{g.\hgp}{g.\hgp \cap \hgp'}{g.\chi} \neq \Res{\hgp'}{g.\hgp \cap \hgp'}{\chi'}$ for every $g$, since the characters concerned are linear.  The result follows.
\end{proof}

Write $\nilgp$ as a quotient $\tfnilgp/\kr$ of a finitely generated torsion-free nilpotent group $\tfnilgp$: for example, we may take $\tfnilgp$ to be the free class $c$ nilpotent group on $N$ generators for appropriate $N$ and $c$.  Let $\proj \colon \tfnilgp \to \nilgp$ be the canonical projection, and let $i \colon \kr \to \tfnilgp$ be inclusion.  Let $\tfpgp$, $\pkr$ be the pro-$p$ completions of $\tfnilgp$, $\kr$ respectively.  Then $\proj$ (respectively $i$) extends to a continuous homomorphism $\pproj \colon \tfpgp \to \pgp$ (respectively $\widehat{i}_p \colon \pkr \to \tfpgp$), and the three groups $\widehat{i}_p(\pkr)$, $\ker{\pproj}$, and the closure of $\kr$ in $\tfpgp$ all coincide (compare \cite[Chapter~1, Ex.~21]{DDMS}; because $\tfnilgp$ is finitely generated nilpotent, it can in fact be shown that $\widehat{i}_p$ is injective, and hence an isomorphism onto its image). Clearly $p$-admissible representations of $\pgp$ correspond bijectively to $p$-admissible representations of $\tfpgp$ that kill $\ker{\pproj}$.  Now $\kr$ is finitely generated (see, e.g., \cite[Lemma~1.2.2]{WagStGp}), so we can choose a Mal'cev basis $\theta_1,\ldots,\theta_s$ for $\kr$.  We identify the $\theta_i$ with their images in $\tfpgp$.

Let $\rts{p^n}$ be the group of all complex $p^n$th roots of unity, and let $\rts{p^\infty}$ be the group of all complex $p$-power roots of unity.
\begin{lemma}
 The groups $\rts{p^\infty}$ and $\Qq_p/\Zz_p$ are isomorphic.
\end{lemma}
\begin{proof}
 Let $p^{-\infty}\Zz \leq \Qq$ be the group of rational numbers of the form $np^{-r}$ for $n \in \Zz$ and $r$ a non-negative integer.  Then $p^{-\infty}\Zz \cap \Zz_p = \Zz$ and $\Zz_p p^{-\infty}\Zz = \Qq_p$, so $\Qq_p/\Zz_p \iso p^{-\infty}\Zz/\Zz$, by one of the standard group isomorphism theorems.  The map $q \mapsto e^{2\pi i q}$ gives an isomorphism from $p^{-\infty}\Zz/\Zz$ to $\rts{p^\infty}$.
\end{proof}
Let $\Phi \colon \rts{p^\infty} \to \Qq_p/\Zz_p$ be the isomorphism described above.  Any $p$-admissible linear character of a pro-$p$ group takes its values in $\rts{p^\infty}$, so we use $\Phi$ to identify $p$-admissible linear characters with $p$-admissible homomorphisms to $\Qq_p/\Zz_p$.  Under this identification, the product $\chi_1\chi_2$ of two $p$-admissible linear characters $\chi_1$ and $\chi_2$ (regarded as functions to ${\mathbb C}^*$) corresponds to their sum (regarded as functions to $\Qq_p/\Zz_p$).

Recall our notation of Eqn.\ \ref{eqn:sbgpcoords}. Let \(a_1,\ldots,a_R\) be a Mal'cev basis of \(\Delta\). Then any subgroup \(\hgp\leq\tfpgp\) has a good basis \(h_1,\ldots,h_R\) and we represent that basis by the tuple \(\lambda_{ij}\in\Zz_p\) such that \(h_i=a_1^{\lambda_{i1}}\cdots a_R^{\lambda_{iR}}\).

\begin{lemma}
\label{lem:par}
 Let ${\mathcal D}_p \subseteq \Zz_p^{R^2} \times \Qq_p^R$ be the set of tuples $(\lambda_{ij},y_k)$, where $1 \leq i,j \leq R$ and $1 \leq k \leq R$, satisfying the following conditions:\\
 (a) the $\lambda_{ij}$ form a good basis $h_1,\ldots,h_R$ for some finite-index subgroup $\hgp$ of $\tfpgp$ such that $\ker{\pproj}\leq \hgp$;\\
 (b) the prescription $h_i \mapsto y_i\ {\mathrm{mod}}\ \Zz_p$ gives a well-defined $p$-admissible homomorphism $\chi \colon \hgp \to \Qq_p/\Zz_p$ that kills $\ker{\pproj}$;\\
 (c) the induced character $\Ind{\tfpgp}{\hgp}{\chi}$ is irreducible.\\
 Then the sets $({\mathcal D}_p)_{p \textrm{ prime}}$ are uniformly $\emptyset$-definable in $p$.  Moreover, $\Ind{\tfpgp}{\hgp}{\chi}$ is a $p$-admissible character of $\tfpgp$ that kills $\ker{\pproj}$ and hence induces a $p$-admissible character of $\pgp$, and every $p$-admissible irreducible character of $\pgp$ arises in this way.
\end{lemma}
\begin{notation}
 Given $(\lambda_{ij},y_k) \in {\mathcal D}_p$, we write $\pair(\lambda_{ij},y_k)$ for the pair $(\hgp,\chi)$.  Since the $h_i$ generate $\hgp$ topologically, the $p$-admissible homomorphism $\chi$ defined by the $y_i$ is unique.
\end{notation}
\begin{proof}
 Condition (a) is uniformly $\emptyset$-definable in $p$, by Proposition~\ref{propn:gddef} (to the formulae that define the set of good bases we add the formulae $(\exists \nu_{1j},\ldots, \nu_{rj}\in \Zz_p)\ \theta_j= h_1^{\nu_{1j}}\cdots h_R^{\nu_{Rj}}$ for $1\leq j\leq s$).  Given that (a) holds, we claim that (b) holds if and only if there exists an $R^2$-tuple $(\mu_{ij})$ such that:\\
 (i) $(\mu_{ij})$ defines a good basis $k_1,\ldots,k_R$ for a finite-index subgroup $\kgp$ of $\tfpgp$;\\
 (ii) $\kgp \unlhd \hgp$;\\
 (iii) $\ker{\pproj} \subseteq \kgp$;\\
 (iv) there exist $y \in \Qq_p$, $r_1,\ldots,r_R \in \Zz_p$, $h \in \hgp$ such that the order of $y$ in $\Qq_p/\Zz_p$ is equal to $|\hgp/\kgp|$ and for every $i$ we have $\ovl{h^{r_i}} = \ovl{h_i}$ and $r_iy=y_i\ {\mathrm{mod}}\ \Zz_p$.  (Here $\ovl{x}$ denotes the image of $x \in \hgp$ under the canonical projection  $\hgp \to \hgp/\kgp$.)\medskip\\
 To see this, note that if (b) holds then $K:= \ker{\chi}$ is a finite-index subgroup of $\hgp$ which satisfies (ii) and (iii).  Take $(\mu_{ij})$ to be any tuple defining a good basis for $\kgp$.  Then $\hgp/\kgp$, being isomorphic to a finite subgroup of $\Qq_p/\Zz_p$, is cyclic, so choose $h \in \hgp$ that generates $\hgp/\kgp$ and choose $y \in \Qq_p$ such that $\chi(h)=y\ {\mathrm{mod}}\ \Zz_p$.  We can choose $r_1,\ldots,r_R \in \Zz$ such that $\ovl{h_i}=\ovl{h}^{r_i}$ for each $i$, and it is easily checked that (iv) holds.

 Conversely, suppose there exists a tuple $(\mu_{ij})$ satisfying (i)--(iv).  The map $\Zz_p \to \hgp$, $\lambda \mapsto h^\lambda$ is continuous because it is polynomial with respect to the Mal'cev coordinates, so there exists an open neighborhood $U$ of $0$ in $\Zz_p$ such that $h^\lambda \in \kgp$ for all $\lambda \in U$.  Since $\Zz$ is dense in $\Zz_p$, we may therefore find $n_1,\ldots,n_R \in \Zz$ such that $\ovl{h_i}=\ovl{h^{n_i}}$ for each $i$.  Hence $H/K$ is cyclic with generator $\ovl{h}$.

 We have a monomorphism $\beta \colon \hgp/\kgp \to \Qq_p/\Zz_p$ given by $\beta(\ovl{h}^n)=ny\ {\mathrm{mod}}\  \Zz_p$.  Let $\chi$ be the composition $\hgp \to \hgp/\kgp \stackrel{\beta}{\to} \Qq_p/\Zz_p$.  The canonical projection $\hgp \to \hgp/\kgp$ is continuous \cite[1.2~Proposition]{DDMS}, so we have $\chi(h^\lambda) = \lambda y\ {\mathrm{mod}}\ \Zz_p$ for every $\lambda \in \Zz_p$.  Condition (iv) implies that $\chi(h_i)=y_i\ {\mathrm{mod}}\ \Zz_p$ for every $i$, as required.  This proves the claim.

 Now condition (i) is uniformly $\emptyset$-definable in $p$, by Proposition~\ref{propn:gddef}.  Condition (iii) can be expressed as
 \begin{equation}
 \label{eqn:killker}
  (\forall \nu_1,\ldots,\nu_s \in \Zz_p) (\exists \sigma_1,\ldots,\sigma_R \in \Zz_p) \, \theta_1^{\nu_1} \cdots \theta_s^{\nu_s}=k_1^{\sigma_1} \cdots k_R^{\sigma_R}.
 \end{equation}
 Eqn.\ (\ref{eqn:killker}) can be expressed in terms of polynomials independent of $p$ over $\Qq$ in the $\mu_{ij}$, the $\nu_k$ and the $\sigma_l$, so condition (iii) is uniformly $\emptyset$-definable in $p$.  (Note that the $\theta_k$ are fixed elements of $\tfnilgp$, so their Mal'cev coordinates are not just elements of $\Zz_p$ but elements of $\Zz$.)

 Similar arguments show that conditions (ii) and (iv) are also uniformly $\emptyset$-definable in $p$.  In (iv), note that the conditions $\ovl{h^{r_i}} = \ovl{h_i}$ imply by the argument above that $\ovl{h}$ is a generator for $\hgp/\kgp$, so the condition that the order of $y$ in $\Qq_p/\Zz_p$ is equal to $|\hgp/\kgp|$ can be expressed as $\left( h^{y^{-1}} \in \kgp \right) \wedge \left( (\forall z\in \Qq_p)\ |z|< |y|\implies h^{z^{-1}} \not \in \kgp \right)$.  This shows that (condition (a))$\wedge$(condition (b)) is uniformly $\emptyset$-definable in $p$.

 Condition (iii) implies that $\chi$ kills $\ker{\pproj}$.  Hence $\Ind{\tfpgp}{\hgp}{\chi}$ kills $\ker{\pproj}$, so $\Ind{\tfpgp}{\hgp}{\chi}$ gives rise to an irreducible $p$-admissible character of $\pgp$.  By Lemma~\ref{lem:induction} (b), irreducibility of the induced character can be written as
  \[ (\forall g \in \tfpgp \sminus \hgp)(\exists h \in \hgp)\;ghg^{-1} \in \hgp\ \mathrm{and}\ \chi(ghg^{-1}) \neq \chi(h). \]
 Writing this in terms of the Mal'cev coordinates, we see that condition (c) is uniformly $\emptyset$-definable in $p$.

 By Lemma~\ref{lem:induction} (a), any $p$-admissible irreducible character $\sigma$ of $\pgp$ is of the form $\Ind{\pgp}{\lgp}{\chi}$ for some finite-index subgroup $\lgp$ of $\pgp$ and some $p$-admissible linear character $\chi$ of $\lgp$.  Let $\hgp$ be the pre-image of $\lgp$ under the canonical projection $\tfpgp \to \pgp$.  Regarding $\sigma,\chi$ as representations of $\tfpgp,\hgp$ respectively, it is easily checked that $\sigma=\Ind{\tfpgp}{\hgp}{\chi}$.  Choose $(\lambda_{ij})$ defining a good basis $h_1,\ldots,h_R$ for $\hgp$, and choose $y_k$ such that $\chi(h_k)=y_k\ \mathrm{mod}\ \Zz_p$ for all $k$.  The above argument shows that $(\lambda_{ij},y_k) \in {\mathcal D}_p$.  This completes the proof.
\end{proof}

Define $f_p \colon {\mathcal D}_p \to \Zz_p$ by $f_p(\lambda_{ij},y_k)=\lambda_{11}\cdots \lambda_{RR}$.  Define an equivalence relation ${\mathcal E}_p$ on ${\mathcal D}_p$ by $(\lambda_{ij},y_k) \sim (\lambda_{ij}',y_k')$ if $\Ind{\tfpgp}{\hgp}{\chi}$ and $\Ind{\tfpgp}{\hgp'}{\chi'}$ are twist-equivalent, where $(\hgp,\chi)=\pair(\lambda_{ij},y_k)$ and $(\hgp',\chi')=\pair(\lambda_{ij}',y_k')$.  The degree of $\Ind{\tfpgp}{\hgp}{\chi}$ equals $|f_p(\lambda_{ij},y_k)|^{-1}$ by Eqn.\ (\ref{eqn:indx}), and likewise for $(\lambda_{ij}',y_k')$, so if $(\lambda_{ij},y_k) \sim (\lambda_{ij}',y_k')$ then $|f_p(\lambda_{ij},y_k)|=|f_p(\lambda_{ij}',y_k')|$.

Construct $D_p$ and $E_p$ from ${\mathcal E}_p$, ${\mathcal D}_p$ and $f_p$ as in the paragraph following Remark~\ref{rem:growth}.  It follows from Lemma~\ref{lem:par} and the definition of ${\mathcal E}_p$ that $D_{p,n}$ is the union of precisely $a_{p^n}$ $E_{p,n}$-equivalence classes (note that if one representation of $\pgp$ is the twist of another by some linear character $\psi$ of $\tfpgp$ then $\psi$ is automatically a character of $\pgp$, by Observation~\ref{obs:twstquot}).
To complete the proof of Theorem~\ref{thm:rattwst_intro}, it suffices by Theorem~\ref{thm:rat} and Remark~\ref{rem:allprimes}.1 to show that $(D_p)_{p \textrm{ prime}}$ and $(E_p)_{p \textrm{ prime}}$ are uniformly $\emptyset$-definable in $p$.  But the sets $(D_p)_{p \textrm{ prime}}$ are uniformly $\emptyset$-definable in $p$ by Lemma~\ref{lem:par}, so it is enough to prove the following result.
\begin{proposition}
 The equivalence relations $({\mathcal E}_p)_{p \textrm{ prime}}$ are uniformly $\emptyset$-definable in $p$.
\end{proposition}
\begin{proof}
 Let $a_1,\ldots, a_R$ be a Mal'cev basis for $\tfnilgp$.  Let ${\mathcal D}_p' \subseteq \Qq_p^R$ be the set of $R$-tuples $(z_1,\ldots,z_R)$ such that the prescription $a_i \mapsto z_i\ \mathrm{mod}\ \Zz_p$ gives a well-defined $p$-admissible linear character of $\tfpgp$ that kills $\ker{\pproj}$.  We denote this character by $\globlin_z$.  Similar arguments to those in the proof of Lemma~\ref{lem:par} show that the sets $({\mathcal D}_p')_{p \textrm{ prime}}$ are uniformly $\emptyset$-definable in $p$.  Let $(z_1,\ldots,z_R) \in {\mathcal D}_p'$, let $(\hgp,\chi)=\pair(\lambda_{ij},y_k)$ and let $h_1,\ldots, h_R$ be the corresponding good basis for $H$.  Then $h_k= a_1^{\lambda_{k1}}\cdots a_R^{\lambda_{kR}}$, so $\globlin_z(h_k)= \lambda_{k1}z_1+\cdots + \lambda_{kR}z_R\ \mod \ \Zz_p$.
Hence $(H, (\Res{\tfpgp}{\hgp}{\globlin_z})\chi) = \pair(\lambda_{ij},y_k+\lambda_{k1}z_1+\cdots + \lambda_{kR}z_R)$.  Applying Lemma~\ref{lem:induction} (c), we see that if $(\hgp',\chi')=\pair(\lambda'_{ij},y'_k)$ then $(\lambda_{ij},y_k) \sim (\lambda_{ij}',y_k')$ if and only if
 \begin{eqnarray*}
  (\exists (z_1,\ldots,z_R) \in {\mathcal D}'_p)\; (\exists g \in \tfpgp)\; (\forall h \in \hgp) & & \\
  ghg^{-1} \in \hgp'\Rightarrow \left( (\Res{\tfpgp}{\hgp}{\globlin_z})\chi\right)(h)= \chi'(ghg^{-1}). & &
 \end{eqnarray*}
Writing this in terms of the Mal'cev coordinates, we obtain an equation independent of $p$ involving ${\mathcal D}_p'$ and absolute values of polynomials over $\Qq$ in the $\lambda_{ij}$, the $y_k$, the $\lambda_{ij}'$, the $y_k'$, the $z_k$, and the Mal'cev coordinates of $g$ and $h$.  We deduce that the equivalence relations $({\mathcal E}_p)_{p \textrm{ prime}}$ are uniformly $\emptyset$-definable in $p$, as required.
\end{proof}

\begin{remark}
\label{rem:multivariate}
 Using the multivariate version of Theorem~\ref{thm:rat}, one can obtain variations on Theorem~\ref{thm:rattwst_intro}: for instance, one can prove uniform rationality for the 2-variable zeta function that counts twist isoclasses of $p^n$-dimensional irreducible representations of $\nilgp$ factoring through a finite quotient of $\nilgp$ of order $p^m$.  We leave the details to the reader.
\end{remark}

 Next we give a variation on Theorem~\ref{thm:rattwst_intro} for nilpotent pro-$p$ groups.  Let $M$ be a topologically finitely generated nilpotent pro-$p$ group for some prime $p$.  Note that since every finite-index subgroup of $M$ is open and has $p$-power index, a representation $\rho\colon M\ra {\rm GL}_n({\Cc})$ is $p$-admissible if and only if it is continuous (with respect to the discrete topology on ${\rm GL}_n({\Cc})$).  Set $a_n= |\pirr{p^n}{M}/({\rm twisting})|$ and set $\zeta_M(s)= \sum_{n= 0}^\infty a_np^{-ns}$.
 
\begin{proposition}
\label{prop:arbnilpt}
 Let $p$, $M$ and $\zeta_M(s)$ be as above.  Then $\zeta_M(s)$ is a rational function of $p^{-s}$ with coefficients in $\Qq$.
\end{proposition}

\begin{proof}
 Let $\Gamma$ be a finitely generated dense subgroup of $M$, and choose an epimorphism $\pi$ from a torsion-free finitely generated nilpotent group $\Delta$ onto $\Gamma$.  Then $\pi$ gives rise to a continuous epimorphism $\widehat{\pi}_p$ from the pro-$p$ completion $\widehat{\Delta}_p$ to $M$.  The kernel $K$ of $\widehat{\pi}_p$ is a closed subgroup of $\widehat{\Delta}_p$, so $K$ is also topologically finitely generated.  Let $\Theta= \langle \theta_1,\ldots, \theta_s \rangle$ be a finitely generated dense subgroup of $K$.  The result now follows from the proof of Theorem~\ref{thm:abscissa_rat_intro} given above (cf.\ the paragraph after the proof of Lemma~\ref{lem:induction}).
\end{proof}

We finish the section by applying our approach to recover some results of Stasinski and Voll on the representation zeta functions of nilpotent groups arising from smooth unipotent group schemes.  Their parametrisation of irreducible representations uses the Kirillov orbit method; it allows one to prove strong uniformity properties of the representation zeta function at the cost of having to discard finitely many primes.  We give a brief summary of the necessary background---see \cite{StaVol} for details.  Let $K$ be a number field with ring of integers ${\mathcal O}$ and let ${\mathbf G}= {\mathbf G}_\Lambda$ be the smooth unipotent group scheme over ${\mathcal O}$ corresponding to a nilpotent Lie lattice $\Lambda$ over ${\mathcal O}$, in the sense of \cite[2.1]{StaVol}.  If $R$ is a ring extension of ${\mathcal O}$ then we denote by ${\mathbf G}(R)$ the group of $R$-points of ${\mathbf G}$.  Note that ${\mathbf G}({\mathcal O})$ is a finitely generated torsion-free nilpotent group; moreover, for any finitely generated nilpotent group $\Gamma$, there exists a smooth unipotent group scheme ${\mathbf H}$ over $\Zz$ such that $\zeta_{\Gamma, p}(s)= \zeta_{{\mathbf H}(\Zz), p}(s)$ for all $p$ sufficiently large.

Let ${\mathfrak p}$ be a nonzero prime ideal of ${\mathcal O}$.  Let $K_{\mathfrak p}$ be the completion of $K$ at ${\mathfrak p}$ and let ${\mathcal O}_{\mathfrak p}$ be the valuation ring of $K_{\mathfrak p}$.  Let $\zeta_{{\mathbf G}({\mathcal O}_{\mathfrak p})}(s): = \sum_{i=0}^\infty \widetilde{a}_{p^i}({\mathbf G}({\mathcal O}_{\mathfrak p}))p^{-is}$ be the zeta function that counts the twist isoclasses of continuous irreducible complex representations of the pro-$p$ group ${\mathbf G}({\mathcal O}_{\mathfrak p})$, where $p$ is the characteristic of the residue field of $K_{\mathfrak p}$.  It follows from Eqn.\ (\ref{eqn:SV}) below that for $p$ sufficiently large, $\widetilde{a}_{p^i}({\mathbf G}({\mathcal O}_{\mathfrak p}))= 0$ unless $p^i$ is a power of $q$, where $q$ is the cardinality of the residue field of $K_{\mathfrak p}$.  There is a ``refined Euler product''

\begin{equation}
 \zeta_{{\mathbf G}({\mathcal O})}(s)= \prod_{\mathfrak p} \zeta_{{\mathbf G}({\mathcal O}_{\mathfrak p})}(s)
\end{equation}

\noindent and the $p$-local representation zeta function is given by the ``mini Euler product''

$$  \zeta_{{\mathbf G}({\mathcal O}),p}(s)= \prod_{{\mathfrak p}| p} \zeta_{{\mathbf G}({\mathcal O}_{\mathfrak p})}(s). $$

Let $L$ be a finite extension of $K$ and let ${\mathcal O}_L$ be the ring of integers of $L$.  Let ${\mathfrak p}$ be a nonzero prime ideal of ${\mathcal O}$ and let ${\mathfrak B}$ be a nonzero prime ideal of ${\mathcal O}_L$ that divides ${\mathfrak p}$.  Let ${\mathfrak o}= {\mathcal O}_{\mathfrak p}$ and let ${\mathfrak D}$ be the valuation ring of the completion $L_{\mathfrak B}$.  Let $p$ be the residue field characteristic of ${\mathfrak o}$ and let $q, q^f$ be the cardinality of the residue field of ${\mathfrak o}$, ${\mathfrak D}$, respectively.  Note that ${\mathbf G}({\mathfrak D})$ is a topologically finitely generated nilpotent pro-$p$ group.  We will show that $\zeta_{{\mathbf G}({\mathfrak D})}(s)$ comes (up to a change of variable) from counting the equivalence classes of a family of equivalence relations that are uniformly $K$-definable over $\mathfrak{F}$ for an appropriate choice of $\mathfrak{F}$.

Let $d$, $k$, $r$ be as defined on \cite[p516]{StaVol}.  Let ${\mathbf Y}= (Y_1,\ldots, Y_d)$ be a tuple\footnote{Here and below we are following the notation of \cite{StaVol} and using bold-face letters to denote tuples.} of indeterminates and define the ${\lfloor r/2 \rfloor}\times {\lfloor r/2 \rfloor}$ commutator matrix ${\mathcal R}({\mathbf Y})$ as in \cite[p516]{StaVol} by choosing a basis for the ${\mathfrak D}$-Lie algebra ${\mathfrak g}$ that is associated to ${\mathbf G}$.   Then for any ${\mathbf y}\in {\mathfrak D}^d$, ${\mathcal R}({\mathbf y})$ is a matrix with entries from ${\mathfrak D}$.  As in the proof of \cite[Thm.~A]{StaVol}, we may choose the data that define ${\mathcal R}({\mathbf Y})$ in a global way and ensure that the quantities $b_i$ that appear in \cite[Eqn.\ (2.6)]{StaVol} are all zero, at the cost of discarding finitely many rational primes.  In particular, for $p$ sufficiently large---say, for $p> p_0$---the linear forms that appear as entries of the matrix ${\mathcal R}({\mathbf Y})$ have coefficients from ${\mathcal O}$, and these coefficients do not depend on $L$, ${\mathfrak B}$, ${\mathfrak p}$ or $p$.  We define the submatrix ${\mathcal S}({\mathbf Y})$ of ${\mathcal R}({\mathbf Y})$ as in \cite[p516]{StaVol}.

Let $\nu$, $\widetilde{\nu}$ be as defined on \cite[p518]{StaVol}.  Let $D$ be the set of tuples $({\mathbf y}, N, {\mathbf a}, {\mathbf c})\in {\mathfrak D}^d\times \Nn\times \Nn^{\lfloor r/2 \rfloor}\times \Nn^k$ such that ${\mathbf y}\neq {\mathbf 0}\ {\rm mod}\ {\mathfrak B}^N$, $\nu(\pi_N({\mathcal R}({\mathbf y})))= {\mathbf a}$ and $\widetilde{\nu}(\pi_N({\mathcal S}({\mathbf y})))= {\mathbf c}$, where $\pi_N$ denotes reduction of the matrix entries mod ${\mathfrak B}^N$.
Define $g(N, {\mathbf a})= \sum_{i= 1}^{\lfloor r/2 \rfloor} (N- a_i)$ and $h(N, {\mathbf c})= \sum_{i= 1}^k (N- c_i)$.  It follows from the definition of $\nu$ and $\widetilde{\nu}$ that if $({\mathbf y}, N, {\mathbf a}, {\mathbf c})\in D$ then $g(N, {\mathbf a}), h(N, {\mathbf c})$ are positive integers, and it is not hard to show using the theory of elementary divisors that $h(N, {\mathbf c})\leq 2g(N, {\mathbf a})$ (recall that ${\mathcal S}({\mathbf y})$ is a submatrix of ${\mathcal R}({\mathbf y})$).  Now define an equivalence relation $E$ on $D$ by
$$ ({\mathbf y}, N, {\mathbf a}, {\mathbf c})E({\mathbf y}', N', {\mathbf a}', {\mathbf c}') \iff $$
$$ N= N', \quad {\mathbf a}= {\mathbf a}', \quad {\mathbf c}= {\mathbf c}'\ \mbox{and } {\mathbf y}= {\mathbf y}'\ {\rm mod}\ {\mathfrak B}^{N+ 2g(N, {\mathbf a})- h(N, {\mathbf c})}. $$
It is easily seen that the functions $\nu$, $\widetilde{\nu}$, $g$ and $h$ are definable over $K$, so $D$ and $E$ are definable over $K$.  Set $D_l= \{({\mathbf y}, N, {\mathbf a}, {\mathbf c})\in D\mid g(N, {\mathbf a})= l\}$ for $l\in \Nn$ and let $E_l$ be the restriction of $E$ to $D_l$.  Let $e_{{\mathfrak D},l}$ be the number of equivalence classes of $E_l$ on $D_l$.

The point of the constructions above is to allow one to count certain coadjoint orbits in the dual of the Lie algebra ${\mathfrak g}$; this yields information about irreducible representations of ${\mathbf G}({\mathfrak D})$ via the Kirillov orbit method (see \cite{StaVol} for further details).  Stasinski and Voll show \cite[Eqn.\ (2.7)]{StaVol} that for $p$ sufficiently large---say, for $p> p_0$---we have

\begin{equation}
\label{eqn:SV}
 \zeta_{{\mathbf G}({\mathfrak D})}(s- 2)= \sum_{l\in \Nn} e_{{\mathfrak D},l}q^{-fls}.
\end{equation}

\noindent (Note that if $p> p_0$ then $\displaystyle e_{{\mathfrak D},l}= \sum_{\{(N, {\mathbf a}, {\mathbf c})\mid g(N,{\mathbf a})= l\}} q^{f(2g(N, {\mathbf a})- h(N, {\mathbf c}))} {\mathcal N}^0_{N, {\mathbf a}, {\mathbf c}}$, where ${\mathcal N}^0_{N,{\mathbf a}, {\mathbf c}}$ is as in \cite[Eqn.\ (2.6)]{StaVol}.  Moreover, although \cite[Eqn.\ (2.7)]{StaVol} is stated only for $L= K$, the equation holds for arbitrary $L$ because the coefficients of the linear forms that appear as entries of ${\mathcal R}({\mathbf Y})$ and ${\mathcal S}({\mathbf Y})$ do not change when one extends the field from $K$ to $L$.)

Now define $\mathfrak{F}_p$ to be empty if $p\leq p_0$ and the set of localisations $L_{\mathfrak B}$ if $p> p_0$, where $L$ runs over all the finite extensions of $K$ and ${\mathfrak B}$ runs over all the nonzero prime ideals of $L$ that divide $p$.  Set $\mathfrak{F}=\bigcup_p \mathfrak{F}_p$.  Then we see that $(D_{L_{\mathfrak B}})_{L_{\mathfrak B}\in {\mathcal L}}$ and $(E_{L_{\mathfrak B}})_{L_{\mathfrak B}\in {\mathcal L}}$ are uniformly $K$-definable in $\mathfrak{F}$; again, the key point is that the entries of ${\mathcal R}({\mathbf Y})$ and ${\mathcal S}({\mathbf Y})$ are linear forms with coefficients from ${\mathcal O}$, and these coefficients do not depend on $L$, ${\mathfrak B}$, ${\mathfrak p}$ or $p$.  Applying Theorem~\ref{thm:rat}, Remark~\ref{rem:allprimes}.3 and Proposition~\ref{prop:arbnilpt}, we obtain the following result.

\begin{theorem}
\label{thm:double}
 Let the notation be as above.  Let $S_{L_{\mathfrak B}}(t)=  \sum_{l\in \Nn} e_{{\mathfrak D},l}t^l$.  Then the power series $(S_{L_{\mathfrak B}}(t))_{L_{\mathfrak B}}$ are uniformly rational for $p\gg 0$.  In particular, $\zeta_{{\mathbf G}({\mathfrak D})}(s)$ has the form

 \begin{equation}
 \label{eqn:double}
  \zeta_{{\mathbf G}({\mathfrak D})}(s)= \frac{\sum_{l=1}^d q_l|X_l(\Ff_{q^f})|q^{-fls}}{q^{fn} \prod_{j=1}^{k} (1- q^{f(a_j- b_js)})}
 \end{equation}

\noindent  for all $p\gg 0$, where the $X_l$, etc., are as in Theorem~\ref{thm:rat}.

Moreover, each $\zeta_{{\mathbf G}({\mathfrak D})}(s)$ is a rational function of $p^{-s}$ with coefficients in $\Qq$.
\end{theorem}

\begin{remark}
\begin{enumerate}
 \item It is not stated explicitly in \cite{StaVol} that the power series $(S_{L_{\mathfrak B}}(t))_{L_{\mathfrak B}}$ are uniformly rational for $p\gg 0$, but this can be seen from the proof of \cite[Thm.~A]{StaVol}; cf.\ \cite[Sec.\ 4]{AKOVDuke}.  The final assertion of Theorem~\ref{thm:double}, however, is new: to prove rationality of $\zeta_{{\mathbf G}({\mathfrak D})}(s)$ for {\em every} ${\mathfrak D}$, we need Proposition~\ref{prop:arbnilpt} (cf.\ the discussion following Theorem~\ref{thm:rattwst_intro} in Section~\ref{sec:intro}).  Note that to apply the Kirillov orbit method, one needs to discard finitely many primes, so Eqn.\ (\ref{eqn:SV}) holds only when $p$ is sufficiently large.
 \item Theorem~\ref{thm:double} shows that $\zeta_{{\mathbf G}({\mathfrak D})}(s)$ depends on ${\mathfrak p}$ and ${\mathfrak D}$ only by way of the residue field of ${\mathfrak D}$.  A different expression for $\zeta_{{\mathbf G}({\mathfrak D})}(s)$ is given in \cite[Thm.~A]{StaVol}; this expression implies very strong uniformity behavior when one varies $L$ and ${\mathfrak B}$ for fixed ${\mathfrak p}$. 
 \item Stasinski and Voll show that $\zeta_{{\mathbf G}({\mathfrak D})}(s)$ satisfies a functional equation \cite[Thm.~A]{StaVol}.  Our methods---which apply to a very general class of problems---do not produce the functional equation that holds in this particular setting; there is no reason to expect the zeta function of an arbitrary definable equivalence relation to satisfy a functional equation.
 \item Dung and Voll show that the abscissa of convergence $\alpha$ of $\zeta_{{\mathbf G}({\mathfrak O})}(s)$ is rational and does not depend on ${\mathcal O}$, and they prove that $\zeta_{{\mathbf G}({\mathfrak O})}(s)$ can be meromorphically continued a short distance to the left of the line ${\rm Re}(s)= \alpha$ \cite[Thm.~A]{DunVol}.  For related results in the context of semisimple arithmetic groups, see \cite{Avni}, \cite{AizAvn1}, \cite{AizAvn2} and \cite{AKOVbasechange}.
\end{enumerate}
\end{remark}

\begin{example}
\label{ex:Heisenberg}
 Let ${\mathbf H}$ be the smooth unipotent group scheme over $\Zz$ corresponding to the Heisenberg group: so for every ring $R$, ${\mathbf H}(R)$ is the group of $3\times 3$ upper unitriangular matrices with entries from $R$.  Then for any number field $L$,

\begin{equation}
\label{eqn:Heisenberg}
 \zeta_{{\mathbf H}({\mathcal O}_L)}(s)= \frac{\zeta_L(s-1)}{\zeta_L(s)},
\end{equation}

\noindent where $\zeta_L(s)$ is the Dedekind zeta function of $L$.  For $L= \Qq$ this follows from results of Nunley and Magid \cite{NunMag}, who explicitly calculated the twist isoclasses of ${\mathbf H}(\Zz)$.  For $L$ a quadratic extension of $\Qq$, Eqn.\ (\ref{eqn:Heisenberg}) follows from work of Ezzat \cite[Theorem~1.1]{Ezz}, while for general $L$, it is a special case of results of Stasinski and Voll \cite[Thm.~B]{StaVol}.

The expression for the subgroup zeta function of ${\mathbf H}(\Zz)$ is more complicated: it is given by

$$ \xi_{{\mathbf H}(\Zz)}(s)= \frac{\zeta(s) \zeta(s-1) \zeta(2s-2) \zeta(2s-3)}{\zeta(3s-3)}, $$

\noindent where $\zeta(s)$ is the Riemann zeta function \cite[Section~1]{duSGruAnn}.  Expressions for $\xi_{{\mathbf H}({\mathcal O}_L)}(s)$ were obtained by Grunewald, Segal and Smith for quadratic and cubic extensions of $\Qq$, but no general formula is known (see \cite[Sec.~8]{GSS})\footnote[1]{Schein and Voll have obtained results on the structure of the normal subgroup zeta function of ${\mathbf H}({\mathcal O})$ \cite{SchVol1}, \cite{SchVol2}.}.  This suggests that the representation zeta function is better behaved than the subgroup zeta function.  The same seems to be true also for semisimple arithmetic groups \cite{LubMar}.
\end{example}

Theorem~\ref{thm:double} (and Example~\ref{ex:Heisenberg}) illustrate a significant difference between the subgroup zeta functions and representation zeta functions of groups of points of smooth unipotent group schemes: the former do not have the same double uniformity properties as the latter.  For instance, let $K= \Qq$ and let ${\mathbf G}$ be the smooth unipotent $\Zz$-scheme ${\mathbb G}_{\rm a}$ (the additive group).  The $p$-local subgroup zeta function of ${\mathbf G}(\Zz)= \Zz$ is given by $\displaystyle \xi_{\Zz,p}(s)= \frac{1}{1- p^{-s}}$.  Now let $L= \Qq(i)$ and let $p$ be any prime such that $p\equiv 3\ {\rm mod}\ 4$.  Let ${\mathfrak B}$ be a prime ideal of $\Zz[i]$ that divides $p$, and let ${\mathfrak D}$ be the valuation ring of the completion of $\Qq[i]$ at ${\mathfrak B}$; note that ${\mathfrak D}$ is isomorphic as an additive group to $\Zz_p^2$.  The residue field of ${\mathfrak D}$ has cardinality $q= p^2$.  Recall from Section~\ref{sec:zetagp} that if $\Gamma$ is a torsion-free finitely generated nilpotent group then subgroups of $\widehat{\Gamma}_p$ are parametrised by good bases, which for $\Gamma= \Zz$ and $\widehat{\Gamma}_p= \Zz_p$ are just 1-tuples of nonzero elements of $\Zz_p$.  But 1-tuples of nonzero elements of $\Zz[i]$ parametrise not finite-index subgroups of $\Zz[i]$ but finite-index subrings of $\Zz[i]$ (cf.\ \cite[Sec.~3]{GSS}), and $\xi_{{\mathfrak D},p}(s)$ is equal not to $\displaystyle \frac{1}{1- q^{-s}}= \frac{1}{1- p^{-2s}}$ but to $\displaystyle \frac{1}{(1- p^{-s})(1- p^{-s+1})}$ (this formula follows from \cite[Thm.~15.1]{LubSeg}).  In the language of Section~\ref{sec:rat}, the definable sets and equivalence relations that we use to parametrise finite-index subgroups via good bases are uniformly $\emptyset$-definable in $p$, but need not be uniformly $\emptyset$-definable in $\mathfrak{F}$ if we take $\mathfrak{F}_p$ to contain more than one extension of $\Qq_p$.

The uniform definability established in Theorem~\ref{thm:double} cannot be seen from our parametrisation of twist isoclasses, which involves good bases: to prove double uniformity one needs the Kirillov orbit formalism of \cite{StaVol}, as sketched above.  Our results give a further illustration of the power of the machinery developed in \cite{Jai}, \cite{VolAnn} and \cite{StaVol}.

\appendix

\section{Rationality results for \texorpdfstring{$p$}{p}-adic subanalytic equivalence relations, by Raf Cluckers}

\hfill\emph{\small Dedicated to Jan Denef, Lou van den Dries, Leonard Lipshitz and Angus Macintyre}

\subsection{Introduction}

One way to understand Denef's rationality results of \cite{Denef} for the generating power series $\sum_{n\geq 0} A_n T^n$ with coefficients 
$$
A_n:= \# \{\underline x \in (\Zz/(p^n\Zz))^d\mid \varphi(x,n),\ \underline x = x\bmod (p^n)\},
$$
for $n\geq 0$, where $\varphi$ is a definable condition on $n$ and on $x\in \Zz_p^d$, goes by writing $A_n$ as an integral
$$
\int_{\Zz_p^d} p^{-f(x,n)}|dx|
$$
for some well-chosen definable function $f$ and by studying the way such integrals may in general depend on the parameter $n$. This has started a vast study of the dependence of such integrals on more general parameters and on $p$, culminating in a way in the theory of motivic integration, see \cite{MacUnif}, \cite{PasEQ}, \cite{DL}, \cite{CLbounded}, \cite{HruKaz}. Most of this study works equally well in the semi-algebraic setting of the main body of the paper as in subanalytic settings, using model-theoretic results from the foundational \cite{Mac}, resp.~\cite{DvdD}.

In this appendix we show that this method also applies to generating power series $\sum_{n\geq 0} a_n T^n$ with coefficients
$$
a_n := \# (X_n/\!\!\sim_n ) ,
$$
where $\sim_n$ is a definable family of equivalence relations with finitely many equivalence classes, depending definably on an integer parameter $n\geq 0$. This is an alternative approach to the rationality result for $S_{L_{p}}(t)$ for each $p$ of Theorem \ref{thm:rat} in Section~\ref{sec:rat} in the case that one uses the semi-algebraic language (also called Macintyre's language) from \cite{Mac}, but the results and method of this appendix differ in two important ways from the main body of the paper. Firstly, our method is very robust in the choice of the language to define the equivalence relations. In particular, the subanalytic language of \cite{DvdD} can be used, or any intermediate structure between the semi-algebraic and this subanalytic language which is given by an analytic structure in the sense of \cite{CLip}. Secondly,
our method derives the rationality result, and more generally explains parameter dependence on arbitrary parameters, without using any form of elimination of imaginaries. Proving elimination of imaginaries is often not easy and seems to be dependent on the language in subtle ways as is shown in \cite{HHMAn}; in particular, in the subanalytic language on $\Qq_p$ the elimination of imaginaries is not yet completely understood. For simplicity of notation we will focus on those settings where elimination of imaginaries is not yet understood: the subanalytic setting on $p$-adic numbers and certain substructures coming from an analytic structure as in \cite{CLip} (which in fact includes the semi-algebraic case).  Our results also hold for many possible other languages allowing a typical kind of cell decomposition for the definable sets, but we leave this generality for the reader to work out. Our method can be adapted to obtain uniformity properties in $p$, both in the semi-algebraic and the subanalytic settings, but we leave this to future work (see the note added in proof below).

Although our arguments go through for any finite field extension of $\Qq_p$, we will work for simplicity with $\Qq_p$ itself.

\subsubsection{} Let us enrich the ring language on $\Qq_p$ with an analytic structure as in Section 4 of \cite{CLip}. As an example of an analytic structure, one may work with the subanalytic language as in \cite{DvdD}, {vdDHM}, where one adds to the ring language a function symbol $f$ for each power series $\sum_{i\in\Nn^n} a_i x^i$ in $n$ variables over $\Zz_p$ for any $n\geq 0$ whose coefficients go to zero as $|i|$ grows, and interpret it by evaluation, as the restricted analytic function
$$
\Qq_p^n\to \Qq_p: x\mapsto \left\{\begin{array}{l}
\sum_{i\in\Nn^n}a_i x^i\, \mbox{ if }\, x\in\Zz_p^n,\\
0\, \mbox{ otherwise.}
\end{array}\right.
$$
Let us further enrich this language by adjoining a sort for the value group $\Zz$, enriched with $\infty$ for the valuation of zero, the valuation map $\ord:\Qq_p\to \Zz\cup\{\infty\}$, a sort for the residue field $\Ff_p$, and a map $\ac:\Qq_p\to\Ff_p$ which sends $0$ to $0$ and nonzero $x$ to $x p^{-\ord x} \bmod (p)$. We denote this three-sorted language by $\cL_{an}$, where the notation refers to the analytic nature of the language.

The first theorem that we present in this introduction is a rather concrete form of Theorem~\ref{intgen} below. 

\begin{theorem} \label{int}
Let $\sim_y$ be an $\cL_{an}$-definable family of equivalence relations on nonempty sets $X_y\subseteq \Qq_p^d$ for some $d>0$, where the family parameters $y$ run over some $\cL_{an}$-definable set $Y$.
Suppose for each $y\in Y$ that each equivalence class of $\sim_y$ has nonempty interior in $\Qq_p^d$.
Then there exist $N>0$ and $\cL_{an}$-definable families of functions $f_y:X_y\to \Zz\cup\{\infty\}$ and $\alpha_y:X_y\to \{1,2,\ldots,N\}$,
such that for each $y\in Y$ and each $a\in X_y$,
$$
\int_{x\sim_y a} \frac{p^{- f_y(x)}}{\alpha_y(x)} |dx|    = 1 ,
$$
where $|dx|$ stands for the Haar measure on $\Qq_p^d$ normalized so that $\Zz_p^d$ has measure $1$, and where $p^{-\infty}$ stands for $0$.
\end{theorem}

By the theorem and with its notation, if moreover each quotient $X_y/\!\!\sim_y$
is finite, say, of size $a_y$, it immediately follows for $y\in Y$ that
\begin{equation}\label{des}
\int_{x\in X_y} \frac{p^{f_y(x)}}{\alpha_y(x)} |dx|    = a_y,
\end{equation}
which follows the philosophy mentioned above of relating finite counting to taking integrals (this philosophy is also followed in Section~\ref{sec:rat} in the semi-algebraic context, via elimination of imaginaries).
The integral description Eqn.\ (\ref{des}) and the more flexible variant Theorem~\ref{intgen} of Theorem~\ref{int} lead in a nowadays standard way to the following rationality result.  Note that a multivariate version of Theorem~\ref{count} (namely  replacing the single variable
$t$ with a tuple, as in Theorem \ref{thm:rat_intro}), as well as other variants, can be obtained by similar arguments.

\begin{theorem}\label{count}
Let $\sim_n$ be an $\cL_{an}$-definable family of equivalence relations on nonempty sets $X_n\subseteq \Qq_p^d$ for some $d>0$, where  $n$ runs over non-negative integers.
Suppose for each $n\geq 0$ that the quotient $X_n/\!\!\sim_n$ is finite, say, of size $a_n$. Then
$$
\sum_{n\geq 0} a_n t^n
$$
is a rational power series over $\Qq$ whose denominator is a product of factors of the form $(1-p^it^j)$ for some integers $i$ and some $j>0$.
\end{theorem}

\subsubsection{Sketch of differences with main body}\label{sec:sketch-app}
Before giving detailed proofs, let us give a sketch of the new ideas and the differences with the main body of the paper. Given a definable equivalence relation $\sim$ on a definable set $X$, in the main body of the paper one performs a definable transformation of the set $X$ to a simpler set $X'\subseteq \Zz_p^k$ for some $k$, with a corresponding equivalence relation $\sim'$ on $X'$, so that the equivalence class $x/\!\!\sim'$ of $x\in X'$ under $\sim'$ has a volume which is an integer power of $p$. Calling this integer exponent $f(x)$, the number of equivalence classes of $\sim$, if finite, equals the integral 
$$
\int_{x\in X'} p^{-f(x)}.
$$
This transformation from $X,\sim$ to $X',\sim'$ is achieved via elimination of imaginaries in the main body of the paper.
In this appendix, the simplification procedure is more elementary: instead of transforming $X$, we construct a definable subset $X''\subseteq X$, so that the intersection of $X''$ with $x/\!\!\sim$ for any $x\in X$ has positive volume $a(x) p^{f(x)}$, where $a(x)$ is an integer between $1$ and $N$ for some $N$, $f(x)$ is an integer, and where $f(x)$ and $a(x)$ depend definably on $x\in X$. Fixing the value of $a(x)$ subsequently for the values $1,\ldots, N$, one gets that the number of equivalence classes of $\sim$, if finite, equals the sum
$$
\sum_{i=1}^N \frac{1}{i} \int_{x\in X'',\ a(x)=i} p^{-f(x)}.
$$

When working out parameter versions of these integrals, rationality follows via either approach.

Finding such a subset $X''$ of $X$ can be done rather elementarily, by decomposing each $x/\!\!\sim$ into cells on the one hand, and, by looking at maximal balls (multi-balls in the general, higher-dimensional case) included in $x/\!\!\sim$ on the other hand. Roughly, the union of all these maximal multi-balls will form $X''$. The factor $a(x)$ is uniformly bounded by the number of cells in a decomposition of the $x/\!\!\sim$ into cells, which is bounded uniformly in $x$ by the cell decomposition result.

\subsection{Proofs via subsets instead of via EI} 

As mentioned in Section \ref{sec:sketch-app}, the proof of rationality given in this appendix relies on choosing simple subsets instead of transforming using EI. To do this, let us recall some aspects of cell decomposition for definable sets.

For integers $m>0$ and $n>0$, write  $Q_{m,n}$ for the set of all $p$-adic numbers of the form $p^{na}(1+p^{m}x)$ with $x\in \Zz_p$ and $a\in\Zz$. 

The following lemma is a direct corollary of cell decomposition results in \cite{Ccell} and \cite[Section~6]{CLip}.
\begin{lem}\label{cd}
For any $\cL_{an}$-definable sets $Y$ and $X\subseteq Y\times \Qq_p$, one can write $X$ as a finite disjoint union of $\cL_{an}$-definable sets of the form
$$
\{(y,x)\in Y\times \Qq_p  \mid  \ord (x-c(y)) \in G_y,\ (x-c(y)) \in \lambda Q_{m,n}     \},
$$
where $c:Y\to \Qq_p$ is an $\cL_{an}$-definable function, $G_y$ is an $\cL_{an}$-definable family of subsets of $\Zz\cup\{\infty\}$ with parameter $y\in Y$, and $\lambda$ lies in $\Qq_p$.
\end{lem}

Note that any set $Q_{m,n}$ equals a finite disjoint union of sets of the form $\lambda P_\ell$ for $\lambda\in \Qq_p$, where $P_\ell$ stands for the nonzero $\ell$th powers in $\Qq_p$, and also the other way around: any set $P_\ell$ equals a finite disjoint union of sets of the form $\lambda Q_{m,n}$ for $\lambda\in \Qq_p$.

The rest of this note is devoted to the proofs of Theorems~\ref{int},
\ref{count} and \ref{intgen}. 
We first give some definitions and lemmas. By a ball we mean a subset $B\subseteq\Qq_p$ of the form
$$
\{x\in\Qq_p\mid \ord (x-c) > g\}
$$
for some $g\in\Zz$ and some $c\in\Qq_p$.

Let $\Vol$ stand for the Haar measure on $\Qq_p$, normalized so that $\Zz_p$ has measure $1$.

\begin{defn}[Multi-balls]\label{defn:multiB}
Let $n\geq 1$, $r_i\geq 0$ for $i=1,\ldots,n$, and let a nonempty set $Y\subseteq \Zz_p^n$ be given.

If $n=1$, then $Y$ is called a \emph{multi-ball of multi-volume $r_1$} if $r_1=\Vol(Y)$ and either $Y$ is a singleton (in which case $r_1=0$), or $Y$ is a ball (in which case $r_1>0$).

If $n\geq 2$, then the set $Y$ is called a \emph{multi-ball of multi-volume $(r_1,\ldots,r_n)$} if and only $Y$ is of the form
$$
\{(x_1,\ldots,x_n)\mid (x_1,\ldots,x_{n-1})\in A,\ x_n\in B_{x_1,\ldots,x_{n-1}}\},
$$
where $A\subseteq \Zz_p^{n-1}$ is a multi-ball of multi-volume $(r_1,\ldots,r_{n-1})$, $B_{x_1,\ldots,x_{n-1}}$ is a subset of $\Zz_p$ which may depend on $(x_1,\ldots,x_{n-1})$,
with $\Vol(B_{x_1,\ldots,x_{n-1}})=r_n$,
and such that $B_{x_1,\ldots,x_{n-1}}$ is either a singleton or a ball. The multi-volume of a multi-ball $Y$ is denoted by $\MVol(Y)$.
\end{defn}

An example of a multi-ball in $\Zz_p^3$ of multi-volume $(1,0,p^{-1})$ is the set
$$
\{(x,y,z)\mid  x\in \Zz_p,\ y=x,\ z\in x + p\Zz_p  \}.
$$

\begin{defn}
Let us put on $\Rr^n$ the reverse lexicographical ordering.
Consider a set $X\subseteq \Zz_p^n$. The multi-box of $X$, denoted by $\MB(X)$, is the union of the multi-balls $Y$ contained in $X$ and with maximal multi-volume $\MVol(Y)$ in $\Rr^n$ (for the reverse lexicographical ordering on $\Rr^n$), where maximality is among all multi-balls contained in $X$. We write $\MVol(X)$ for $\MVol(Y)$ for any multi-ball $Y$ contained in $X$ with maximal multi-volume.
\end{defn}

For a set $X\subseteq \Zz_p^n$, we next define, by induction on $n$, an $\Nn$-valued function $\MNum_X$ on $X$ called the multinumber function of $X$. 

\begin{defn}
For a set $X\subseteq \Zz_p$, let $\MNum_X$ be the constant function on $X$ taking as value the number of distinct multi-balls $Y$ contained in $X$ with maximal multi-volume if this is finite, and taking the value $+\infty$ otherwise.

For a set $X\subseteq \Zz_p^n$ with $n>1$, let $p:\Zz_p^n\to \Zz_p^{n-1}$ be the projection on the first $n-1$ coordinates. We define $\MNum_X:X\to \Nn$ as the function sending $x=(p(x),x_n)$ to the product
$$
\MNum_{p(X)}(p(x)) \cdot  \MNum_{X_{p(x)}}(x_n),
$$
where $X_{p(x)}\subseteq \Qq_p$ is the fiber above $p(x)$ under the projection map $X\to p(X)$. Here, the product of $+\infty$ with any $a>0$ is set to be $+\infty$.
\end{defn}

The following two simple lemmas are key.

\begin{lem}\label{MBcell}
Let $X$ be a nonempty 
subset of $\Zz_p$ satisfying $X=\MB(X)$ and let $N\geq 1$ be an integer. Suppose that $X$ can be written as the disjoint union of $N$
sets of the form
\begin{equation}\label{cel}
A_j = \{x\in \Zz_p  \mid  \ord (x-c_j) \in G_j,\ (x-c_j) \in \lambda_j Q_{m_j,n_j}     \},
\end{equation}
for $j=1,\ldots,N$, where $c_j$ and $\lambda_j$ lie in $\Qq_p$, $G_j$ is a
subset of $\Zz\cup\{\infty\}$, and $m_j,n_j\geq 1$.
Then one has for $x\in X$ that
$$
\MNum_X(x) \leq N.
$$
\end{lem}
\begin{proof}
If $X$ is a finite set, then the $A_j$ are of size at most $1$,
and then the bound is clear.
Hence, we may and do suppose that $X$ is infinite. Then at least one of the sets $A_j$ is infinite, and since any infinite set of the form Eqn.\ (\ref{cel}) contains
a ball, it follows that $X$ contains at least one ball of maximal size. Since $\Zz_p$ has finite measure and since $X=\MB(X)$, $X$ equals a finite union of balls of the same volume, and hence, $\MNum_X(x)$ is finite, nonzero, and moreover constant since $n=1$. Write $s$ for $\MNum_X(x)$. The set $X$ thus equals a disjoint union of balls $B_i$ for $i=1,\ldots,s$ all of equal maximal volume $V$ (where maximality is among the balls contained in $X$).  
By the simple form of Eqn.\ (\ref{cel}), each of the sets $A_j$ for $j=1,\ldots,N$ contains at most one ball of maximal volume among all the balls included in $A_j$ (obtained by replacing $G_j$ with its minimum). Write $B_{A_j}$ for this ball of maximal volume contained in $A_j$ if it exists, and otherwise let $B_{A_j}$ be the empty set. If the volume of $B_{A_j}$ equals $V$, then $B_{A_j}$ equals one of the $B_i$, and we can replace $X$ with $X\setminus B_i$ and $A_j$ by $A_j\setminus B_{A_j}$ and prove the lemma for this new situation (with $N$ replaced by $N-1$ if $A_j\setminus B_{A_j}$  is empty, and with $N$ unchanged if $A_j\setminus B_{A_j}$  is nonempty). Hence, it is enough to prove the lemma when for each $j=1,\ldots,N$ we have
\begin{equation}\label{1bis} 
\Vol(B_{A_j}) \leq V/p.
\end{equation}
Further, by the simple form of Eqn.\ (\ref{cel}), one has for each $j$ that
\begin{equation}\label{1}
\Vol(A_j) \leq \frac{p}{p-1} \cdot \Vol(B_{A_j}).
\end{equation}
Indeed, writing $g_j$ for the minimum of $G_j$, if $B_{A_j}$ is nonempty then $B_{A_j}$ equals
$$
\{x\in \Zz_p  \mid  \ord (x-c_j) = g_j,\ (x-c_j) \in \lambda_j Q_{m_j,n_j} \}
$$
and the set $A_j$ is clearly contained in
$$
\{x\in \Zz_p  \mid  \ord (x-c_j) \geq  g_j,\ (x-c_j) \in \lambda_j Q_{m_j,1}\},
$$ 
whose volume equals $\frac{p}{p-1} \cdot \Vol(B_{A_j})$.
We calculate, by finite additivity of $\Vol$,
\begin{equation}\label{2}
sV = s \Vol(B_1)  =  \sum_{i=1}^s\Vol(B_i) = \Vol(\bigcup_{i=1}^s B_i) = \Vol(X) = \sum_{j=1}^N \Vol(A_j).\\
\end{equation}
Combining Equations (\ref{1}), (\ref{1bis}) and (\ref{2}), the lemma follows.
\end{proof}

\begin{lem}\label{MB}
Let  $\sim_y$ be an $\cL_{an}$-definable family of equivalence relations on $\Zz_p^n$ for the family parameter $y$ running over some $\cL_{an}$-definable set $Y$.
For $x\in \Zz_p^n$, write $x/\!\!\sim_y$ to denote the equivalence class of $x$ modulo $\sim_y$. We regard $x/\!\!\sim_y$ as a subset of $\Zz_p^n$.
Then the following properties hold.\\
The union
$$
\bigcup_{x\in \Zz_p^n} \MB(x/\!\!\sim_y)
$$
is an $\cL_{an}$-definable family of subsets of $\Zz_p^n$ with parameter $y\in Y$. There exists
an $\cL_{an}$-definable family of functions $g_y:\Zz_p^n\to (\Zz\cup\{\infty\})^n$ such that
$(p^{- g_{y,i}(x)})_{i=1}^n$ equals $\MVol(x/\!\!\sim_y)$ for each $x$ in $\Zz_p^n$. Finally, $x\mapsto \MNum_{ \MB(x/\sim_y) }(x)$ has uniformly bounded range (uniformly bounded in $x\in \Zz_p^n$ and in $y\in Y$), and depends definably on $x$ and $y$.
\end{lem}

\begin{proof}
Clearly the condition on $x\in\Zz_p^n$ to lie inside $\MB( x/\!\!\sim_y )$ is an $\cL_{an}$-definable condition, and also the existence of the $\cL_{an}$-definable family of functions $g_y$ is immediate.

We now show the finiteness of $\MNum_{ \MB(x/\sim_y) }$ and that it is uniformly bounded in $x$ and $y$.
It is enough, by induction on $n$ and by the definition 
of $\MNum$ as a product, 
to consider the case that $n=1$.
Let us thus assume that $n=1$.
By Lemma~\ref{cd}, applied to the family of subsets
$$
\MB(x/\!\!\sim_y)\subseteq \Zz_p
$$
with family parameter $(y,x)$, there exists $N\geq 1$ such that any set $\MB(x/\!\!\sim_y)$ equals a finite disjoint union of at most $N$ 
definable sets of the form in Eqn.\ (\ref{cel}) of Lemma~\ref{MBcell}. Applying that lemma to our family yields that
$$
\MNum_{ \MB( x/\sim_y) }(x)\leq N,
$$
for all $x$ and $y$. This proves that $\MNum_{\MB( x/\sim_y ) }$ has a uniformly bounded range, uniformly in $x$ and $y$. Having such a uniformly bounded range, the definability of $\MNum_{\MB( x/\sim_y ) }$ on $x$ and $y$ becomes an exercise.
\end{proof}

Let $I$ be a subset of $\{1,2,\ldots, d\}$ for some $d\geq 1$.
Let $\mu_{I,d}$ be the measure on $\Qq_p^d$ which is the product measure of the following measures on the $d$ Cartesian factors of $\Qq_p^d$: the normalized Haar measure on the $i$th factor $\Qq_p$ of $\Qq_p^d$ for $i\in I$, and the counting measure on the $j$th factor $\Qq_p$ of $\Qq_p^d$ for $j\not\in I$. 

The following proposition is a close variant of the well-known rationality result from \cite{DvdD}.
\begin{proposition}
\label{Denef}
 Let $f_n:\Zz_p^d\to \Zz \cup\{\infty\}$ be an $\cL_{an}$-definable family of functions, with an integer parameter $n\geq 0$.
Suppose that, for each $n\geq 0$, the function $x\mapsto p^{-f_n(x)}$ is integrable for the measure $\mu_{I,d}$, with $I$ a subset of $\{1,2,\ldots, d\}$.
Then the generating power series
$$
\sum_{n\geq 0} X_n t^n
$$
with
$$
X_n = \int_{x\in \Zz_p^d} p^{-f_n(x)} \mu_{I,d}
$$
is a rational power series over $\Qq$, with denominator a product of factors of the form $(1-p^it^j)$ for some integers $i$ and some $j>0$.
\end{proposition}
\begin{proof}
By Lemma~\ref{cd}, by reordering the coordinates so that $I=\{1,\ldots,a\}$ for some $a\geq 0$, and by finite additivity of the integral operator, one reduces to the case that the set $\{x\mid p^{-f_n(x)}\not=0\}$ is contained in the graph of an $\cL_{an}$-definable function
$$
\Zz_p^a \to \Zz_p^b
$$
for $b$ with $a+b=d$. But then one may suppose that $I=\{1,\ldots,d\}$, by replacing $d$ with $a$. Now the result is a standard variant of the rationality result for $p$-adic integrals from \cite{DvdD} (where the slightly more general integrability condition has been brought into the picture more recently, see \cite[Section 3]{CGH}).
\end{proof} 
Proposition~\ref{Denef} has several generalizations. For example, parameter integrals of a more general type and with more general parameters for any of the sorts $\Zz,\Qq_p,\Ff_p$, as well as uniformity in $p$, are well understood, see, e.g., \cite{PasEQ}, \cite{CLbounded}. We will not need more general results of this type here, and can come directly to the main result.

\begin{theorem}\label{intgen}
Let $\sim_y$ be an $\cL_{an}$-definable family of equivalence relations on nonempty sets $X_y\subseteq \Qq_p^d$ for some $d>0$, where the family parameters $y$ run over some $\cL_{an}$-definable set $Y$.
Then there exist 
$N>0$ and $\cL_{an}$-definable families of functions $f_{I,y}:X_y\to \Zz\cup\{\infty\}$ and $\alpha_{y}:X_y\to \{1,\ldots,N\}$,
such that for each $y\in Y$ and each $a\in X_y$,
\begin{equation}\label{sumI}
\sum_I \int_{x\sim_y a} \frac{p^{-f_{I,y}(x)}}{\alpha_{y}(x)} \mu_{I,d}(x)    = 1 ,
\end{equation}
where the sum runs over the subsets $I$ of $\{1,\ldots, d\}$.
\end{theorem}

\begin{proof}
Clearly we may suppose that the sets $X_y$ are subsets of $\Zz_p^d$, by replacing $d$ with $2d$ and by applying coordinate-wise the map sending  $w\in \Qq_p$ to $(w,0)\in\Zz_p^2$ if $|w|\leq 1$ and to $(0,w^{-1})\in\Zz_p^2$ if $|w|> 1$ and by replacing the sets $X_y$ correspondingly. 
Apply Lemma~\ref{MB} to the family $\sim_y$ to find an $\cL_{an}$-definable family of functions $g_y=(g_{y,i})_{i=1}^d$. Now, given $I\subseteq \{1,2,\ldots, d\}$, one can take for $f_{I,y}$ the function that maps $x$ to the sum of the finite component functions
$$
\sum_{i,\ g_{y,i}(x)\not= \infty} g_{y,i}(x)
$$
if $x$ lies in $\MB(x/\!\!\sim_y)$ and $\MB(x/\!\!\sim_y)$ has nonzero and finite $\mu_{I,d}$-measure, and to $\infty$ in all other cases. For $\alpha_y(x)$ one takes $\MNum_{\MB( x/\sim_y ) }(x)$ if $x$ lies in $\MB(x/\!\!\sim_y)$, and zero if $x$ lies in $x/\!\!\sim_y$ but outside $\MB(x/\!\!\sim_y)$. The $\alpha_y$ are an $\cL_{an}$-definable family of functions with finite range by Lemma~\ref{MB}. Clearly Eqn.\ (\ref{sumI}) holds for all $y\in Y$ and $a\in X_y$, as desired.
\end{proof}

We can now prove the rationality result of Theorem~\ref{count}.

\begin{proof}[Theorem~\ref{count}]
Consider $N$, $f_{I,n}$ for each $I\subseteq \{1,\ldots,d\}$ and $\alpha_{n}$ as given by Theorem~\ref{intgen}, with $Y$ the set of nonnegative integers $n$. For each integer $i$ with $1\leq i\leq  N$, 
let $X_{n,i}$ be the subset of $X_n$ on which $\alpha_{n}$ takes the value $i$. Let $a_{n,i}$ be number of equivalence classes of the restriction of $\sim_n$ to $X_{n,i}$ if $X_{n,i}$ is nonempty, and let $a_{n,i}$ be zero otherwise.
Since clearly $a_n= \sum_{i=1}^N a_{n,i}$ for all $n\geq 0$, one has  
$$
\sum_{n\geq 0} a_n T^n=\sum_{i=1}^N\sum_{n\geq 0 } a_{n,i} T^n.
$$
Also, for each $n\geq 0$ and each $i \in \{1,2,\ldots,N\}$, 
\begin{equation}\label{eq:app:final}
i a_{n,i} = \sum_{I\subseteq \{1,\ldots,d\} }\int_{x\in X_{n,i}} p^{ - f_{I,n}(x) } \mu_{I,d}(x).
\end{equation}
Now we are done since for each $i\in \{1,2,\ldots,N\}$, the integer multiple $i\sum_{n\geq 0}a_{n,i}T^n$ of $\sum_{n\geq 0}a_{n,i}T^n$ is rational and of the desired form
by Eqn.\ (\ref{eq:app:final}) and Proposition \ref{Denef}. 
\end{proof}

\subsubsection*{Acknowledgments}

The authors wish to thank Thomas Rohwer, Deirdre Haskell, Dugald Macpherson and Elisabeth Bouscaren for their comments on earlier drafts of this work, Martin Hils for suggesting that the proof could be adapted to finite extensions and Zo\'e Chatzidakis for pointing out an error in how constants were handled in earlier versions. The second author is grateful to Jamshid Derakhshan, Marcus du Sautoy, Andrei Jaikin-Zapirain, Angus Macintyre, Dugald Macpherson, Mark Ryten, Alexander Stasinski, Christopher Voll and Michele Zordan for helpful conversations.
  We are grateful to Alex Lubotzky for suggesting studying representation growth; several of the ideas in Section~\ref{sec:twst} are due to him.

The first author was supported by the European Research Council under the European Union's Seventh Framework Programme (FP7/2007-2013) / ERC Grant agreement no.\ 291111/ MODAG, the second author was supported by a Golda Meir Postdoctoral Fellowship at the Hebrew University of Jerusalem and the third author was partly supported by ANR MODIG (ANR-09-BLAN-0047) Model Theory and Interactions with Geometry.

The author of the appendix would like to thank M. du Sautoy, C. Voll, and Kien Huu Nguyen for interesting discussions on this and related subjects.  He was partially supported by the European Research Council under the European Community's Seventh Framework Programme (FP7/2007-2013) with ERC Grant Agreement nr. 615722
MOTMELSUM and he thanks the Labex CEMPI  (ANR-11-LABX-0007-01).

We are grateful to the referee for their careful reading of the paper and for their many comments, corrections and suggestions for improving the exposition.

\subsubsection*{Added in proof}

After this paper was submitted, we learnt that the method of the appendix has been generalized by Kien Huu Nguyen \cite{KHNguyen} to the uniform $p$-adic and uniform $\Ff_q((t))$ cases with the subanalytic languages.


\providecommand{\bysame}{\leavevmode\hbox to3em{\hrulefill}\thinspace}
\providecommand{\MR}{\relax\ifhmode\unskip\space\fi MR }
\providecommand{\MRhref}[2]{%
  \href{http://www.ams.org/mathscinet-getitem?mr=#1}{#2}
}
\providecommand{\href}[2]{#2}

\bigskip
\noindent \textsc{Institute of Mathematics,
Hebrew University,
Jerusalem 91904,
Israel, and Mathematical Institute,
University of Oxford,
Andrew Wiles Building,
Oxford,
OX2 6GG,
United Kingdom} \\
{\em Email address:} {\tt ehud@math.huji.ac.il }

\bigskip
\noindent \textsc{{Department of Mathematics,
University of Aberdeen,
King's College,
Fraser Noble Building,
Aberdeen AB24 3UE,
United Kingdom}
} \\
{\em Email address:} {\tt b.martin@abdn.ac.uk}

\bigskip
\noindent \textsc{University of California, Berkeley,
Mathematics Department,
Evans Hall,
Berkeley, CA, 94720-3840, United States} \\
{\em Email address:} {\tt silvain.rideau@berkeley.edu}

\bigskip
\noindent \textsc{Universit\'e Lille 1, Laboratoire Painlev\'e, CNRS -
UMR 8524, Cit\'e Scientifique, 59655 Villeneuve d'Ascq Cedex, France, and
KU Leuven, Department of Mathematics, Celestijnenlaan 200B, B-3001
Leuven, Belgium} \\
{\em Email address:} {\tt Raf.Cluckers@math.univ-lille1.fr}

\end{document}